\documentclass[11pt]{article}
\usepackage[utf8]{inputenc}
\usepackage{blindtext}
\usepackage{soul}
\usepackage{latexsym} 
\usepackage{amsmath,amssymb,amsfonts,amsthm} 
\usepackage{graphics,psfrag} 
\usepackage{fancyhdr,lastpage} 
\usepackage{color,verbatim} 
\usepackage[hidelinks]{hyperref} 
\usepackage{hyperref} 
\usepackage[table]{xcolor} 
\usepackage{metalogo}
\usepackage{etoolbox}
\patchcmd{\thebibliography}{\section*{\refname}}{}{}{}

\newtheorem{theorem}{Theorem}[section]

\newtheorem{lemma}[theorem]{Lemma}
\newtheorem{assumption}[theorem]{Assumption}
\newtheorem{definition}[theorem]{Definition}

\hypersetup{
    colorlinks=false,
    pdfborder={0 0 0},
}

\usepackage{calc}

\usepackage[a4paper,top=1in,bottom=0.8in,right=0.8in,left=0.8in]{geometry}
\let\oldbibliography\thebibliography
\renewcommand{\thebibliography}[1]{%
  \oldbibliography{#1}%
  \setlength{\itemsep}{-1.5mm}%
}

\def\R{\mathbb{R}}

\newcommand{\be}{\begin{equation}}
\newcommand{\ee}{\end{equation}}
\newcommand{\bea}{\begin{eqnarray}}
\newcommand{\eea}{\end{eqnarray}}
\newcommand{\beann}{\begin{eqnarray*}}
\newcommand{\eeann}{\end{eqnarray*}}
\newcommand{\benn}{\begin{equation*}}
\newcommand{\eenn}{\end{equation*}}

\fancypagestyle{plain}{%
  \fancyhf{}%
  \fancyfoot[C]{\footnotesize Page \thepage\ of \pageref{LastPage}}%
}

\title{Stochastic Rotating Waves}
\author{Christian Kuehn and James MacLaurin and Giulio Zucal}
\date{\today}

\begin{document}
\maketitle

\begin{abstract}
Stochastic dynamics has emerged as one of the key themes ranging from models in applications to theoretical foundations in mathematics. One class of stochastic dynamics problems that has received considerable attention recently are travelling wave patterns occurring in stochastic partial differential equations (SPDEs), i.e., how deterministic travelling waves behave under stochastic perturbations. In this paper, we start the mathematical study of related class of problems: stochastic rotating waves generated by SPDEs. We combine deterministic dynamics PDE techniques with methods from stochastic analysis. We establish two different approaches, the variational phase and the approximated variational phase, for defining stochastic phase variables along the rotating wave, which track the effect of noise on neutral spectral modes associated to the special Euclidean symmetry group of rotating waves. Furthermore, we prove transverse stability results for rotating waves showing that over certain time scales and for small noise, the stochastic rotating wave stays close to its deterministic counterpart.         
\end{abstract}

\section{Introduction}
\label{ch:introduction}

Pattern formation is a very active research topic in applied mathematics and natural sciences. Nonlinear patterns and waves describe many physical and biological phenomena and many advances in this field have been made in recent years. As some well-established references for using partial differential equations (PDEs) to model and analyze pattern formation, we refer to  \cite{PatternForm1,PatternForm2,PatternForm4}. 

However, for many processes arising in the natural sciences and engineering, one has to take into account stochastic effects for PDEs caused, e.g., by intrinsic finite-size fluctuations or external stochastic forcing. For finite-dimensional stochastic systems given by stochastic ordinary differential equations (SODEs), the mathematical analysis of patterns has progressed quite significantly~\cite{Khasminskii1,ArnoldSDE,BerglundGentz,LaingLord}. For spatial stochastic systems modelled by stochastic partial differential equations (SPDEs), less is known regarding the influence of noise on pattern formation~\cite{GarciaOjalvoSancho,Sagus2007SpatiotemporalOO}. At this point, several pattern formation scenarios have been studied quite a bit for SPDEs. For example, these include the emergence of patterns from homogeneous states near local bifurcation points with noise~\cite{Lythe,CaraballoLangaRobinson,Bloemker,Bloemker1,GowdaKuehn,KuehnRomano1,BlumenthalEngelNeamtu,mckane2014stochastic} and noise-induced metastability between deterministically locally stable steady states~\cite{BerglundGentz10,Barret,NewhallVandenEijnden}. Of course, going beyond steady states is important and in the last few decades, the cases of travelling waves and moving interfaces for SPDEs have taken center stage. The research on this topic started in the early 1980s~\cite{SchimanskyGeierMikhailovEbeling} and it has been growing quickly in recent years particularly for the Fisher-KPP equation~\cite{mueller1995random,DOERING2003243,conlon2005travelling,BrunetDerridaMuellerMunier,MuellerMytnikQuastel}, Nagumo-type SPDEs~\cite{krugerStochBist,Cartwright2019ACC,Hamster2019TravellingWF,geldhauser2020travelling,hamster2020stability,maclaurin2020metastability}, neural field integro-differential equations~\cite{HuttLongtin_additiveTuring,Kilpatrick_WanderingSNF,NeurField_Thul,lang2015multiscale,NeuralFields_MACLAURIN}, ecology \cite{constable2013stochastic} as well as regarding associated computational tools~\cite{LordThuemmler,Tuckwell1,SauerStannat}. For more detailed surveys including a larger-scale view of the literature on the effect of noise on travelling waves we refer to~\cite{Panja2003EffectsOF,Sagus2007SpatiotemporalOO,StochasticTravellingKuehn}.  

Although combining deterministic PDE theory for travelling waves with stochastic analysis techniques has been a clear success story, there is surprisingly little rigorous mathematical analysis for other time-dependent PDE patterns under the influence of noise. An obvious candidate to investigate further are rotating and spiral waves. For the deterministic PDE case, a lot of detailed analysis exists, which one can hope to leverage in an SPDE context. For stochastic rotating/spiral waves, several numerical results exist but a more detailed mathematical analysis is completely missing. In this work, we provide a first starting point to aim to close this gap focusing on the case of stochastic rotating waves. A key challenge we address is to understand the main differences between a simple symmetry group for the stochastic pattern, such as the translation group $\R$ for travelling waves, and more complex symmetry groups, such as the special Euclidean group $\mathrm{SE}\left(2\right)$ for spiral/rotating waves.    

Inspired by \cite{inglis2016general,maclaurin2020metastability}, our strategy to analyze stochastic rotating waves aims to adapt and extend a well-developed framework developed for stochastic travelling waves. In this strategy, one uses the travelling wave of the deterministic PDE as the structure around which a small noise analysis is performed. The main step starts with a linearization around a deterministic and locally stable travelling wave, which yields a linear SPDE. Then one splits the linearized SPDE dynamics into a part along the wave, which is governed for a travelling wave by a one-dimensional SODE, while the part transverse to the wave is described by an infinite-dimensional SPDE. Yet, the infinite-dimensional SPDE can often be shown to have good dissipation, or stability, properties over long time scales as one can combine the deterministic stability result with the assumption of sufficiently small noise. In summary, the SODE characterizes the phase diffusion along the wave, while the transverse linear infinite-dimensional SPDE can be used to establish stochastic stability results showing that stochastic waves do not leave a suitable neighbourhood of the deterministic wave over a long time scale. 

One can now hope that the same approach works as long as one can find suitable replacements for each step and still leverage deterministic stability. There are well-established results establishing stability results for rotating/spiral waves for PDEs~\cite{SANDSTEDE1997122,GolubitskyLeBlancMelbourne,SandstedeScheelWulff2,nonlinearstab_Rotwaves}, and one can obviously just still make the assumption of small noise. However, in this work we uncover two major challenges, when one tries to analyze stochastic rotating waves: (a) the SODE along the rotating wave becomes three-dimensional, and (b) the linearized SPDE dynamics contains a completely new class of algebraic terms arising due to the non-commutativity of $\mathrm{SE}\left(2\right)$ that interacts with the noise term in a non-trivial way. The challenge (a) is somewhat expected as the symmetry group of a pattern effectively encodes its neutral directions corresponding to spectrum on the imaginary axis for the deterministic PDE linearization. The challenge (b) is more subtle as one is very much used to employing commutativity for the translation group in the context of travelling waves for SPDEs, which substantially simplifies the interaction between noise and spatial differential operators. In this work, we study stochastic rotating waves, provide the adapted setup for the linearized stochastic analysis and obtain a first result on the behaviour of their dynamics. We achieve this using two different frameworks, the approximated variational phase, which works in a fixed reference frame and yields a local-in-time stability result, as well as the variational phase, which adapts in time to the pattern and hence yields a more global-in-time stability result. Although the second approach leads to a technically stronger result, the algebraic complexity is higher, so we decided to also present the technically more straightforward approximate variational phase. In summary, our results are probably not optimal yet as they are the first ones for stochastic rotating waves, but we expect that they provide a good baseline as well as a clear starting point for even further rigorous analysis of stochastic rotating/spiral waves.\medskip

The paper is structured as follows: In Section~\ref{sec:Setup}, we provide background and set up the notation. In Section~\ref{sec:rotwaves}, we summarize the necessary results from the deterministic PDE case regarding rotating waves and their stability. The new results are contained in~\ref{sec:stochRot}, where we derive the linearized dynamics for SPDEs around a rotating wave revealing new non-commutative terms as well as the reduced three-dimensional SODE along the pattern corresponding to the neutral spectrum. A key technical point in this context is to study the different possible notions of phase along the spiral. Finally, we also prove stability theorems for stochastic rotating waves. We conclude in Section~\ref{sec:conclusion} with a summary and an outlook to future challenges.  

\section{Background and Notation}
\label{sec:Setup}

\subsection{Stochastic Partial Differential Equations}

First of all we present some notions for SPDEs following~\cite{NeamtuStochDyn,daprato_zabczyk_2014,SPDERockner}. Let $\left(\mathcal{H},\langle\cdot, \cdot\rangle\right)$ be an infinite-dimensional separable Hilbert space with the induced norm $\|\cdot\|_{\mathcal{H}}$. We denote with $\mathcal{L}_{2}(\mathcal{H})$ the space of Hilbert-Schmidt operators on $\mathcal{H}$ and with $\mathcal{L}(\mathcal{H})$ the space of linear bounded operators on $\mathcal{H}$. 

Furthermore, we consider a complete probability space $(\Omega, \mathcal{F}, \mathbb{P})$, a filtration $\left(\mathcal{F}_{t}\right)_{t \geq 0}$ and we indicate with $\mathcal{B}(\cdot)$ the Borel-sigma algebra on a certain space, with $\otimes$  the product $\sigma$-algebra as well as the product measure, and with $\mathbb{E}$ the mathematical expectation. $\mathbb{P}-$a.s.~and a.e.~are respectively abbreviations for almost surely and almost everywhere. Additionally, for $1<p<\infty$ we denote with $L^{p}(D, \mathcal{H})$, where $D$ is a subset of $ \mathbb{R}^{d}$ for $d \geq 1$,  the space of $p$-integrable (over $D$) $\mathcal{H}$-valued functions and with $L^{p}(\Omega, \mathcal{F}, \mathbb{P} ; \mathcal{H})$ the space of $\mathcal{H}$-valued random variables on the probability space $(\Omega, \mathcal{F}, \mathbb{P})$ with finite moments of order $p .$  $W^{k, p}(D, \mathcal{H})$ denote, for $k \in \mathbb{R}$, the usual Sobolev spaces (i.e.~the spaces of functions in $L^{p}(D, \mathcal{H})$ with weak derivatives up to order $k$  again in $\left.L^{p}(D, \mathcal{H})\right)$. Moreover, we define $H^{k}(D, \mathcal{H})=W^{k, 2}(D, \mathcal{H})$.

On the stochastic basis $\left(\Omega, \mathcal{F},\left(\mathcal{F}_{t}\right)_{t \in[0, T]}, \mathbb{P}\right)$ we consider the SPDE
\begin{equation}\label{SPDEbasicCauchy}
\begin{cases}
\mathrm{d} u(t)=[A u(t)+f(t, u(t))] \mathrm{d} t+B(t, u(t)) \mathrm{d} W_{\mathcal{U}}(t), \quad t \in[0, T] \\
u(0)=u_{0} \in \mathcal{H}
\end{cases}
\end{equation}where the initial condition $u_{0}$ is $\mathcal{F}_{0}$-measurable, the linear closed densely defined operator $A$ generates a $C_{0}-$ semigroup $(\exp(tA))_{t \geq 0}$ on $\mathcal{H}$, the nonlinear terms $f:[0, T] \times \Omega \times \mathcal{H} \rightarrow \mathcal{H}$, $B:[0, T] \times \Omega \times \mathcal{H} \rightarrow \mathcal{L}_{2}(\mathcal{U}, \mathcal{H})$ are $\left(\mathcal{B}([0, t]) \otimes \mathcal{F}_{t} \otimes \mathcal{B}(\mathcal{H}) ; \mathcal{B}(\mathcal{H})\right)$-measurable, resp.~$\left(\mathcal{B}([0, t]) \otimes \mathcal{F}_{t} \otimes \mathcal{B}(\mathcal{H})\right. ; \left.\mathcal{B}\left(\mathcal{L}_{2}(\mathcal{U}, \mathcal{H})\right)\right)$-measurable for every $t \in [0, T]$ and $\left(W_{\mathcal{U}}(t)\right)_{t \in[0, T]}$ is an $\mathcal{U}$-cylindrical Brownian motion. For notational simplicity, the $\omega$-dependence of $f$ and $B$ has been omitted. A typical example is $u=(u(t))_{t \in[0, T]}=(u(t, x))_{t \in[0, T], x \in D}$ (where the space-dependence is usually dropped whenever there is no risk of confusion) where $A=\Delta$  is the Laplacian and $f$ and $B$ are sufficiently regular maps. To simplify the notation, we will often replace $u(t)$ by $u_t$ in the following.

In the remaining part of this section we provide an overview of mild solutions for SPDEs, an important solution notion for SPDEs of the form~\eqref{SPDEbasicCauchy}. 
\begin{definition}
 An $\mathcal{H}$-valued adapted process $(u(t))_{t \in[0, T]}$ is called a mild solution for \eqref{SPDEbasicCauchy} if
\begin{equation}\label{conditionStochMildsol}
\mathbb{P}\left(\int_{0}^{T}\|u(s)\|_{\mathcal{H}}^{2} \mathrm{~d} s<\infty\right)=1
\end{equation}and for $t \in[0, T]$ the variation of constants formula holds true:
\begin{equation}\label{stochMildsol}
\begin{aligned}
u(t)=&\exp(tA)u_{0}+\int_{0}^{t}\exp((t-s) A)f(s, u(s)) \mathrm{d} s\\
&+\int_{0}^{t}\exp((t-s) A)B(s, u(s)) \mathrm{d} W_{U}(s),\quad \mathbb{P}-a . s
\end{aligned}
\end{equation}
The first integral in \eqref{stochMildsol} will be referred to as deterministic convolution and the second one will be called stochastic convolution. 
\end{definition}

Using a classical fixed-point argument, one can show the following existence result for mild solutions of \eqref{SPDEbasicCauchy}, see~\cite[Thm.~7.2]{daprato_zabczyk_2014}.

\begin{theorem}\label{ExistenceUniquenessSPDE}
 Let $f$ and $B$ additionally satisfy the following Lipschitz and growth boundedness assumptions:
 \begin{enumerate}
 \item there exists a constant $L>0$ such that for all $u, v \in \mathcal{H}, t \in[0, T]$ and almost all $\omega \in \Omega$ we have:
$$
\|f(t,  u,\omega)-f(t,v,\omega)\|_{\mathcal{H}}+\|B(t, u,\omega)-B(t, v, \omega)\|_{\mathcal{L}_{2}(\mathcal{U}, \mathcal{H})} \leq L\|u-v\|_{\mathcal{H}}
$$
\item there exists a constant $l>0$ such that for all $u \in \mathcal{H}, t \in[0, T]$ and almost all $\omega \in \Omega$ we have:
$$
\|f(t, u, \omega)\|_{\mathcal{H}}^{2}+\|B(t, u, \omega)\|_{\mathcal{L}_{2}(\mathcal{U}, \mathcal{H})}^{2} \leq l^{2}\left(1+\|u\|_{\mathcal{H}}^{2}\right) .
$$
\end{enumerate}
Then there exists a unique (up to equivalence) mild solution of \ref{stochMildsol}.
\end{theorem}

Naturally, similar results can be obtained under weaker assumptions on the coefficients (e.g.~local Lipschitz continuity, dissipativity, local monotonicity etc.) using, for example, cut-off and localization techniques, see~\cite{daprato_zabczyk_2014,SPDERockner}.

\subsection{Special Euclidean Group}

In this section, we briefly introduce a particular type of Lie group, the special Euclidean group. We base our presentation on~\cite[Ch.~3]{skriptBeyn_Otten}. More detailed references and background can be found in~\cite{LieGroupsWarner,LinearGroupsLie}.

We start by introducing the Euclidean group. We have that every affine transformation $T: \mathbb{R}^{d} \rightarrow \mathbb{R}^{d}$ that preserves the Euclidean distance (i.e.~$\left.\left\|T p_{1}-T p_{2}\right\|_{2}=\left\|p_{1}-p_{2}\right\|_{2} \forall p_{1}, p_{2} \in \mathbb{R}^{d}\right)$
is of the form
$$
T p=Q p+b \quad \text { for some } Q \in \mathrm{O}\left(d\right), b \in \mathbb{R}^{d}\text {. }
$$
where $\mathrm{O}\left(d\right)$ is the orthogonal group. The composition of two transformations is
$$
\left(T_{1} \circ T_{2}\right) p=Q_{1}\left(T_{2} p\right)+b_{1}=Q_{1}\left(Q_{2} p+b_{2}\right)+b_{1}=Q_{1} Q_{2} p+Q_{1} b_{2}+b_{1}
$$
hence we have the group operation 
\begin{equation}\label{groupOperSE}
\left(Q_{1}, b_{1}\right) \circ\left(Q_{2}, b_{2}\right)=\left(Q_{1} Q_{2}, Q_{1} b_{2}+b_{1}\right)
\end{equation}
These transformations form the Euclidean group
$$
\mathrm{E}\left(d\right)=\left\{(Q, b): Q \in \mathrm{O}\left(d\right), b \in \mathbb{R}^{d}\right\}
$$
Taking $Q \in \mathrm{SO}\left(d\right)$, where $\mathrm{SO}\left(d\right)$ is the special orthogonal group, one obtains the subgroup $\mathrm{SE}\left(d\right)$, called the special Euclidean group. The group operation \eqref{groupOperSE}  differs from the standard group operation for products and for this reason the resulting group is called a semidirect product and written as
$$
\mathrm{E}\left(d\right)=\mathrm{O}\left(d\right) \ltimes \mathbb{R}^{d}, \quad \mathrm{SE}\left(d\right)=\mathrm{SO}\left(d\right) \ltimes \mathbb{R}^{d}
$$
It is possible to represent the group $\mathrm{SE}\left(d\right)$ as a subgroup of $\mathrm{GL}\left(\mathbb{R}^{d+1}\right)$ in the following way
\begin{equation}\label{SpecialEuclGroup}
\mathrm{SE}\left(d\right)=\left\{\left(\begin{array}{cc}
Q & b \\
0 & 1
\end{array}\right): Q \in \mathrm{SO}\left(d\right), b \in \mathbb{R}^{d}\right\}
\end{equation}
For the isomorphism between the two representations, note that the group operation \eqref{groupOperSE}  is represented by  the following matrix multiplication
$$
\left(\begin{array}{cc}
Q_{1} & b_{1} \\
0 & 1
\end{array}\right)\left(\begin{array}{cc}
Q_{2} & b_{2} \\
0 & 1
\end{array}\right)=\left(\begin{array}{cc}
Q_{1} Q_{2} & Q_{1} b_{2}+b_{1} \\
0 & 1
\end{array}\right)
$$
Notice also that this Lie group is a linear group and is not Abelian for every $d\geq 1$. A chart covering all of $\mathrm{SE}\left(2\right)$ is 
\begin{equation}
\label{chartSE}
\begin{aligned}[t]
\gamma\colon \mathbb{R}/(2\pi\mathbb{Z})\times \mathbb{R}^2&\longrightarrow \mathrm{SE(2)}\\
(\theta,\bar{b}) &\longmapsto \left(\begin{bmatrix}
\cos \theta & -\sin \theta \\
\sin \theta & \cos \theta
\end{bmatrix},\bar{b}\right)
\end{aligned}
\end{equation}
where $\bar{b}=(b_1,b_2)$.

\section{Deterministic Rotating Waves}
\label{sec:rotwaves}

Rotating waves are non-stationary solutions of PDEs with particular dynamical properties. We will consider here the particular case of rotating waves for reaction-diffusion systems of the form

\begin{equation}\label{react_diff2D}
  	\partial_t u(\bar{x},t)=D \Delta u(\bar{x},t)+f(u(\bar{x},t)), \quad \bar{x}=\begin{bmatrix}
  	x_1 \\
    x_2
  	\end{bmatrix}\in \mathbb{R}^2,t>0
  	\end{equation}
where $D$ is a diagonal matrix with positive entries on the diagonal and $\Delta$ is the two dimensional Laplacian, that is $\Delta=\partial_{x_1 x_1}+\partial_{x_2 x_2}$. In the following we will assume
\begin{assumption}\label{reg_f_rot}
$f\in C^4(\mathbb{R}^N,\mathbb{R}^N)$ with first, second and third derivatives bounded.
\end{assumption}

In Euclidean coordinates, a rotating wave is a solution of \eqref{react_diff2D} of the form
\begin{equation}
    u(\bar{x},t)=u_*\left(R_{-\omega_*t}\bar{x}\right)
\end{equation}
for some $\omega_*\in \mathbb{R}$, where we have the rotation matrix
\begin{equation*}
    R_{\theta}=\begin{bmatrix}
    \cos{\theta} & -\sin{\theta}\\
    \sin{\theta} & \cos{\theta}
    \end{bmatrix}.
\end{equation*}
It is important to observe that if $u_*\left(R_{-\omega_*t}(\cdot)\right)$ is a rotating wave then also
\begin{equation}
u_*\left(R_{-\omega_*t}R_{-\theta}(\cdot-\bar{b})\right)=u_*\left(R_{-\omega_*t-\theta}(\cdot-\bar{b})\right), \text{   } \theta \in \mathbb{R}/(2\pi\mathbb{Z}) \text{ and } \bar{b}\in \mathbb{R}^2
\end{equation}is still a rotating wave solution of \eqref{react_diff2D}. Therefore, translating and rotating a rotating wave will still give us a rotating wave. With this observation we have an entire family of rotating waves parametrized by the special Euclidean group $SE(2)$.

It will be convenient in some cases to consider the equation \eqref{react_diff2D} in polar coordinates that is 
\begin{equation}\label{react_diff2D_polar}
  	\partial_t u(r,\phi,t)=D \Delta u(r,\phi,t)+f(u(r,\phi,t)), \quad r>0,\phi\in  \mathbb{R}/(2\pi\mathbb{Z}) \text{ and } t>0,
  	\end{equation}
 where the Laplacian in polar coordinates is $\Delta=\frac{1}{r}\partial_r+\partial_{rr}+\frac{1}{r^2}\partial_{\psi \psi}$. In this new coordinates system, rotating waves are solutions of the form
 \begin{equation}
     u(r,\phi,t)=u_*(r,\phi-\omega_*t)
 \end{equation}
 We now consider what in physics is called the change to a co-rotating reference frame, that in Euclidean coordinates  is the change of variables 
$$
\bar{y}=\begin{bmatrix}
 y_1\\
 y_2
\end{bmatrix}=R_{-\omega_*t}\begin{bmatrix}x_1\\
 x_2
 \end{bmatrix}
$$
 Equivalently, we can consider polar coordinates and in this case, the change to the co-rotating reference frame is the change of variables
 $$
\begin{bmatrix}
r\\
 \psi
 \end{bmatrix}=\begin{bmatrix}
 r\\
 \phi-\omega_*t
 \end{bmatrix}
 $$
 Passing to a co-rotating reference frame equation \eqref{react_diff2D_polar} becomes 
 \begin{equation}
 \label{react_diff2D_corotat}
  	\partial_t u(r,\psi,t)=D \Delta u(r,\psi,t)+\omega_*\partial_{\psi}u(r,\psi,t)+f(u(r,\psi,t))
  	\end{equation}
 for $ r>0$, $\psi\in \mathbb{R}/(2\pi\mathbb{Z})$ and $t\geq 0$. It is important to notice that this change of reference frame changes the regularity properties of the PDE \eqref{react_diff2D_polar}. In fact, in Euclidean coordinates, we have by the Leibniz rule
 \begin{equation}
 \label{tangential_der}
     \partial_{\psi}=y_1 \partial_{y_2}-y_2 \partial_{y_1}
 \end{equation} and therefore $D\Delta+\omega_*\partial_{\psi}$ is a second-order operator but with unbounded coefficients. We will see in the following chapters that, in fact, $D\Delta+\omega_*\partial_{\psi}$ is not a sectorial operator and therefore does not generate an analytic semigroup. For a rotating wave $u_*$ we have
  \begin{equation}\label{rot_wav_eq}
      0=D \Delta u_*(r,\psi)+\omega_*\partial_{\psi}u_*(r,\psi)+f(u_*(r,\psi))
  \end{equation}
It follows that rotating waves are stationary solutions of \eqref{react_diff2D_corotat}, i.e., they are stationary solutions of the reaction-diffusion equation in the co-rotating frame. We will focus in the following on a particular type of rotating waves, asymptotically constant rotating waves, and discuss some stability properties of these type of rotating waves. We recall the following definition:

\begin{definition}\label{DefAsyRotWave}
A rotating wave such that  for some constant vector $u_{\infty}\in \mathbb{R}^N$  
\begin{equation}
    \sup_{|x|\geq R}|u_*(x)-u_{\infty}|\rightarrow 0 \text{  as  } R\rightarrow+\infty
\end{equation}
is called asymptotically constant rotating wave with asymptotic constant $u_{\infty}$.
\end{definition}

In the following, we will assume:
\begin{assumption}\label{regularity_rotWaves}
$u_*$ is twice differentiable. Additionally, $$u_*-u_{\infty}\in H^2(\mathbb{R}^2,\mathbb{R}^N)$$
and $\partial_{y_1}u_*$,$\partial_{y_2}u_*$ and $\partial_{\psi}u_*\in H^4_{eucl}(\mathbb{R}^2,\mathbb{R}^N)$.
\end{assumption}

In addition, we assume also:

\begin{assumption}\label{NegDefMat}
The matrix $f'(u_{\infty})$  is negative definite, where $f'$ is the derivative of the non-linearity $f$ in \eqref{react_diff2D} and $u_{\infty}\in \mathbb{R}^N$ is the asymptotic constant from Definition \ref{DefAsyRotWave}.
\end{assumption}

Without loss of generality, we will always consider in the following $u_{\infty}=0$.

\subsection{Linearization Operator}

For a stationary solution $u_*$ of a (partial) differential equation it is natural to consider the linearization at $u_*$ to study stability. In particular, if $u_*$ is rotating wave then it is a stationary solution of \eqref{react_diff2D_corotat} we consider the Cauchy problem
\begin{equation}
\label{react_diff_comoving_pert}
\begin{cases}
\partial_t u(t)=D\Delta u(t)+\omega_*\partial_{\psi} u(t)+f(u(t)), \quad y,t\in \mathbb{R} \\
u_0=u_*+v_0
\end{cases}
\end{equation}where $v_0$ is a general small perturbation in $L^2(\mathbb{R},\mathbb{R})$ for the moment. Now we obtain that the difference $w(t)=u(t)-u_*$ between the solution to the Cauchy problem with perturbed initial condition \eqref{react_diff_comoving_pert} and the travelling wave solution satisfies the equation
$$
\begin{aligned}
    \partial_t w(t)
    =&D\Delta u(t)+\omega_*\partial_{\psi} u(t)+f(u(t))-\left(D\Delta u_*+\omega_*\partial_{\psi} u_*+f(u_*)\right)\\
    =& D\Delta w(t)+\omega_*\partial_{\psi} w(t)+f(u_*+w(t))-f(u_*)\\
    =&  D\Delta w(t)+\omega_*\partial_{\psi} w(t)+f'(u_*)w(t)+O\left(||w(t)||^2\right)
\end{aligned}
$$
using \eqref{rot_wav_eq} and \eqref{react_diff_comoving_pert}. $w_t$ is the evolution of the initial perturbation over time. Therefore, we obtain that the linearized PDE is 
\begin{equation}
   \partial_t w(t)=\mathcal{L}_*w(t)
\end{equation}
where $\mathcal{L}_*$ is the linearization operator that, in this case, is 
\begin{equation}\label{linearizationRotat}
    \mathcal{L}_*=D\Delta+\partial_{\psi}+f'(u_*)
\end{equation}
where $u_*$ is the asymptotically constant rotating wave. The linearization operator $\mathcal{L}_*$ is a closed and densely defined operator on $H^2\left(\mathbb{R}^2,\mathbb{R}^N\right)$ with domain 
\benn
H^4_{eucl}\left(\mathbb{R}^2,\mathbb{R}^N\right)=\{u\in H^4\left(\mathbb{R}^2,\mathbb{R}^N\right): \partial_{\psi}u\in H^2\left(\mathbb{R}^2,\mathbb{R}^N\right)\}.
\eenn
We notice via direct calculations (e.g.~using $\partial_{\psi}=y_1\partial_{y_2}-y_2\partial_{y_1}$) that 
\benn
\mathcal{L}_{*}\left(\partial_{y_1} u_*\right)  =-\omega_{*}\partial_{y_2}u_{*},\quad 
\mathcal{L}_{*}\left(\partial_{y_2} u_*\right)
=\omega_{*}\partial_{y_1}u_{*},\quad\mathcal{L}_{*}\left(\partial_{\psi} u_*\right)=0.
\eenn
From the previous equalities it follows easily that $0$, $+\mathrm{i}\omega_*$ and $-\mathrm{i}\omega_*$  are eigenvalues with eigenfunctions $\partial_{\psi} u_*$, $\left(\partial_{y_1}+ \mathrm{i}\partial_{y_2}\right)u_{*}$ and $\left(\partial_{y_1}- \mathrm{i}\partial_{y_2}\right)u_{*}$, e.g., one checks
\benn
\mathcal{L}_{*}\left(\left(\partial_{y_1}+ \mathrm{i}\partial_{y_2}\right)u_{*}\right) 
=\mathrm{i}\omega_{*}\Big(\left(\partial_{y_1}+ \mathrm{i}\partial_{y_2}\right)u_{*}\Big)
\eenn
The functions $\partial_{\psi}u_{*}$, $\left(\partial_{y_1}+ \mathrm{i}\partial_{y_2}\right)u_{*}$ and $\left(\partial_{y_1}- \mathrm{i}\partial_{y_2}\right)u_{*}$ are indeed eigenfunctions as from \eqref{regularity_rotWaves} it follows that they are in $H^4_{eucl}(\mathbb{R}^2,\mathbb{C}^N)$. We need now some additional spectral stability assumptions: 

\begin{assumption}\label{SpectAsyRot}
The point spectrum $\sigma_{pt}\left(\mathcal{L}_{*}\right)$ of $\mathcal{L}_{*}$ on $H^{2}\left(\mathbb{R}^2, \mathbb{R}^N\right) $ is such that
\begin{equation}\label{spectrum_rot}
\sigma_{pt}\left(\mathcal{L}_{*}\right) \subset\{\lambda \in \mathbb{C}: \mathfrak{R e}(\lambda) \leq-b\} \cup\{0,+\mathrm{i}\omega_*,-\mathrm{i}\omega_*\}
\end{equation}for some positive constant $b$. Moreover, we assume that the eigenvalue $+\mathrm{i}\omega_*,-\mathrm{i}\omega_*$ and $0$ have all multiplicity $1$ and the corresponding eigenfunctions are $\left(\partial_{y_1}+ \mathrm{i}\partial_{y_2}\right)u_{*}$, $\left(\partial_{y_1}- \mathrm{i}\partial_{y_2}\right)u_{*}$ and $\partial_{\psi}u_{*}$.
\end{assumption}
From Assumption \ref{SpectAsyRot} we already have, without loss of generality, that the essential spectrum of $\mathcal{L}_*$ is contained in $\{\lambda \in \mathbb{C}: \mathfrak{R e}(\lambda) \leq-b\}$, see \cite{nonlinearstab_Rotwaves}, and therefore we have for the entire spectrum that  
\begin{equation}
\sigma\left(\mathcal{L}_{*}\right) \subset\{\lambda \in \mathbb{C}: \mathfrak{R e}(\lambda) \leq-b\} \cup\{0,+\mathrm{i}\omega_*,-\mathrm{i}\omega_*\}
\end{equation}

As the strongly continuous semigroup $\exp(\mathcal{L}_*t)$ generated by $\mathcal{L}_*$ is not an analytic semigroup, the spectrum does not necessarily determine exponential decay bounds for the semigroup as in the traveling waves case. We will see in the following sections that, in this case, spectral stability \eqref{SpectAsyRot} is still enough to obtain exponential decay bounds.

\subsection{Projections on the Center Space}

It is helpful~\cite{nonlinearstab_Rotwaves} to introduce the adjoint operator 
\begin{equation*}
    \mathcal{L}_*^{ad}=D\Delta-\omega_*\partial_{\psi}+\left(f'(u_*)\right)^\top.
\end{equation*}
Assumption \ref{NegDefMat} implies that the essential spectrum of $\mathcal{L}_*$ lies in the subset $\{\lambda \in \mathbb{C}:$ $\mathfrak{R e}(\lambda) \leq-b\}$. This means that $\mathcal{L}_{*}$ is Fredholm of zero index, and therefore for the Fredholm alternative the kernel of the operators $\mathcal{L}_{*}^{ad}$, $\mathcal{L}_{*}^{ad}-i\omega_{*}=\left(\mathcal{L}_{*}+i\omega_{*}\right)^{ad}$ and $\mathcal{L}_{*}^{ad}+\mathrm{i}\omega_{*}=\left(\mathcal{L}_{*}-\mathrm{i}\omega_{*}\right)^{ad}$ is $1$ -dimensional. Therefore, we have one eigenfunction $\phi^{0}$ for  $\mathcal{L}_{*}^{ad}$ for the eigenvalue $0$, one eigenfunction $\phi^{\mathrm{i}\omega_*}$ for the eigenvalue $\mathrm{i}\omega_*$ and one eigenfunction $\phi^{-\mathrm{i}\omega_*}$ for the eigenvalue $-\mathrm{i}\omega_*$. 
Without loss of generality (in case divide the eigenfunction of the adjoint by the scalar product of the eigenfunction of the adjoint with the respective eigenfunction of the original linearization operator) we can assume 
\benn
\left\langle\phi^{0}, \partial_{\psi}u_* \right\rangle=1,\quad 
\left\langle\phi^{\mathrm{i}\omega_*}, \left(\partial_{y_1}-\mathrm{i}\partial_{y_2}\right)u_*\right\rangle=1,\quad 
\left\langle\phi^{-\mathrm{i}\omega_*}, \left(\partial_{y_1}+ \mathrm{i}\partial_{y_2}\right)u_* \right\rangle=1 
\eenn
and it is easy to show that
$$
\begin{aligned}
0=&\left\langle\phi^{\mathrm{i}\omega_*}, \partial_{\psi}u_* \right\rangle=\left\langle\phi^{-\mathrm{i}\omega_*}, \partial_{\psi}u_* \right\rangle=\left\langle\phi^{\mathrm{i}\omega_*}, \left(\partial_{y_1}+\mathrm{i}\partial_{y_2}\right)u_*\right\rangle\\
=&\left\langle\phi^{-\mathrm{i}\omega_*}, \left(\partial_{y_1}-\mathrm{i}\partial_{y_2}\right)u_*\right\rangle=\left\langle\phi^{0}, \left(\partial_{y_1}-\mathrm{i}\partial_{y_2}\right)u_* \right\rangle=\left\langle\phi^{0}, \left(\partial_{y_1}+ \mathrm{i}\partial_{y_2}\right)u_* \right\rangle
\end{aligned}
$$
Hence, one can decompose the spaces $H^{2}\left(\mathbb{R}^2, \mathbb{R}^N\right)$ and $H^{4}_{eucl}\left(\mathbb{R}^2, \mathbb{R}^N\right)$ in the following way. We define 
\begin{align*}
\Phi &=\text{span}\{\partial_{\psi} u_*,\left(\partial_{y_1}+ \mathrm{i}\partial_{y_2}\right)u_{*},\left(\partial_{y_1}-\mathrm{i}\partial_{y_2}\right)u_{*}\}=N(\mathcal{L}_*),\\
\Psi &=\text{span} \{\phi^{0},\phi^{\mathrm{i}\omega_*},\phi^{-\mathrm{i}\omega_*}\}\text{ and }\\
W&=\Psi^{\perp}\\
W_{eucl}&=H^{4}_{eucl}\left(\mathbb{R}^2, \mathbb{R}^N\right)\cap W
\end{align*}
where the orthogonal complement $\Psi^{\perp}$ of $\Psi$ is taken in $H^{2}\left(\mathbb{R}^2, \mathbb{R}^N\right)$ and $N(\mathcal{L}_*)$ is the kernel of $\mathcal{L}_*$. In this setting we have the following decompositions
\begin{align*}
H^{2}\left(\mathbb{R}^2, \mathbb{R}^N\right)=\Phi\oplus W\qquad 
H^{4}_{eucl}\left(\mathbb{R}^2, \mathbb{R}^N\right)=\Phi\oplus W_{eucl}
\end{align*}In fact, given $h\in H^{2}\left(\mathbb{R}^2, \mathbb{R}^N\right)$ we can define 
\begin{equation*}
w=h-\left\langle \phi^{0},h\right \rangle\partial_{\psi} u_*-\left\langle \phi^{-\mathrm{i}\omega_*},h\right \rangle\left(\partial_{y_1}+ \mathrm{i}\partial_{y_2}\right)u_{*}-\left\langle \phi^{\mathrm{i}\omega_*},h\right \rangle\left(\partial_{y_1}- \mathrm{i}\partial_{y_2}\right)u_{*}
\end{equation*}and we have that $w\in W=\Psi^{\perp} $ as $\left\langle \phi^{0},w\right \rangle=\left\langle \phi^{\mathrm{i}\omega_*},w\right \rangle=\left\langle \phi^{-\mathrm{i}\omega_*},w\right \rangle=0$. Moreover, assuming $w\in W=\Psi^{\perp} $ we have that 
\begin{equation}\label{indepEq}
0=w+a_1\partial_{\psi} u_*+a_2\left(\partial_{y_1}+ \mathrm{i}\partial_{y_2}\right)u_*+a_3\left(\partial_{y_1}-\mathrm{i}\partial_{y_2}\right)u_*
\end{equation}
where $a_1,a_2$ and $a_3\in \mathbb{C}$ implies $a_1=a_2=a_3=0$. This follows easily applying the scalar product with the eigenfunctions $\phi^{0},\phi^{\mathrm{i}\omega_*}$ and $\phi^{-\mathrm{i}\omega_*}$ to both sides of equation \eqref{indepEq}. The same argument works for the decomposition of $H^{4}_{eucl}\left(\mathbb{R}^2, \mathbb{R}^N\right)$. We can now define
\begin{align}\label{functionsEta}
\eta^{0}&=\phi^{0}\\
\eta^{1}&=\frac{1}{2}\left(\phi^{\mathrm{i}\omega_*}+\phi^{-\mathrm{i}\omega_*}\right)\\
\eta^{2}&=\frac{1}{2}\mathrm{i}\left(\phi^{\mathrm{i}\omega_*}-\phi^{-\mathrm{i}\omega_*}\right)
\end{align}
and we get, adapting the notation defining $\partial_0=\partial_{\psi}, \partial_1=\partial_{y_1}$ and $\partial_2=\partial_{y_2}$, 
\benn
\left\langle\eta^{i}, \partial_ju_*\right\rangle=0 \text { if } i \neq j,\qquad \left\langle\eta^{i}, \partial_iu_*\right\rangle=1
\eenn
where $i,j \in \{0,1,2\}$. The next natural step is to define the projection on the subspace $\Phi$ spanned by the eigenspaces of the eigenvalues $0, +\mathrm{i}\omega_*$ and $-\mathrm{i}\omega_*$ of the linearization operator $\mathcal{L}_*$
\begin{eqnarray}\label{ProjectRotat}
&P :  H^{2}\left(\mathbb{R}^2, \mathbb{R}^N\right)\longrightarrow \Phi \\
&h \mapsto  Ph=\sum_{i=0}^2 \left\langle\eta^{i} , h\right\rangle\partial_iu_*
\end{eqnarray}
Notice in fact that $$Ph=\left\langle \phi^{0},h\right \rangle\partial_{\psi} u_*+\left\langle \phi^{-\mathrm{i}\omega_*},h\right \rangle\left(\partial_{y_1}+ \mathrm{i}\partial_{y_2}\right)u_{*}+\left\langle \phi^{\mathrm{i}\omega_*},h\right \rangle\left(\partial_{y_1}- \mathrm{i}\partial_{y_2}\right)u_{*}.$$
We consider also the associated projection on $W=\Psi^{\perp}$.
\begin{eqnarray*}
&I-P :  H^{2}\left(\mathbb{R}^2, \mathbb{R}^N\right)\longrightarrow W \\
&h \mapsto  (I-P)h=h-\sum_{i=0}^2 \left\langle\eta^{i} , h\right\rangle\partial_iu_*
\end{eqnarray*}
In the case of asymptotically rotating waves we also have that the semigroup $\exp(\mathcal{L}_*t)$, generated by the operator $\mathcal{L}_*$ satisfies an exponential decay bound. In fact, we have the following result from \cite{nonlinearstab_Rotwaves}:

\begin{theorem}
\label{expBoundAsyRot}
Under Assumptions \ref{NegDefMat} and \eqref{spectrum_rot} the operator $\mathcal{L}_*$ generates a strongly continuous semigroup on $ H^{2}\left(\mathbb{R}^2, \mathbb{R}^N\right)$. For any $0\leq \omega< b$ (with $b$ from Assumption \ref{SpectAsyRot}) there exists a constant $C=C(\tau)$ such that 
\begin{equation*}
    \|\left.\exp(\mathcal{L}_*t)\right|_{ R(I-P)}h\|_{H^{2}}\leq C\exp(-\omega t)||h||_{H^{2}}
\end{equation*} for all $h\in R(I-P)=W$ the range of $I-P$.
\end{theorem}

We additionally notice that 
\begin{equation}
\label{modifiedEigen}
\mathcal{L}_*^{ad}\eta_0=\omega_*\eta_0,\qquad 
\mathcal{L}_*^{ad}\eta_1=\omega_*\eta_2,\qquad
\mathcal{L}_*^{ad}\eta_2=-\omega_*\eta_1.
\end{equation}

\subsection{Lie Group Representation}

We want now to introduce the family isomorphisms $(\mathcal{T}_{\gamma})_{\gamma\in SE(2)}$. An element $\gamma$ of $SE(2)$ can be represented as $\gamma(\theta,b_1,b_2)=\left(R_{\theta},\bar{b}\right)$ where $R_{\theta}\in SO(2)$ is a rotation matrix and $\bar{b}=\begin{bmatrix}b_1 ~ b_2\end{bmatrix}^\top\in \mathbb{R}^2$. Then one defines
\begin{eqnarray*}
\mathcal{T}_{\gamma} : H^2(\mathbb{R}^2,\mathbb{R}^N) &\rightarrow& H^2(\mathbb{R}^2,\mathbb{R}^N) \\
u &\mapsto& \mathcal{T}_{\gamma}u = u\left(R_{-\theta}\left(\cdot-\bar{b}\right)\right)
\end{eqnarray*}
From $1=\sin(\theta)^2+\cos(\theta)^2$, the following expressions
$$
\begin{aligned}
\partial_{y_1}\mathcal{T}_{\gamma}u_*(\bar{y})
=&\partial_{y_1}u_*(R_{\theta}^{-1}(\bar{y}-\bar{b}))
=\partial_{y_1}u_* \left(\begin{bmatrix}
\cos(\theta) & \sin(\theta) \\
-\sin(\theta) & \cos(\theta)
\end{bmatrix} (\bar{y}-\bar{b})\right)\\
=&\left(\nabla u_* \left(\begin{bmatrix}
\cos(\theta) & \sin(\theta) \\
-\sin(\theta) & \cos(\theta)
\end{bmatrix} (\bar{y}-\bar{b})\right)\right)^T\begin{bmatrix}
\cos(\theta) & \sin(\theta) \\
-\sin(\theta) & \cos(\theta)
\end{bmatrix}\begin{bmatrix}
1 \\
0
\end{bmatrix}\\
=&\left(\mathcal{T}_{\gamma}\nabla u_* \left(\bar{y}\right)\right)^T\begin{bmatrix}
\cos(\theta) & \sin(\theta) \\
-\sin(\theta) & \cos(\theta)
\end{bmatrix}\begin{bmatrix}
1 \\
0
\end{bmatrix}\\
\end{aligned}
$$and similarly
$$
\begin{aligned}
\partial_{y_2}\mathcal{T}_{\gamma}u_*(\bar{y})
=&\left(\mathcal{T}_{\gamma}\nabla u_* \left(\bar{y}\right)\right)^T\begin{bmatrix}
\cos(\theta) & \sin(\theta) \\
-\sin(\theta) & \cos(\theta)
\end{bmatrix}\begin{bmatrix}
0 \\
1
\end{bmatrix}\\
\end{aligned}
$$
and similar computations for the second derivatives it follows that $\mathcal{T}_{\gamma}$ are isometries of $H^2(\mathbb{R}^2,\mathbb{R}^N)$. This family of isometries can also be seen as the Lie group representation 
\begin{eqnarray*}
\mathcal{T}_{(\cdot)} : SE(2) &\rightarrow& GL\left(H^2(\mathbb{R}^2,\mathbb{R}^N)\right) \\
\gamma &\mapsto& \left[u\mapsto\mathcal{T}_{\gamma}u= u\left(R_{-\theta}\left(\cdot-\bar{b}\right)\right)\right]
\end{eqnarray*}
Given a function $u$ we can define the map
\begin{eqnarray}
\label{IsometrySpeciaEucl}
\mathcal{T}_{(\cdot)}u : \mathbb{R} &\rightarrow&H^2(\mathbb{R}^2,\mathbb{R}^N) \\
\gamma &\mapsto& \mathcal{T}_{\gamma}u= u\left(R_{-\theta}\left(\cdot-\bar{b}\right)\right)
\end{eqnarray}
A particularly interesting case of the previous map will be $\mathcal{T}_{(\cdot)}u_*$ where $u_*$ is a rotating wave profile. In fact, we have that $\mathcal{T}_{\gamma}u_*(R_{-\omega t}(\cdot ))=u_*(R_{-\omega t}R_{-\theta}\left(\cdot-\bar{b}\right))$ is still a rotating wave as mentioned before. If we assume the function $u$ to be twice continuously differentiable (like in the case of rotating waves) we can also take derivatives of the map~\eqref{IsometrySpeciaEucl}. We have 
$$
\begin{aligned}
\partial_{b_1}\mathcal{T}_{\gamma}u_*(\bar{y})
=&\partial_{b_1}u_*(R_{\theta}^{-1}(\bar{y}-\bar{b}))
=\partial_{b_1}u_* \left(\begin{bmatrix}
\cos(\theta) & \sin(\theta) \\
-\sin(\theta) & \cos(\theta)
\end{bmatrix} (\bar{y}-\bar{b})\right)\\
=&-\left(\nabla u_* \left(\begin{bmatrix}
\cos(\theta) & \sin(\theta) \\
-\sin(\theta) & \cos(\theta)
\end{bmatrix} (\bar{y}-\bar{b})\right)\right)^T\begin{bmatrix}
\cos(\theta) & \sin(\theta) \\
-\sin(\theta) & \cos(\theta)
\end{bmatrix}\begin{bmatrix}
1 \\
0
\end{bmatrix}\\
=&-\left(\mathcal{T}_{\gamma}\nabla u_* \left(\bar{y}\right)\right)^T\begin{bmatrix}
\cos(\theta) & \sin(\theta) \\
-\sin(\theta) & \cos(\theta)
\end{bmatrix}\begin{bmatrix}
1 \\
0
\end{bmatrix}\\
\end{aligned}
$$
and similarly 
$$
\begin{aligned}
\partial_{b_2}\mathcal{T}_{\gamma}u_*(\bar{y})
=&-\left(\mathcal{T}_{\gamma}\nabla u_* \left(\bar{y}\right)\right)^T\begin{bmatrix}
\cos(\theta) & \sin(\theta) \\
-\sin(\theta) & \cos(\theta)
\end{bmatrix}\begin{bmatrix}
0 \\
1
\end{bmatrix}\\
\end{aligned}
$$
Moreover, we have
$$
\begin{aligned}
\partial_{\theta}\mathcal{T}_{\gamma}u_*(\bar{y})
=&\partial_{\theta}u_*(R_{\theta}^{-1}(\bar{y}-b))
=\partial_{\theta}u_* \left(\begin{bmatrix}
\cos(\theta) & \sin(\theta) \\
-\sin(\theta) & \cos(\theta)
\end{bmatrix} (\bar{y}-b)\right)\\
=&\nabla u_* \left(\begin{bmatrix}
\cos(\theta) & \sin(\theta) \\
-\sin(\theta) & \cos(\theta)
\end{bmatrix} (\bar{y}-b)\right)\begin{bmatrix}
0 & 1 \\
-1 & 0
\end{bmatrix}\begin{bmatrix}
\cos(\theta) & \sin(\theta) \\
-\sin(\theta) & \cos(\theta)
\end{bmatrix}(\bar{y}-b)\\
=&\mathcal{T}_{\gamma}\left(\nabla u_*\left(\bar{y}\right)\begin{bmatrix}
0 & 1 \\
-1 & 0
\end{bmatrix}(\bar{y})\right)\\
=&-\mathcal{T}_{\gamma}\left(y_1\partial_{y_2}-y_2\partial_{y_1}\right)u_* \left(\bar{y}\right)=-\mathcal{T}_{\gamma}\left(\partial_\psi u_*(\bar{y})\right)
\end{aligned}
$$
We additionally notice that the integration by parts formulas hold 
\begin{equation}\label{integByparts}
\begin{array}{l}
\left\langle\partial_{b_1}\mathcal{T}_{\gamma}u(\bar{y}), \mathcal{T}_{\gamma}v(\bar{y})\right\rangle=-\left\langle\mathcal{T}_{\gamma}u(\bar{y}), \partial_{b_1}\mathcal{T}_{\gamma}v(\bar{y})\right\rangle
\\
\left\langle\partial_{b_2}\mathcal{T}_{\gamma}u(\bar{y}), \mathcal{T}_{\gamma}v(\bar{y})\right\rangle=-\left\langle\mathcal{T}_{\gamma}u(\bar{y}), \partial_{b_2}\mathcal{T}_{\gamma}v(\bar{y})\right\rangle
\\
\left\langle\partial_{\theta}\mathcal{T}_{\gamma}u(\bar{y}), \mathcal{T}_{\gamma}v(\bar{y})\right\rangle=-\left\langle\mathcal{T}_{\gamma}u(\bar{y}), \partial_{\theta}\mathcal{T}_{\gamma}v(\bar{y})\right\rangle\\
\end{array}
\end{equation}
This trivially follows from the classical integration by parts formulas 
$$\begin{aligned}
&\left\langle\partial_{y_1}u(\bar{y}), v(\bar{y})\right\rangle=-\left\langle u(\bar{y}), \partial_{y_1}v(\bar{y})\right\rangle\\
&\left\langle\partial_{y_2}u(\bar{y}), v(\bar{y})\right\rangle=-\left\langle u(\bar{y}), \partial_{y_2}v(\bar{y})\right\rangle \\
&\left\langle\partial_{\psi}u(\bar{y}), v(\bar{y})\right\rangle=-\left\langle u(\bar{y}), \partial_{\psi}v(\bar{y})\right\rangle
\end{aligned}
$$ 
We can now define the operator   
\begin{eqnarray*}
F : H_{eucl}^4(\mathbb{R}^2,\mathbb{R}^N) &\rightarrow& H^2(\mathbb{R}^2,\mathbb{R}^N) \\
u &\mapsto& D\Delta u+f(u).
\end{eqnarray*}
The operator $F$ is equivariant with respect to the action $\mathcal{T}_{(\cdot)}$ of $SE(2)$, that is 
\begin{equation*}
F(\mathcal{T}_{\gamma}u)=\mathcal{T}_{\gamma}F(u)
\end{equation*}
for $\gamma\in SE(2)$ and $u\in H^4_{eucl}(\mathbb{R}^2,\mathbb{R}^N) $. This easily follows from the fact that the Laplacian operator is rotation invariant (see Lemma 32.3 in \cite{10.5555/3360180}). We can rewrite \eqref{react_diff2D} in the form
\begin{equation*}
\partial_t u(t)=F(u(t)).
\end{equation*}
For rotating waves, we have only that translations of a rotating wave are rotating waves too. These new rotating waves are stationary points with respect to a different co-rotating reference frame, with a different center of rotation with respect the original rotating wave. Therefore, they are not solutions of \eqref{react_diff2D_corotat}. However, in the co-rotating reference frame with respect to the original rotating wave $\bar{y}=R_{-\omega t}\bar{x}$ we have that $\mathcal{T}_{\gamma}u_*$ is a stationary solution of
\begin{equation}\label{react_diff2D_corotat_Gamma}
  	\partial_t u(r,\psi,t)=D \Delta u(r,\psi,t)+\omega_*\mathcal{T}_{\gamma}\partial_{\psi}\mathcal{T}_{-\gamma}u(r,\psi,t)+f(u(r,\psi,t)) 
  	\end{equation}
where  $r>0,\psi\in \mathbb{R}/(2\pi\mathbb{Z})$  and $t\geq0$.
\begin{equation}\label{eq: off-center stationary solution}
D \Delta  \mathcal{T}_{\gamma}u_*+\omega_*\mathcal{T}_{\gamma}\partial_{\psi}\mathcal{T}_{-\gamma} \mathcal{T}_{\gamma}u_*+f(\mathcal{T}_{\gamma}u_*)=\mathcal{T}_{\gamma}\left(D\Delta u_*+\omega_* \partial_{\psi} u_*+f(u_*)\right)=0
\end{equation}
Therefore, we can consider the linearization operator at $\mathcal{T}_{\gamma}u_*$ that is 
\be
\label{linearizationOperRotatGamma}
\mathcal{L}_{\gamma}=D \Delta +\omega_*\mathcal{T}_{\gamma}\partial_{\psi}\mathcal{T}_{-\gamma}+f'(\mathcal{T}_{\gamma}u_*)=
\mathcal{T}_{\gamma}\left(D\Delta u_*+\omega_* \partial_{\psi} u_*+f'(u_*)\right)\mathcal{T}_{-\gamma}=
\mathcal{T}_{\gamma}\mathcal{L}_*\mathcal{T}_{-\gamma}
\ee
where $\mathcal{L}_*$ is the linearization operator at $u_*$ defined in \eqref{linearizationRotat}. From \eqref{linearizationOperRotatGamma} it is clear that $+\mathrm{i}\omega_*,-\mathrm{i}\omega_*$ and $0$ are eigenvalues of $\mathcal{L}_{\gamma}$ with eigenfunctions 
\benn
\mathcal{T}_{\gamma}\left(\partial_{y_1}+ \mathrm{i}\partial_{y_2}\right)u_{*},\quad \mathcal{T}_{\gamma}\left(\partial_{y_1}-\mathrm{i}\partial_{y_2}\right)u_{*} \quad \text{ and } \quad 
\mathcal{T}_{\gamma}\partial_{\psi}u_*.
\eenn
Moreover, we notice that $\mathcal{L}_{\gamma}^{ad}$, the adjoint of $\mathcal{L}_{\gamma}$, is 
\begin{equation*}
\mathcal{L}^{ad}_{\gamma}=\mathcal{T}_{\gamma}\mathcal{L}_*^{ad}\mathcal{T}_{-\gamma}
\end{equation*}where $\mathcal{L}_*^{ad}$ is the adjoint of $\mathcal{L}_*$ and $+\mathrm{i}\omega_*,-\mathrm{i}\omega_*$ and $0$ are eigenvalues of $\mathcal{L}^{ad}_{\gamma}$ with eigenfunctions $\mathcal{T}_{\gamma}\phi^{+\mathrm{i}\omega_*},\mathcal{T}_{\gamma}\phi^{-\mathrm{i}\omega_*}$ and $\mathcal{T}_{\gamma}\phi^0$ where $\phi^{+\mathrm{i}\omega_*},\phi^{-\mathrm{i}\omega_*}$ and $\phi^{0}$ are the eigenfunctions of $\mathcal{L}_*^{ad}$ relative to the eigenvalues $+\mathrm{i}\omega_*,-\mathrm{i}\omega_*$ and $0$. It is important to notice that the spectrum of $\mathcal{L}_{\gamma}$ is the same as the spectrum of $\mathcal{L}_*$ as 
\begin{equation*}
\lambda I-\mathcal{L}_{\gamma}=\mathcal{T}_{\gamma}\left(\lambda I-\mathcal{L}_*\right)\mathcal{T}_{-\gamma}
\end{equation*}
and therefore also $\mathcal{L}_{\gamma}$ satisfies the spectral stability property (Assumption \ref{SpectAsyRot}). It follows that for every $\gamma\in SE(2)$ we can decompose $H^{2}\left(\mathbb{R}^2, \mathbb{R}^N\right)$ as in the previous section and we obtain the projections for $h\in H^{2}\left(\mathbb{R}^2, \mathbb{R}^N\right)$
\begin{equation}\label{ProjectRotatGamma}
P_{\gamma}h=\sum_{i=0}^2 \left\langle\mathcal{T}_{\gamma}\eta^{i} , h\right\rangle\mathcal{T}_{\gamma}\partial_iu_*=\mathcal{T}_{\gamma}P\mathcal{T}_{-\gamma}h
\end{equation}where $P$ is the projection \eqref{ProjectRotat} and $\eta^0,\eta^1$ and $\eta^2$ are the functions defined in \eqref{functionsEta}. Moreover as we have $\mathcal{T}_{\gamma}H^4_{eucl}\left(\mathbb{R}^2,\mathbb{R}^N\right)= H^4_{eucl}\left(\mathbb{R}^2,\mathbb{R}^N\right)$ from Proposition 3.43 in \cite{skriptBeyn_Otten} and therefore from the definition of generator it follows that the semigroup satisfies $\exp(\mathcal{L}_{\gamma}t)=\mathcal{T}_{\gamma}\exp(\mathcal{L}_{*}t)\mathcal{T}_{-\gamma}$. This yields that the strongly continuous semigroup $\exp(\mathcal{L}_{\gamma}t)$ generated by $\mathcal{L}_{\gamma}$ satisfies the bound
\begin{equation}\label{GrowthBoundGamma}
    \|\left.\exp(\mathcal{L}_{\gamma}t)\right|_{ R(I-P_{\gamma})}\|\leq C\exp(-\omega t)
\end{equation}
for any $0\leq \omega<b$ where $b>0$ is the spectral bound in \eqref{SpectAsyRot}.

\subsection{Deterministic Stability Results}

We recall now the PDE stability results for asymptotically constant rotating waves. Consider the following Cauchy problem
\begin{equation}\label{react_diff_corot_pert}
\begin{cases}
\partial_t u(t)=D\Delta u(t)+\omega_*\partial_{\psi} u(t)+f(u(t)), \quad  t\geq 0 \\
u_0=u_*+v_0
\end{cases}
\end{equation}
where $D\Delta+\omega_*\partial_{\psi}$ is an operator defined on $H^2(\mathbb{R}^2,\mathbb{R}^N)$ with domain $H^4_{eucl}(\mathbb{R}^2,\mathbb{R}^N)$ and where $u_*$ is an asymptotically constant rotating wave and $v_0\in H^2(\mathbb{R}^2,\mathbb{R}^N)$ is a small perturbation.
We have the following result:

\begin{theorem} \label{nonlStabRotWavesTHM}(Theorem 1.1 in \cite{nonlinearstab_Rotwaves})
Under Assumptions \ref{reg_f_rot}, \ref{regularity_rotWaves}, \ref{NegDefMat} and \ref{SpectAsyRot} there exist constants $\varepsilon>0$ and $C>0$ so that the following holds for $\|v_0\|_{H^2}<\varepsilon$ 
\begin{enumerate}
\item  The mild solution $u(y,t)$ of \ref{react_diff_corot_pert} exists for all $t\geq0$.
\item  The solution $u(y,t)$ can be written in the form \begin{equation}
u(y,t)=u_*\left(R_{-\theta(t)}\left(\bar{y}-\bar{b}(t)\right)\right)+z(x,t)
\end{equation} 
where $\theta \in C^1([0,+\infty),\mathbb{R}/(2\pi\mathbb{Z}))$ and $\bar{b}\in C^1([0,+\infty],\mathbb{R}^2)$.
\item We have the stability estimate
\begin{equation}
\|z(\cdot,t)\|_{H^2}\leq C\exp\left(-\frac{b}{2} t\right)\|v_0\|_{H^2}
\end{equation}
\item The initial conditions satisfy 
\begin{equation*}
|\theta(0)|+|\bar{b}(0)|\leq C\|v_0\|_{H^2}.
\end{equation*}
\item  There exists $\theta_{\infty}\in \mathbb{R}/(2\pi\mathbb{Z})$ and $\bar{b}_{\infty}\in \mathbb{R}^2$ depending on $v_0$, such that \begin{equation}
|\theta(t)-\theta_{\infty}|+|\bar{b}(t)-R_{-ct}\bar{b}_{\infty}|\leq C\exp(-\beta t)\|v_0\|_{H_2}.
\end{equation}
\end{enumerate}
\end{theorem}

In the proof of the previous theorem the identity  
\begin{equation}
\label{noncommutativityIdent}
\partial_{\psi}\mathcal{T}_{\gamma}w=\mathcal{T}_{\gamma}\partial_{\psi}w+\mathcal{T}_{\gamma} \left(\left(\nabla w\right)^T R_{-\theta+\frac{\pi}{2}}\begin{bmatrix}b_1 \\
 b_2\end{bmatrix}\right)
\end{equation}
is needed and will be very important also in the stochastic case. The identity \eqref{noncommutativityIdent} follows from the following computation
\begin{equation}
\label{computatDerivRot}
\begin{aligned}
\partial_{\psi}\mathcal{T}_{\gamma}w=&\partial_{\psi}w(R_{\theta}^{-1}(\bar{y}-b))
=
\left(y_1\partial_{y_2}-y_2 \partial_{y_1}\right)w(R_{-\theta}(\bar{y}-b))\\
=&\left(	\mathcal{T}_{\gamma}\nabla w(\bar{y})\right)^\top R_{-\theta+\frac{\pi}{2}} \bar{y}
\\=&
\left(	\mathcal{T}_{\gamma}\nabla w(\bar{y})\right)^\top R_{-\theta+\frac{\pi}{2}} \left(\bar{y}-\bar{b}\right)+\mathcal{T}_{\gamma} \left(\nabla w\right)^\top R_{-\theta+\frac{\pi}{2}}\begin{bmatrix}b_1 \\
 b_2\end{bmatrix}
\\
=&
\mathcal{T}_{\gamma}\partial_{\psi}w+\mathcal{T}_{\gamma} \left(\nabla w\right)^T R_{-\theta+\frac{\pi}{2}}\begin{bmatrix}b_1 \\
 b_2\end{bmatrix}\\
\end{aligned}.
\end{equation}
For our work below, we shall also need a certain rotational derivative in different reference frames. That is, we are going to need the following object
\begin{align*}
\mathcal{T}_{\gamma}\partial_{\psi}\mathcal{T}_{-\gamma} =\mathcal{T}_{\gamma}\left( y_1 \partial_{y_2}-y_2 \partial_{y_1}\right)\mathcal{T}_{-\gamma}
\end{align*}
Applying this to a smooth function $u$ and applying the chain rule we obtain
\begin{align}
\mathcal{T}_{\gamma}\partial_{\psi}\mathcal{T}_{-\gamma}u(\bar{y})&=\mathcal{T}_{\gamma}\left( y_1 \partial_{y_2}-y_2 \partial_{y_1}\right)u(R_{\theta}\bar{y}+\bar{b})\\=
& \mathcal{T}_{\gamma}\left( \left(	\mathcal{T}_{-\gamma}\nabla u(\bar{y})\right)^\top R_{\theta+\frac{\pi}{2}}\bar{y}\right)\\=
& \left(\nabla u(\bar{y})\right)^\top R_{\theta+\frac{\pi}{2}}R_{-\theta}\left(\bar{y}-\bar{b}\right)\\=
&\partial_{\psi}u-\left(\nabla u(\bar{y})\right)^\top R_{\frac{\pi}{2}}\bar{b}\\=
&\partial_{\psi}u-\left( b_1 \partial_{y_2}-b_2 \partial_{y_1}\right)
\end{align}This shows that 
\begin{align}\label{RotDerDiffRefFrame}
\mathcal{T}_{\gamma}\partial_{\psi}\mathcal{T}_{-\gamma}=\partial_{\psi}+\left( b_1 \partial_{y_2}-b_2 \partial_{y_1}\right)\\
=\left( y_1+b_1\right) \partial_{y_2}-\left(y_2+b_2\right) \partial_{y_1}
\end{align}Therefore, $\mathcal{T}_{\gamma}\partial_{\psi}\mathcal{T}_{-\gamma}$ is exactly the rotational derivative centered in $\left(-b_1,-b_2\right)$ as we expected as the rotation does not affect the rotational derivative.

As outlined in the introduction, we have now completed the relevant preparation to set up the notation and to recall the required stability results for rotating waves in the deterministic PDE case. Next, we proceed to the novel SPDE situation.

\section{Stochastic Rotating Waves}
\label{sec:stochRot}

In this section we want to study the dynamics of rotating waves in the stochastic setting. We, therefore, consider in the setting of Section~\ref{sec:rotwaves} the reaction-diffusion equation on the plane perturbed with small scale noise and an asymptotically constant rotating wave profile as initial condition. That is, we are considering the Cauchy problem 
\begin{equation}
\begin{cases}
d u_{t} =\left[D \Delta u_t+f\left(u_{t}\right)\right] d t+\varepsilon C\left(t,u_{t}\right) d W_{t}\\
u(\bar{x},0)=u_*(\bar{x})
\end{cases}
\end{equation}
where $u_*$ is an asymptotically constant rotating wave profile, $\varepsilon>0$ is a small scale parameter and $W_t$ is a (cylindrical) Wiener process (to simplify the notation we will assume, without loss of generality, that the covariance matrix is the identity) and we are in the setting of Chapter \ref{ch:introduction}. Changing to the co-rotating reference frame, we obtain the SPDE on the stochastic basis $\left(\Omega, \mathcal{F},\left(\mathcal{F}_{t}\right)_{t \in[0, T]}, \mathbb{P}\right)$ and on the Hilbert space $\mathcal{H}=H^2\left(\mathbb{R}^2,\mathbb{R}^N\right)$
\begin{equation}\label{stochReactDiffRotat}
\begin{cases}
d u_{t} =\left[D \Delta u_t+\omega_{*} \partial_{\psi}u_t+f\left(u_{t}\right)\right] d t+\varepsilon B\left(t,u_{t}\right) d W_{t}\\
u(\bar{y},0)=u_*(\bar{y})
\end{cases}
\end{equation} with the closed densely defined linear operator $A=D\Delta+\omega_*\partial_{\psi}$, defined on $H^2(\mathbb{R}^2,\mathbb{R}^N)$ with domain $H_{eucl}^4(\mathbb{R}^2,\mathbb{R}^N)$ where $B(t,z)w(\cdot)=C(t,z)w(R_{\omega_*t}(\cdot))$. Assuming
\begin{assumption}\label{Reg_stochRot}
$f$ and $B$ are like in Theorem \ref{ExistenceUniquenessSPDE} and, additionally, we require the uniform bound $B(t,x)<C_B$ for every $x\in H^2(\mathbb{R}^2,\mathbb{R}^N)$ and $t\geq 0$.
\end{assumption}
Theorem \ref{ExistenceUniquenessSPDE} guarantees the existence of a unique mild solution to \eqref{stochReactDiffRotat}. We now want to understand the dynamics of this solution. In order to do this, we will split the dynamics into an SDE identifying an element of a three dimensional manifold and a fluctuation term. We will do this in two different ways in the following sections. In the next sections we will identify a chart of $SE(2)$ with $\mathbb{R}^3$. This will not be a problem because the stochastic processes we will consider are continuous.

\subsection{Variational phase SDE}

We are looking for a splitting of the form
\begin{equation*}
u(\bar{y},t)=\mathcal{T}_{\beta_t}u_*(\bar{y})+w(\bar{y},t)
\end{equation*}
where with $\mathcal{T}_{\beta_t}$ we mean $\mathcal{T}_{\gamma(\beta^0_t,\beta^1_t,\beta^2_t)}$ and where $\mathcal{T}_{(\cdot)}$ is the operator defined in \eqref{IsometrySpeciaEucl} and $\gamma$ is the "chart" of $\mathrm{SE(2)}$ defined in \eqref{chartSE}. In particular, we want $w(\cdot,t)\in R(I-P_{\beta_t})$ where $R(I-P_{\beta_t})$ is the range of the projection $P_{\gamma}$ defined in \eqref{ProjectRotatGamma} with $\gamma=\beta_t$ at every time $t>0$. In fact, in this way we can consider the linearization operator at $\mathcal{T}_{\beta_{\tilde{t}}}u_*$ for a fixed time $\tilde{t}>0$, that is $\mathcal{L}_{\beta_{\tilde{t}}}$ (see \eqref{linearizationOperRotatGamma}).
From \eqref{GrowthBoundGamma} it follows that for the strongly continuous semigroup $\exp\left(\mathcal{L}_{\beta_{\tilde{t}}}t\right)$ generated by $\mathcal{L}_{\beta_{\tilde{t}}}$ we have
\begin{equation}
||\exp\left(\mathcal{L}_{\beta_{\tilde{t}}}t\right)w(\cdot,\tilde{t})||_{H^2}\leq C\exp(-\omega t)||w(\cdot,\tilde{t})||_{H^2}
\end{equation}for every $0<\omega<b$ where $b>0$ is the spectral bound in Section~\ref{SpectAsyRot}. To have $w(\cdot,t)\in R(I-P_{\beta_t})$, we have to impose 
\begin{equation*}
P_{\beta_t}w(\cdot,t)=0
\end{equation*}and this condition is equivalent (see \eqref{ProjectRotatGamma}) to requiring that $\left\langle w(\cdot,t),\mathcal{T}_{\beta_t}\eta^i\right\rangle=0$ for $i=0,1,2$. Therefore, we impose the condition 
\begin{equation}
\label{VariationalPhaseRotating}
0=\left\langle u(\cdot,t)-\mathcal{T}_{\beta_t}u_*,\mathcal{T}_{\beta_t}\eta^i\right\rangle
\end{equation}
for $i=0,1,2$. For this reason, we define the variational phase $\beta_{t}$ as the stochastic process satisfying the following conditions
\begin{enumerate}
\item  It is continuous for all $t<\tau$ where $\tau$ is the stopping time \begin{equation*}
\tau=\inf \left\{t \geq 0: \operatorname{det}\left(\mathcal{M}\left(u_{t}, \beta_{t}\right)\right)=0\right\}
\end{equation*}where the matrix $\mathcal{M}(z, \boldsymbol{\alpha})$ is defined in \eqref{matrixDerRotatingVar}.
\item  For all $t<\tau$, it solves the identity\begin{equation}\label{variationalCondition}
0=\mathcal{G}(u_t,\beta_t)
\end{equation}where $$
\begin{aligned}
\mathcal{G}: \left( H^2(\mathbb{R}^2,\mathbb{R}^N)\right)&\times \mathbb{R}^3  \longrightarrow\mathbb{R}^3  \\
\left(z,\gamma\right)\mapsto\mathcal{G}_i(z, \gamma) &:=\left\langle z-\mathcal{T}_{\gamma}u_*, \mathcal{T}_{\gamma}\eta^i\right\rangle
\end{aligned}
$$
for $i=0,1,2$ where $\eta^i$ are as defined in \eqref{functionsEta}.
\end{enumerate}
The first condition ensures, thanks to the implicit function theorem, that the identity $\mathcal{G}(u_t,\beta_t)$ can be uniquely inverted for $\beta_t$ in a small neighborhood. We emphasize that - except for highly contrived examples - in most stochastic systems one must expect that eventually any local definition of a phase will break down, if rare fluctuations force the system a long way from the manifold $\lbrace \mathcal{T}_{\gamma}u_* \rbrace$. We define also the $3 \times 3$ square matrix $\mathcal{M}(z, \boldsymbol{\alpha})$ with elements
\begin{equation}
\label{matrixDerRotatingVar}
\mathcal{M}_{i j}(z, \gamma)=-\partial_{\gamma_{j} }\mathcal{G}_{i}(z, \gamma)
\end{equation}
 Using the differentiation rules we can compute the elements of the matrix more explicitly 
$$
\begin{aligned}
\mathcal{M}_{ij}(z, \alpha)&=\left\langle \partial_{\gamma^j}\mathcal{T}_{\gamma}u_*, \mathcal{T}_{\gamma}\eta^{i}\right\rangle-\left\langle z-\mathcal{T}_{\gamma}u_*, \partial_{\gamma^j}\mathcal{T}_{\gamma}\eta^{i}\right\rangle\\
\end{aligned}
$$
We now define the matrix $\mathcal{M}^*(z,\gamma)$ with elements
\begin{equation}
\label{Mstar}
\mathcal{M}_{pq}^*(z,\gamma)=\left\langle z-\mathcal{T}_{\gamma}u_*, \mathcal{T}_{\gamma}\partial_p\eta^{q}\right\rangle
\end{equation}
where $\partial_0=\partial_{\psi}$, $\partial_1=\partial_{y_1}$ and $\partial_2=\partial_{y_2}$ and we notice 
\begin{equation}
\label{formMRot}
\mathcal{M}(z,\gamma)=\begin{bmatrix} 1 & 0\\
0 & -R_{-\gamma^0}
\end{bmatrix}-\begin{bmatrix} 1 & 0\\
0 & -R_{-\gamma^0}
\end{bmatrix}\mathcal{M}^*(z,\gamma)\\
=\begin{bmatrix} 1& 0\\
0 &-R_{-\gamma^0}
\end{bmatrix}\left(I-\mathcal{M}^*(z,\gamma)\right)
\end{equation}
From the implicit function theorem, we have that we can always solve \eqref{variationalCondition} for some neighborhood of $\left(u_{t}, \beta_{t}\right)$ as long as the matrix $\mathcal{M}\left(u_{t}, \beta_{t}\right)$ is invertible. Therefore, assuming that $t<\tau$ we always have a local solution for $\beta_{t}$ in terms of $u_{t}$. In fact, for $t=0$ we have $u_0=u_*$ and therefore $\beta_0=0$ is the solution to \eqref{variationalCondition} and
\begin{equation*}
\mathcal{M}(u_*,0)=\begin{bmatrix} 1&0\\
0 &-I
\end{bmatrix}
\end{equation*}
is invertible. We notice from \eqref{Mstar} using the Cauchy inequality that the elements of $\mathcal{M}^*\left(u_{t}, \beta_{t}\right)$ are small as long as $\|u_t-\mathcal{T}_{\beta_t}u_*\|$ is small enough. From \eqref{formMRot} we notice that $\mathcal{M}(z,\gamma)$ is always invertible as long as the eigenvalues of $\mathcal{M}^*(u_t,\beta_t)$ are small. Therefore, as long as $\|u_t-\mathcal{T}_{\beta_t}u_*\|$ is small enough we have that $\mathcal{M}(z,\gamma)$ is invertible. We want now to determine the dynamics of $\beta_{t}$  finding (at least formally) an $\mathrm{SDE}$ of the form
$$
d \beta_{t}=\mathcal{V}\left(u_{\mathrm{t}}, \beta_{t}\right) d t+\varepsilon \mathcal{Y}\left(t, u_{t}, \beta_{t}\right) d W_{t}
$$
for functions $\mathcal{V}: \mathbb{R}_{+}\times H^2(\mathbb{R}^2,\mathbb{R}^N) \times \mathbb{R}^3 \rightarrow \mathbb{R}^3$ and $\mathcal{Y}: \mathbb{R}_{+}\times  H^2(\mathbb{R}^2,\mathbb{R}^N)\times \mathbb{R}^3\rightarrow \mathcal{L}\left(H^2(\mathbb{R}^2,\mathbb{R}^N), \mathbb{R}^3\right)$ to be determined below. All of the following considerations of this section will be formal as It{\^o} formula is not directly available for unbounded operators in infinite dimension.

To satisfy \eqref{variationalCondition} we impose  $d \mathcal{G}\left(u_{t}, \beta_{t}\right)=0$. Applying formally Ito formula we obtain
\begin{equation}\label{ItoGRot}
d \mathcal{G}_{i}=\left\langle d u_{t}, \mathcal{T}_{\beta_t}\eta^{i}\right\rangle+\sum_{j=0}^{2} \partial_{\beta_{t}^{j}} \mathcal{G}_{i} d\beta_{t}^{j}+\frac{1}{2} \sum_{j, k=0}^{2}  \partial_{\beta_{t}^{j} \beta_{t}^{k} }\mathcal{G}_{i}d \beta_{t}^{j} d \beta_{t}^{k}+\sum_{j=0}^{2}\left\langle d u_{t},  \partial_{\beta_t^j}\mathcal{T}_{\beta_t}\eta^{i}\right\rangle d \beta_{i}^{j},
\end{equation}
Here the respective covariations of the processes are written as $d \beta_{t}^{j} d \beta_{t}^{k}$ and $d u_{t} d \beta_{t}^{j} .$ Using the fact that $D \Delta u_* +\omega_{*} \partial_{\psi}u_*+f\left(u_*\right)=0$ and the computation \eqref{computatDerivRot}, we get that 
$$
\begin{aligned}
 &D\Delta\mathcal{T}_{\beta_{t}}u_*+\omega_{*} \partial_{\psi}\mathcal{T}_{\beta_{t}}u_*+f\left(\mathcal{T}_{\beta_{t}}u_*\right)\\=&\mathcal{T}_{\beta_{t}} \left(\Delta u_* +\omega_{*} \partial_{\psi}u_*+f\left(u_*\right)\right)+\mathcal{T}_{\beta_{t}}\left( \left(\nabla u_*\right)^T R_{-\beta_t^0+\frac{\pi}{2}}\begin{bmatrix}\beta_t^1 \\
 \beta_t^2\end{bmatrix}\right)\\
 =&\mathcal{T}_{\beta_t} \left(\nabla u_*\right)^T R_{-\beta_t^0+\frac{\pi}{2}}\begin{bmatrix}\beta_t^1 \\
 \beta_t^2\end{bmatrix}\\
\end{aligned}
$$ 
where the last term means the application of the operator $\mathcal{T}_{\beta_t}$ to the function obtained by the matrix operations, see also \eqref{noncommutativityIdent}. Now if $u_{t}$ were in $ H_{eucl}^4(\mathbb{R}^2,\mathbb{R}^N) $, then, it would hold that
$$
\begin{aligned}
\left\langle d u_{t}, \mathcal{T}_{\beta_t}\eta^{i}\right\rangle &=\left\langle\left\{ D  \Delta u_t +\omega_{*} \partial_{\psi} u_{t}+f\left(u_{t}\right)\right\} d t+\varepsilon B\left(t, u_{t}\right) d W_{t}, \mathcal{T}_{\beta_t}\eta^{i}\right\rangle \\
&= \left\langle \left\{D \Delta u_t +\omega_{*} \partial_{\psi} u_{t}-D \Delta\mathcal{T}_{\beta_t}u_*-\omega_{*} \partial_{\psi} \mathcal{T}_{\beta_t}u_*+f\left(u_{t}\right)-f\left(\mathcal{T}_{\beta_t}u_*\right)\right. \right. \\
& \left. \left. +\mathcal{T}_{\beta_{t}} \left(\nabla u_*\right)^T R_{-\beta_t^0+\frac{\pi}{2}}\begin{bmatrix}\beta_t^1 \\
 \beta_t^2\end{bmatrix}\right\} d t+\varepsilon B\left(t, u_{t}\right) d W_{t}, \mathcal{T}_{\beta_t}\eta^{i}\right\rangle \\
&=\left\langle u_{t}-\mathcal{T}_{\beta_t}u_*, D \Delta \mathcal{T}_{\beta_t}\eta^{i}-\omega_{*}\partial_{\psi}\mathcal{T}_{\beta_t}\eta^{i}\right\rangle d t+\left\langle f\left(u_{t}\right)-f\left(\mathcal{T}_{\beta_t}u_*\right), \mathcal{T}_{\beta_t}\eta^{i}\right\rangle d t\\
&+\varepsilon\left\langle B\left(t, u_{t}\right) d W_{t}, \mathcal{T}_{\beta_t}\eta^{i}\right\rangle+\left\langle\mathcal{T}_{\beta_{t}}\left(\nabla u_*\right)^T R_{-\beta_t^0+\frac{\pi}{2}}\begin{bmatrix}\beta_t^1 \\
 \beta_t^2\end{bmatrix}, \mathcal{T}_{\beta_t}\eta^{i}\right\rangle dt.
\end{aligned}
$$
Obviously the previous computation is only formal in fact we do not have so much regularity for $u_t$. However, notice that if $u_{t}$ is not in $ H^4_{eucl}(\mathbb{R}^2,\mathbb{R}^N)$ the expression in the last line of the previous equation is still well-defined as we are assuming that $\eta^i\in H^4_{eucl}(\mathbb{R}^2,\mathbb{R}^N)$ for $i=0,1,2$.
Matching the stochastic terms (the coefficients of $d W_{t}$) in \eqref{ItoGRot}, we find that
\begin{equation*}
\varepsilon\left\langle B\left(t, u_{t}\right) d W_{t}, \mathcal{T}_{\beta_t}\eta^{i}\right\rangle-\varepsilon \sum_{j=0}^{2} \mathcal{M}_{i j}\left(u_{t}, \beta_t\right) \mathcal{Y}_{j}\left(t, u_{t}, \beta_t\right) d W_{t}=0
\end{equation*}
Inverting the previous equation, we obtain that the linear operator $\mathcal{Y}_{j}\left(t, u_{t}, \beta_t\right)$ must be characterized as follows for each $z \in H$,
$$
\begin{aligned}
\mathcal{Y}_{i}\left(t, u_{t}, \beta_{t}\right) \cdot z &=\sum_{j=0}^{2} \mathcal{N}_{i j}\left(u_{t}, \beta_t\right)\left\langle B\left(t, u_{t}\right) z, \mathcal{T}_{\beta_t}\eta^{j}\right\rangle \text { where } \\
\mathcal{N}\left(u_{t}, \beta_{t}\right) &=\mathcal{M}\left(u_{t}, \beta_{t}\right)^{-1} \text { and } \mathcal{N}\left(u_{t}, \beta_{t}\right)=\left(\mathcal{N}_{i j}\left(u_{t}, \beta_{t}\right)\right)_{1 \leq i, j \leq m}
\end{aligned}
$$
denoting as $\mathcal{M}\left(u_{t}, \beta_t\right)^{-1}$ the matrix inverse of $\mathcal{M}\left(u_{t}, \beta_t\right)$. The inverse matrix exists as  from the definition of the stopping time $\tau$  we have that $\mathcal{M}\left(u_{t}, \beta_t\right)$ is invertible for $t<\tau$. We thus obtain that the covariations are of the form (using standard theory for stochastic integrals with respect to infinite dimensional Wiener Processes, see Chapter 4.3 in \cite{daprato_zabczyk_2014})
$$
\begin{aligned}
&d \beta_{t}^{j} d \beta_{t}^{k} =\varepsilon^{2} \sum_{p, q=0}^{2} \mathcal{N}_{j p}\left(u_{t}, \beta_t\right) \mathcal{N}_{k q}\left(u_{t}, \beta_t\right)\left\langle B^{ad}\left(t, u_{t}\right) \mathcal{T}_{\beta_t}\eta^{p}, B^{ad}\left(t, u_{t}\right) \mathcal{T}_{\beta_t}\eta^{q}\right\rangle d t \\
&\left\langle d u_{t}, \partial_{\beta^j_t}\mathcal{T}_{\beta_t}\eta^{i}\right\rangle d \beta_{t}^{j} =\varepsilon^{2} \sum_{p=0}^{2} \mathcal{N}_{j p}\left(u_{t}, \beta_t\right)\left\langle B^{ad}\left(t, u_{t}\right) \mathcal{T}_{\beta_t}\eta^{p}, B^{ad}\left(t, u_{t}\right) \partial_{\beta^j_t}\mathcal{T}_{\beta_t}\eta^{i}\right\rangle d t,
\end{aligned}
$$
recalling that $B^{ad}(t,u_t)$ is the adjoint of $B$. We can now compare the drift terms (all of the coefficients of $d t$ terms) in \eqref{ItoGRot} to formally compute $\mathcal{V}\left(u_{t}, \beta_t\right)$. Observe that
$$
\begin{aligned}
& \partial_{ \gamma^{j}\gamma^{k}}\mathcal{G}_{i}=\\
&\left\langle z-\mathcal{T}_{\gamma}u_*, \partial_{\gamma^j\gamma^k}\mathcal{T}_{\gamma}\eta^i\right\rangle-\left\langle \partial_{\gamma^j}\mathcal{T}_{\gamma}u_*,  \partial_{\gamma^k}\mathcal{T}_{\gamma}\eta^i\right\rangle\\
&=\left\langle z-\mathcal{T}_{\gamma}u_*, \partial_{\gamma^j\gamma^k}\mathcal{T}_{\gamma}\eta^i\right\rangle+\left\langle\mathcal{T}_{\gamma}u_*,  \partial_{\gamma^j\gamma^k}\mathcal{T}_{\gamma}\eta^i\right\rangle \\
&=\left\langle z, \partial_{\gamma^j\gamma^k}\mathcal{T}_{\gamma}\eta^i\right\rangle
 \end{aligned}
$$
using the integration by parts formula \eqref{integByparts}. Thus we find that
$$
\begin{aligned}
&-\sum_{j=0}^{2} \mathcal{M}_{i j}\left(u_{t}, \beta_t\right) \mathcal{V}_{j}\left(u_{t}, \beta_t\right)+\\
&\varepsilon^{2} \sum_{j, p=0}^{2} \mathcal{N}_{j p}\left(u_{t}, \beta_t\right)\left\langle B^{ad}\left(t, u_{t}\right) \mathcal{T}_{\beta_t}\eta^{p}, B^{ad}\left(t, u_{t}\right)\partial_{\beta^j_t}\mathcal{T}_{\beta_t}\eta^{i}\right\rangle +\\
&\frac{\varepsilon^{2}}{2} \sum_{j, k, p, q=0}^{2}\left\langle u_{t}, \partial_{\beta_t^j\beta_t^k}\mathcal{T}_{\beta_t}\eta^i\right\rangle \mathcal{N}_{j p}\left(u_{t}, \beta_t\right) \mathcal{N}_{k q}\left(u_{t}, \beta_t\right)\left\langle B^{ad}\left(t, u_{t}\right) \mathcal{T}_{\beta_t}\eta^{p}, B^{ad}\left(t, u_{t}\right) \mathcal{T}_{\beta_t}\eta^{q}\right\rangle \\
&+\left\langle u_{t}-\mathcal{T}_{\beta_t}u_*, D\Delta\mathcal{T}_{\beta_t}\eta^{i}-\omega_* \partial_{\psi}\mathcal{T}_{\beta_t}\eta^{i}\right\rangle+\left\langle f\left(u_{t}\right)-f\left(\mathcal{T}_{\beta_t}u_*\right), \mathcal{T}_{\beta_t}\eta^{i}\right\rangle +\\
&\left\langle\mathcal{T}_{\beta_t} \left(\nabla u_*\right)^T R_{-\beta_t^0+\frac{\pi}{2}}\begin{bmatrix}\beta_t^1 \\
 \beta_t^2\end{bmatrix}, \mathcal{T}_{\beta_t}\eta^{i}\right\rangle=0 .
\end{aligned}
$$
We now invert the matrix $\mathcal{M}\left(u_{t}, \beta_t\right)$ and we find that for $0 \leq r \leq 2$,
$$
\begin{aligned}
&\mathcal{V}_{r}\left(t, u_{t}, \beta_t\right)=\\
&\sum_{i=0}^{2} \mathcal{N}_{r i}\left(u_{t}, \beta_t\right)\left\{\varepsilon^{2} \sum_{j, p=0}^{2} \mathcal{N}_{j p}\left(u_{t}, \beta_t\right)\left\langle B^{ad}\left(t, u_{t}\right) \mathcal{T}_{\beta_t}\eta^{p}, B^{ad}\left(t, u_{t}\right) \partial_{\beta^j_t}\mathcal{T}_{\beta_t}\eta^{i}\right\rangle\right.+ \\
&\frac{\varepsilon^{2}}{2} \sum_{j, k, p, q=0}^{2}\left\langle u_{t}, \partial_{\beta_t^j\beta_t^k}\mathcal{T}_{\beta_t}\eta^i\right\rangle \mathcal{N}_{j p}\left(u_{t}, \beta_{t}\right) \mathcal{N}_{k q}\left(u_{t}, \beta_t\right)\left\langle B^{ad}\left(t, u_{t}\right) \mathcal{T}_{\beta_t}\eta^{p}, B^{ad}\left(t, u_{t}\right) \mathcal{T}_{\beta_t}\eta^{q}\right\rangle \\
+&\left\langle u_{t}-\mathcal{T}_{\beta_t}u_*, D\Delta \mathcal{T}_{\beta_t}\eta^{i}-\omega_* \partial_{\psi} \mathcal{T}_{\beta_t}\eta^{i}\right\rangle+\left\langle f\left(u_{t}\right)-f\left(\mathcal{T}_{\beta_t}u_*\right),  \mathcal{T}_{\beta_t}\eta^{i}\right\rangle+\\
&\left\langle\mathcal{T}_{\beta_{t}} \left(\nabla u_*\right)^T R_{-\beta_t^0+\frac{\pi}{2}}\begin{bmatrix}\beta_t^1 \\
 \beta_t^2\end{bmatrix},  \mathcal{T}_{\beta_t}\eta^{i}\right\rangle\}
\end{aligned}
$$
Now defining $\sigma:\{0,1,2\}\rightarrow \{0,1,2\}$ such that $\sigma(0)=0, \sigma(1)=2$ and $\sigma(2)=0$ and $\lambda_0=0, \lambda_1=\omega_*$ and $\lambda_2=-\omega_*$ we have
$$
\begin{aligned}
&\left\langle u_{t}-\mathcal{T}_{\beta_t}u_*, D\Delta\mathcal{T}_{\beta_t}\eta^{i}-\omega_* \partial_{\psi}\mathcal{T}_{\beta_t}\eta^{i}\right\rangle+\left\langle f'\left(\mathcal{T}_{\beta_t}u_*\right) \cdot\left(u_{t}-\mathcal{T}_{\beta_t}u_*\right), \mathcal{T}_{\beta_t}\eta^{i}\right\rangle\\
=&\left\langle u_{t}-\mathcal{T}_{\beta_t}u_*, \mathcal{T}_{\beta_t}^{i}\left(D\Delta\eta^{i}-\omega_*\partial_{\psi}\eta^{i}+ f'(u_*)^T\eta^{i}\right)\right\rangle+\\
&\left\langle u_{t}-\mathcal{T}_{\beta_t}u_*, \mathcal{T}_{\beta_{t}} \left(\nabla\eta^i\right)^T R_{-\beta_t^0+\frac{\pi}{2}}\begin{bmatrix}\beta_t^1 \\
 \beta_t^2\end{bmatrix}\right\rangle\\
=&\lambda_i\left\langle u_{t}-\mathcal{T}_{\beta_t}u_*,\mathcal{T}_{\beta_t}\eta^{\sigma(i)}\right\rangle+\\
&\left\langle u_{t}-\mathcal{T}_{\beta_t}u_*, \mathcal{T}_{\beta_{t}} \left(\nabla\eta^i\right)^T R_{-\beta_t^0+\frac{\pi}{2}}\begin{bmatrix}\beta_t^1 \\
 \beta_t^2\end{bmatrix}\right\rangle ,
\end{aligned}
$$
since  $\mathcal{L}_{\beta_{t}}^{ad}\mathcal{T}_{\beta_t}\eta^i=\lambda_i\mathcal{T}_{\beta_t}\eta^{\sigma(i)}$, see \eqref{modifiedEigen}. We thus find that
$$
\begin{aligned}
&\mathcal{V}_{r}\left(t, u_{t}, \beta_t\right)=\\
&\sum_{i=0}^{2} \mathcal{N}_{r i}\left(u_{t}, \beta_t\right)\left\{\varepsilon^{2} \sum_{j, p=0}^{2} \mathcal{N}_{j p}\left(u_{t}, \beta_t\right)\left\langle B^{ad}\left(t, u_{t}\right) \mathcal{T}_{\beta_t}\eta^{p}, B^{ad}\left(t, u_{t}\right) \partial_{\beta^j_t}\mathcal{T}_{\beta_t}\eta^{i}\right\rangle+\right. \\
&\left.\frac{\varepsilon^{2}}{2} \sum_{j, k, p, q=0}^{2}\left\langle u_{t}, \partial_{\beta_t^j\beta_t^k}\mathcal{T}_{\beta_t}\eta^i\right\rangle \mathcal{N}_{j p}\left(u_{t}, \beta_t\right) \mathcal{N}_{k q}\left(u_{t}, \beta_t\right)\left\langle B^{ad}\left(t, u_{t}\right) \mathcal{T}_{\beta_t}\eta^{p}, B^{ad}\left(t, u_{t}\right) \mathcal{T}_{\beta_t}\eta^{q}\right\rangle \right.\\
&\left.+\left\langle f\left(u_{t}\right)-f\left(\mathcal{T}_{\beta_t}u_*\right)-f'\left(\mathcal{T}_{\beta_t}u_*\right) \left(u_{t}-\mathcal{T}_{\beta_t}u_*\right), \mathcal{T}_{\beta_t}\eta^{i}\right\rangle+\lambda_i\left\langle u_{t}-\mathcal{T}_{\beta_t}u_*,\mathcal{T}_{\beta_t}\eta^{\sigma(i)}\right\rangle \right.\\
&\left.+\left\langle u_{t}-\mathcal{T}_{\beta_t}u_*, \mathcal{T}_{\beta_{t}} \left(\nabla\eta^i\right)^T R_{-\beta_t^0+\frac{\pi}{2}}\begin{bmatrix}\beta_t^1 \\
 \beta_t^2\end{bmatrix}\right\rangle+\left\langle\mathcal{T}_{\beta_t} \left(\nabla u_*\right)^T R_{-\beta_t^0+\frac{\pi}{2}}\begin{bmatrix}\beta_t^1 \\
 \beta_t^2\end{bmatrix}, \mathcal{T}_{\beta_t}\eta^{i}\right\rangle\right\}
\end{aligned}
$$

Notice that if we have that $0=\mathcal{G}_{i}\left(u_{t}, \beta_{t}\right) =\left\langle z-\mathcal{T}_{\beta_t}u_*, \mathcal{T}_{\beta_t}\eta^{i}\right\rangle$ for every $i=0,1,2$ we obviously have $0=\lambda_i\left\langle u_{t}-\mathcal{T}_{\beta_t}u_*,\mathcal{T}_{\beta_t}\eta^{\sigma(i)}\right\rangle$ and therefore it follows that
$$
\begin{aligned}
&\mathcal{V}_{r}\left(t, u_{t}, \beta_t\right)=\\
&\sum_{i=0}^{2} \mathcal{N}_{r i}\left(u_{t}, \beta_t\right)\left\{\varepsilon^{2} \sum_{j, p=0}^{2} \mathcal{N}_{j p}\left(u_{t}, \beta_t\right)\left\langle B^{ad}\left(t, u_{t}\right) \mathcal{T}_{\beta_t}\eta^{p}, B^{ad}\left(t, u_{t}\right) \partial_{\beta^j_t}\mathcal{T}_{\beta_t}\eta^{i}\right\rangle+\right. \\
&\left.\frac{\varepsilon^{2}}{2} \sum_{j, k, p, q=0}^{2}\left\langle u_{t}, \partial_{\beta_t^j\beta_t^k}\mathcal{T}_{\beta_t}\eta^i\right\rangle \mathcal{N}_{j p}\left(u_{t}, \beta_t\right) \mathcal{N}_{k q}\left(u_{t}, \beta_t\right)\left\langle B^{ad}\left(t, u_{t}\right) \mathcal{T}_{\beta_t}\eta^{p}, B^{ad}\left(t, u_{t}\right) \mathcal{T}_{\beta_t}\eta^{q}\right\rangle \right.\\
&\left.+\left\langle f\left(u_{t}\right)-f\left(\mathcal{T}_{\beta_t}u_*\right)-f'\left(\mathcal{T}_{\beta_t}u_*\right) \left(u_{t}-\mathcal{T}_{\beta_t}u_*\right), \mathcal{T}_{\beta_t}\eta^{i}\right\rangle\right.\\
&\left.+\left\langle u_{t}-\mathcal{T}_{\beta_t}u_*, \mathcal{T}_{\beta_{t}} \left(\nabla\eta^i\right)^T R_{-\beta_t^0+\frac{\pi}{2}}\begin{bmatrix}\beta_t^1 \\
 \beta_t^2\end{bmatrix}\right\rangle+\left\langle\mathcal{T}_{\beta_t} \left(\nabla u_*\right)^T R_{-\beta_t^0+\frac{\pi}{2}}\begin{bmatrix}\beta_t^1 \\
 \beta_t^2\end{bmatrix}, \mathcal{T}_{\beta_t}\eta^{i}\right\rangle\right\}
\end{aligned}
$$
We now define the fluctuation term $w_t=u_t-\mathcal{T}_{\beta_t}u_*$ and we want to understand its dynamics. We observe this satisfies the equation
$$
\begin{aligned}
d w_{t}=&\left(D\Delta u_{t}+\omega_*\partial_{\psi} u_{t}+f\left(u_{t}\right)\right) d t+\varepsilon B\left(t, u_{t}\right) d W_{t}\\
&-\sum_{i=0}^{2} \partial_{\beta_t^i}\mathcal{T}_{\beta_t}u_* d \beta_t^{i}-\frac{1}{2} \sum_{j, k=0}^{2}\partial_{\beta_t^j \beta_t^k}\mathcal{T}_{\beta_t}u_* d \beta_{t}^{j} d \beta_{t}^{k} \\
=&D\Delta u_{t}+\omega_*\partial_{\psi} u_{t} +f\left(u_{t}\right)-\mathcal{T}_{\beta_t}\left(D\Delta u_*+\omega_*\partial_{\psi}u_*+f\left(u_*\right)\right) d t\\
&+\varepsilon B\left(t, u_{t}\right) d W_{t}-\sum_{i=0}^{2} \partial_{\beta_t^i}\mathcal{T}_{\beta_t}u_* d \beta_{t}^{i}-\frac{1}{2} \sum_{j, k=0}^{2}\partial_{\beta_t^j \beta_t^k}\mathcal{T}_{\beta_t}u_* d \beta_{t}^{j} d \beta_{t}^{k} \\
&\left(w_{t}+f^{\prime}\left(\mathcal{T}_{\beta_a}u_*\right) \cdot w_{t} +f\left(u_{t}\right)-f\left(\mathcal{T}_{\beta_t}u_*\right)\right. \\
&\left.-f^{\prime}\left(\mathcal{T}_{\beta_a}u_*\right) \cdot w_{t}+\mathcal{T}_{\beta_t} \left(\nabla u_*\right)^T R_{-\beta_t^0+\frac{\pi}{2}}\begin{bmatrix}\beta_t^1 \\
 \beta_t^2\end{bmatrix}\right) d t+\varepsilon B\left(t, u_{t}\right) d W_{t}\\
 &-\sum_{i=0}^{2} \partial_{\beta_t^i}\mathcal{T}_{\beta_t}u_* d \beta_{t}^{i} -\frac{1}{2} \sum_{j, k=0}^{2}\partial_{\beta_t^j \beta_t^k}\mathcal{T}_{\beta_t}u_* d \beta_{t}^{j} d \beta_{t}^{k}=\\
=&D\Delta u_{t}+\omega_*\partial_{\psi} u_{t} +f\left(u_{t}\right)\\
&-\left( D\Delta\mathcal{T}_{\beta_t}u_* +\omega_{*} \partial_{\psi}\mathcal{T}_{\beta_t}u_*+f\left(\mathcal{T}_{\beta_t}u_*\right)-\mathcal{T}_{\beta_t} \left(\nabla u_*\right)^T R_{-\beta_t^0+\frac{\pi}{2}}\begin{bmatrix}\beta_t^1 \\
 \beta_t^2\end{bmatrix}\right) d t\\
 &+\varepsilon B\left(t, u_{t}\right) d W_{t}-\sum_{i=0}^{2} \partial_{\beta_t^i}\mathcal{T}_{\beta_t}u_* d \beta_{t}^{i}-\frac{1}{2} \sum_{j, k=0}^{2}\partial_{\beta_t^j \beta_t^k}\mathcal{T}_{\beta_t}u_* d \beta_{t}^{j} d \beta_{t}^{k} \\
\end{aligned}
$$
where we have used the identity \eqref{noncommutativityIdent} again. Therefore, at least formally, we can do the following computation
$$
\begin{aligned}
d w_{t}=&\left(D\Delta w_{t}+\omega_*\partial_{\psi} w_{t}+f^{\prime}\left(\mathcal{T}_{\beta_a}u_*\right) \cdot w_{t} +f\left(u_{t}\right)-f\left(\mathcal{T}_{\beta_t}u_*\right)\right. \\
&\left.-f^{\prime}\left(\mathcal{T}_{\beta_a}u_*\right) \cdot w_{t}+\mathcal{T}_{\beta_t} \left(\nabla u_*\right)^T R_{-\beta_t^0+\frac{\pi}{2}}\begin{bmatrix}\beta_t^1 \\
 \beta_t^2\end{bmatrix}\right) d t+\varepsilon B\left(t, u_{t}\right) d W_{t}\\
 &-\sum_{i=0}^{2} \partial_{\beta_t^i}\mathcal{T}_{\beta_t}u_* d \beta_{t}^{i} -\frac{1}{2} \sum_{j, k=0}^{2}\partial_{\beta_t^j \beta_t^k}\mathcal{T}_{\beta_t}u_* d \beta_{t}^{j} d \beta_{t}^{k}=\\
&\left(D\Delta w_{t}+\omega_*\mathcal{T}_{\boldsymbol{\beta}_{a}}\partial_{\psi}\mathcal{T}_{-\boldsymbol{\beta}_{a}} w_{t}-\left( \beta_a^1 \partial_{y_2}-\beta_a^2 \partial_{y_1}\right)w_t+f^{\prime}\left(\mathcal{T}_{\beta_a}u_*\right) \cdot w_{t} \right. \\
&\left.+f\left(u_{t}\right)-f\left(\mathcal{T}_{\beta_t}u_*\right)-f^{\prime}\left(\mathcal{T}_{\beta_a}u_*\right) \cdot w_{t}+\mathcal{T}_{\beta_t} \left(\nabla u_*\right)^T R_{-\beta_t^0+\frac{\pi}{2}}\begin{bmatrix}\beta_t^1 \\
 \beta_t^2\end{bmatrix}\right) d t\\
 &+\varepsilon B\left(t, u_{t}\right) d W_{t}-\sum_{i=0}^{2} \partial_{\beta_t^i}\mathcal{T}_{\beta_t}u_* d \beta_{t}^{i} -\frac{1}{2} \sum_{j, k=0}^{2} \partial_{\beta_t^i\beta_t^j}\mathcal{T}_{\beta_t}u_* d \beta_{t}^{j} d \beta_{t}^{k}=\\
 &\left(\mathcal{L}_{\beta_a}w_t-\left( \beta_a^1 \partial_{y_2}-\beta_a^2 \partial_{y_1}\right)w_t+\right. \\
&\left.+f\left(u_{t}\right)-f\left(\mathcal{T}_{\beta_t}u_*\right)-f^{\prime}\left(\mathcal{T}_{\beta_a}u_*\right) \cdot w_{t}+\mathcal{T}_{\beta_t} \left(\nabla u_*\right)^T R_{-\beta_t^0+\frac{\pi}{2}}\begin{bmatrix}\beta_t^1 \\
 \beta_t^2\end{bmatrix}\right) d t\\
 &+\varepsilon B\left(t, u_{t}\right) d W_{t}-\sum_{i=0}^{2} \partial_{\beta_t^i}\mathcal{T}_{\beta_t}u_* d \beta_{t}^{i} -\frac{1}{2} \sum_{j, k=0}^{2} \partial_{\beta_t^i\beta_t^j}\mathcal{T}_{\beta_t}u_* d \beta_{t}^{j} d \beta_{t}^{k}\\
\end{aligned}
$$

where we used the identity $$\partial_{\psi}=\mathcal{T}_{\beta_{a}}\partial_{\psi}\mathcal{T}_{-\beta_{a}}-\left( \beta_a^1 \partial_{y_2}-\beta_a^2 \partial_{y_1}\right)$$ see \eqref{RotDerDiffRefFrame}. We remark that the previous equation is only formal as $w_t$ is not a smooth function. 

Therefore, we have
 $$
\begin{aligned}
d w_{t}=&
\left(D\Delta w_{t}+\omega_*\mathcal{T}_{\boldsymbol{\beta}_{a}}\partial_{\psi}\mathcal{T}_{-\boldsymbol{\beta}_{a}} w_{t}+f^{\prime}\left(\mathcal{T}_{\beta_a}u_*\right) \cdot w_{t}-\left( \beta_a^1 \partial_{y_2}-\beta_a^2 \partial_{y_1}\right)w_t\right.\\
& \left. +f\left(u_{t}\right)-f\left(\mathcal{T}_{\beta_t}u_*\right)-f^{\prime}\left(\mathcal{T}_{\beta_a}u_*\right) \cdot w_{t}+\mathcal{T}_{\beta_t} \left(\nabla u_*\right)^T R_{-\beta_t^0+\frac{\pi}{2}}\begin{bmatrix}\beta_t^1 \\
 \beta_t^2\end{bmatrix}\right) d t\\
 &+\varepsilon B\left(t, u_{t}\right) d W_{t}-\sum_{i=0}^{2} \partial_{\beta_t^i}\mathcal{T}_{\beta_t}u_* \mathcal{Y}_{i}\left(t, u_{t}, \beta_t\right) dW_t-\sum_{i=0}^{2} \partial_{\beta_t^i}\mathcal{T}_{\beta_t}u_* \mathcal{V}_{i}\left(t, u_{t}, \beta_t\right)dt \\
 &-\frac{\varepsilon^{2}}{2} \sum_{j, k, p, q=0}^{2}\partial_{\beta_t^j \beta_t^k}\mathcal{T}_{\beta_t}u_* \mathcal{N}_{j p}\left(u_{t}, \beta_t\right) \mathcal{N}_{k q}\left(u_{t}, \beta_t\right)\left\langle B^{ad}\left(t, u_{t}\right) \psi_{\beta_t}^{p}, B^{ad}\left(t, u_{t}\right) \psi_{\beta_t}^{q}\right\rangle
\end{aligned}
$$

We observe that the linear correction term $$-\left( \beta_a^1 \partial_{y_2}-\beta_a^2 \partial_{y_1}\right)$$ is not a stochastic term and does not depend on the scale parameter $\varepsilon$. For this reason, it can be cancelled from $\sum_{i=0}^{2} \partial_{\beta_t^i}\mathcal{T}_{\beta_t}u_* \mathcal{V}_{i}\left(t, u_{t}, \beta_t\right)dt$ only.
 However, this is in general not possible, since $\partial_{\beta_t^i}\mathcal{T}_{\beta_t}u_* $ is an element of $span(\mathcal{T}_{\beta_t}\partial_0 u_*,\mathcal{T}_{\beta_t}\partial_1 u_*,\mathcal{T}_{\beta_t}\partial_2 u_*)$ but $w_t$  lives in $R(I-P_{\beta_t})$  by construction. 
For completeness, we notice that this correction term is a first order operator with bounded coefficients.

Therefore, instead of $\mathcal{L}_{\beta_a}$, we have a new linear operator. The linear correction term emerges in fact from the "manifold" structure of the special Euclidean group $\mathrm{SE(2)}$. To our best knowledge, this is the first time that this linear correction term is pointed out. This is an important starting observation for an understanding of stochastic pattern with more complex manifold-like underlying geometry. In fact, this correction term does not appear for the more studied stochastic patterns with a simpler vector-like underlying geometry, like for example stochastic travelling waves (see \cite{maclaurin2020metastability}) where the isometries commute with the differential operators.

Of course, this additional correction term makes it more challenging applying the variational phase algebraically to study the dynamics of stochastic rotating waves. Therefore, in this work, we also develop an approximated variational phase as a possible alternative, which is algebraically easier but is going to yield less strong stability results. In fact, the variational phase characterizes more precisely, not only locally, the dynamics of the stochastic rotating wave in comparison to the approximated variational phase. 

\subsection{Rigorous Proof of the Variational Phase SDE}
Up to this point, our calculations have only been formal. We now rigorously prove that the phase satisfies the above SDE. There does not exist a general version of It\^{o}'s Lemma for SPDEs with unbounded operators. However, although our system is driven by the diffusion operator (which is unbounded) we are saved by the fact that the manifold of shifted / rotated solutions that we are projecting onto lies in the domain of the operator. Thus when performing the change of variable, one can always take the adjoint of the diffusion operator, and avoid the problem of the diffusion operator causing some terms to blowup. Our approach, following \cite{inglis2016general,maclaurin2020metastability}, which in turn is an adaptation of classical finite-dimensional proofs of It\^{o}'s Lemma \cite{karatzas1998brownian}, is to (i) discretize time (ii) take a second-order Taylor expansion over each time-step, (iii) move the unbounded operators onto the smooth manifold and (iv) take the width of the time intervals to zero.

We start by defining $\beta_t$ to be the solution of the SDE
\begin{equation}\label{eq: beta t pathwise}
d\beta_t = \mathcal{V}(t,\beta_t,u_t) + \varepsilon \mathcal{Y}(t,\beta_t,u_t) dW_t.
\end{equation}
Since the coefficients $\mathcal{V}$ are locally Lipschitz for all times before the stopping time $\tau$, one easily checks that there exists a unique strong solution to the above SDE for all times less than $\tau$. The `strong' property of the solution implies that the identity \eqref{eq: beta t pathwise} holds for almost-every realization of the Brownian Motions $W_t$ (see the discussion in \cite[Section 5]{karatzas1998brownian} of strong and weak solutions of finite-dimensional SDEs). This is a stronger result than merely identifying the `statistics' (i.e. the probability law) of the process, because it ensures that the phase SDE $\beta_t$ and the original SPDE $u_t$ are anchored in the same probability space. In particular, it allows us to directly compare solutions of \eqref{eq: beta t pathwise} with the original spatially-extended system $u_t$. 

 Let $\Pi = (t_i)_{i=1}^M$ be a partition of $[0,t]$ for some $t\geq0$. It is assumed that, as $M\to\infty$,
 \begin{equation}
 \sup_{1\leq k \leq M-1} \big| t_{k+1} - t_k \big| \to 0. \label{eq: Partition to zero}
 \end{equation}
 We then take the mild solution for $u_t$ over each time interval, i.e. for $t\in [t_k , t_{k+1}]$, and writing $A = D \Delta  +\omega_*\partial_{\psi}$ for the linear operator driving the dynamics, 
\begin{equation}\label{eq: SDE v t}
u_t = \exp\big(A(t-t_k) \big)u_{t_k} + \int_{t_k}^t \exp\big(A(t-s) \big)f(u_s) ds  
+ \varepsilon \int_{t_k}^t  \exp\big(A(t-s) \big) B(s,u_s,\beta_s)dW_s.
\end{equation}
One easily checks that the above solution is consistent with the mild solution in \eqref{stochMildsol}. Next, we introduce an additional stopping time that is necessary for us to control the size of terms in the coming Taylor expansion. We thus define, for $n \gg 1$,
\begin{multline}
\xi_n := \inf \bigg\lbrace t\in[0, \tau]: \det(\mathcal{M}(u_{t},\beta_t)) = n^{-1}\text{ or } \bigg\| \int_0^t B(s,u_s)dW_s \bigg\| \geq n\\ \text{ or }\sup_{1\leq i \leq m} \|\int_0^t \mathcal{Y}_i(s,u_s,\beta_s)dW_s\| \geq n \text{ or }\sup_{1\leq i \leq m}  |\beta^i_t | \geq n \bigg\rbrace.
\end{multline}
 It may be seen that $(\xi_n)_{n\geq1}$ is nondecreasing, and that $\lim_{n\to\infty}\xi_n = \tau$ a.s. Define for any $t\geq0$ $\beta^{n}_t = \beta_{t\wedge \xi_n}$ and $u^{n}_t = u_{t\wedge \xi_n}$. Although the It\^{o} formula in \eqref{ItoGRot} is not necessarily valid (one cannot \textit{a priori} assume an It\^{o} formula for unbounded operators), we may perform a second order Taylor expansion over the time interval $[t_k ,t_{k+1}]$ (which is valid). To this end, for some family $\lbrace \theta_k\rbrace_{k=1}^{M-1} \subset [0,1]$, set $w_k = \theta_k u^{n}_{t_{k}} + (1-\theta_k) u^{n}_{t_{k+1}}$ and $\beta^n_{t_k} = \theta_k \beta^{n}_{t_{k}} + (1-\theta_k) \beta^{n}_{t_{k+1}}$. A second-order Taylor Expansion of $\mathcal{G}_i$ implies that
\begin{multline}
\label{taylor}
\mathcal{G}_i(u^{n}_{t_k},\beta^{n}_{t_k}) -\mathcal{G}_i(u_0,\beta_0) = \sum_{k=1}^{M-1}\bigg\lbrace \sum_{j=0}^2 \big( -\mathcal{M}_{ij}(u^n_{t^k}, \beta_{t^k}) \big) (\beta^{n,j}_{t_{k+1}} - \beta^{n,j}_{t_k})+ \langle u^{n}_{t_{k+1}} -u^{n}_{t_k},\mathcal{T}_{\beta^n_k}\eta_i\rangle  \\
 +\sum_{p=0}^2\bigg(  \frac{1}{2}\sum_{q=0}^2\frac{\partial^2 \mathcal{G}_i}{\partial \beta_t^p\partial \beta_t^q}\bigg|_{\beta^n_{t_k},u^n_{t_k}} (\beta^{n,p}_{t_{k+1}} - \beta^{n,p}_{t_k})  (\beta^{n,q}_{t_{k+1}} - \beta^{n,q}_{t_k}) + \langle u^{n}_{t_{k+1}}-u^{n}_{t_k} ,\mathcal{T}_{\beta^n_{t_k}}\eta_i \rangle (\beta^{n,p}_{t_{k+1}} - \beta^{n,p}_{t_k})  \bigg)\bigg\rbrace ,
\end{multline}
for some choice of the parameters $\lbrace \theta_k\rbrace_{k=1}^{M-1} \subset [0,1]$. We now study the convergence of each of the terms as the partition size shrinks to zero. Our goal is to prove that the formal Ito identity \eqref{ItoGRot} is in fact valid.

\textit{Step 1: First-Order Terms}

\vspace{0.25 cm}
Since $\mathcal{M}_{ij}(\cdot,\cdot)$ is locally-Lipschitz for all times $t \leq \xi_n$, one obtains that as $M\to \infty$,
\begin{equation}
 \sum_{k=1}^{M-1}  \sum_{j=0}^2 \big( -\mathcal{M}_{ij}(u^n_{t^k}, \beta_{t^k}) \big) (\beta^{n,j}_{t_{k+1}} - \beta^{n,j}_{t_k})  \to  - \sum_{j=0}^2 \int_0^t \mathcal{M}_{ij}(u^n_{s}, \beta^n_{s}) d\beta^{n,j}_s .
\end{equation}
Substituting in the mild solution \eqref{eq: SDE v t} over each time interval for $u(t)$, we obtain that
\begin{multline}
 \sum_{k=1}^{M-1} \langle u^{n}_{t_{k+1}} -u^{n}_{t_k},\mathcal{T}_{\beta^{n}_{t_k}}\eta^i\rangle 
 = \sum_{k=1}^{M-1}\bigg\lbrace \big\langle  \exp\big((t_{k+1}-t_k)A \big)u_{t_k} - u_{t_k}, \mathcal{T}_{\beta^{n}_{t_k}}\eta^i \big\rangle \\+ \int_{t_k}^{t_{k+1}} \big\langle  \exp\big((t_{k+1}-s)A \big) f(u_{s}), \mathcal{T}_{\beta^{n}_{t_k}}\eta^i \big\rangle ds \\
  + \int_{t_k}^{t_{k+1}} \big\langle  \exp\big((t_{k+1}-s) A \big) B_s dW_s, \mathcal{T}_{\beta^{n}_{t_k}}\eta^i \big\rangle ds \bigg\rbrace .\label{eq: SDE v t}
\end{multline}
Now writing the adjoint of $A$ as $A^{ad}$,
\begin{align}
\sum_{k=1}^{M-1}\big\langle  \exp\big((t_{k+1}-t_k)A \big)u_{t_k} - u_{t_k}, \mathcal{T}_{\beta^{n}_{t_k}}\eta^i \big\rangle = \sum_{k=1}^{M-1} \big\langle  u_{t_k}, \exp\big((t_{k+1}-t_k)A^{ad} \big)\mathcal{T}_{\beta^{n}_{t_k}}\eta^i - \mathcal{T}_{\beta^{n}_{t_k}}\eta^i\big\rangle .
\end{align}
By assumption, $ \mathcal{T}_{\beta^{n}_{t_k}}\eta^i$ is in the domain of $A^{ad}$. This means that, if $t_k,t_{k+1} \to v \in [0,\tau)$ as $M \to \infty$, then
\begin{equation}
(t_{k+1} - t_k)^{-1}\big\lbrace \exp\big((t_{k+1}-t_k)A^{ad} \big)\mathcal{T}_{\beta^{n}_{t_k}}\eta^i - \mathcal{T}_{\beta^{n}_{t_k}}\eta^i \big\rbrace
\to A^{ad}\mathcal{T}_{\beta^{n}_{v}}\eta^i .
\end{equation}
We thus obtain from the dominated convergence theorem that, as the size of the partition asymptotes to zero,
\begin{align}
\sum_{k=1}^{M-1}\big\langle  \exp\big((t_{k+1}-t_k)A \big)u_{t_k} - u_{t_k}, \mathcal{T}_{\beta^{n}_{t_k}}\eta^i \big\rangle \to \int_0^t \big\langle u_s , A^{ad}  \mathcal{T}_{\beta_s}\eta^i \big\rangle ds .
\end{align}
Since $\exp\big((t_{k+1}-s)A \big) \to \mathbf{I}$ as $|t_{k+1} - s| \to 0$, we obtain from the dominated convergence theorem that
\begin{align}
\sum_{k=1}^{M-1}\int_{t_k}^{t_{k+1}} \big\langle  \exp\big((t_{k+1}-s)A \big) f(u_{s}), \mathcal{T}_{\beta^{n}_{t_k}}\eta^i \big\rangle ds \to \int_0^t  \big\langle   f(u_{s}), \mathcal{T}_{\beta^{n}_{s}}\eta^i \big\rangle ds .
\end{align}
Similarly, regularity properties of the stochastic integral imply that as the partition size goes to zero 
\begin{align}
\sum_{k=1}^{M-1}\int_{t_k}^{t_{k+1}} \big\langle  \exp\big((t_{k+1}-s) A \big) B_s dW_s, \mathcal{T}_{\beta^{n}_{t_k}}\eta^i \big\rangle ds \to \int_0^t \big\langle \mathcal{T}_{\beta^{n}_{s}}\eta^i , B_s dW_s \big\rangle.
\end{align}

\textit{Step 2: Second Order Terms}

\vspace{0.25 cm}

Substituting the mild solution for $u$ in \eqref{eq: SDE v t}, we obtain that
\begin{multline}
 \sum_{k=1}^{M-1} \big\langle u^{n}_{t_{k+1}}-u^{n}_{t_k} ,\mathcal{T}_{\beta^n_{t_k}}\eta_i \big\rangle (\beta^{n,p}_{t_{k+1}} - \beta^{n,p}_{t_k})    =
  \sum_{k=1}^{M-1} \bigg\langle\exp\big((t_{k+1}-t_k)A \big)u^n_{t_k} - u^n_{t_k} \\+ \int_{t_k}^{t_{k+1}}  \exp\big((t_{k+1}-s)A \big) f(u^n_{s})ds + \exp\big((t_{k+1}-s) A \big) B_s dW_s ,\mathcal{T}_{\beta^n_{t_k}}\eta_i \bigg\rangle (\beta^{n,p}_{t_{k+1}} - \beta^{n,p}_{t_k})=\\
  \sum_{k=1}^{M-1}\bigg\lbrace \big\langle u^n_{t_k}, \exp\big((t_{k+1}-t_k)A^{ad} \big)\mathcal{T}_{\beta^n_{t_k}}\eta_i-\mathcal{T}_{\beta^n_{t_k}}\eta_i \big\rangle + \int_{t_k}^{t_{k+1}} \big\langle f(u^n_{s}),  \exp\big((t_{k+1}-s)A^{ad} \big) \mathcal{T}_{\beta^n_{t_k}} \big\rangle ds \\+ \int_{t_k}^{t_{k+1}}\bigg\langle \exp\big((t_{k+1}-s) A \big) B_s dW_s ,\mathcal{T}_{\beta^n_{t_k}}\eta_i \bigg\rangle \bigg\rbrace (\beta^{n,p}_{t_{k+1}} - \beta^{n,p}_{t_k})
\end{multline}
Now  employing the fact that $\mathcal{T}_{\beta^n_{t_k}}\eta_i$ is in the domain of $A^{ad}$, it must be that there exists a constant $C$ such that for all choices of the partition and all $k$,
\begin{equation} 
 \sup_{M \in \mathbb{Z}^+} \sup_{1\leq k \leq M-1}(t_{k+1} - t_k)^{-1}\| \exp\big((t_{k+1}-t_k)A^{ad} \big)\mathcal{T}_{\beta^n_{t_k}}\eta_i - \mathcal{T}_{\beta^n_{t_k}}\eta_i\| \leq C,
 \end{equation}
 $\mathbb{P}$-almost-surely. 
Since the coefficients in the SDE for $\beta^n_t$ are Lipschitz and bounded (noting the definition of the stopping time $\xi_n$), standard properties of stochastic integrals \cite{karatzas1998brownian} imply that for any $r \in (0, 1/2)$,
\begin{equation}\label{eq: bound beta fluctuations}
\lim_{M\to\infty}\sup_{1\leq k \leq M-1} (t_{k+1}-t_k)^{-r}  \big| \beta^{n,p}_{t_{k+1}} - \beta^{n,p}_{t_k} \big| = 0.
\end{equation}
The above two identities, together with \eqref{eq: Partition to zero}, imply that $\mathbb{P}$-almost-surely,
\begin{align}
\lim_{M\to\infty} \bigg| \sum_{k=1}^{M-1}\big\langle u_{t_k}, \exp\big((t_{k+1}-t_k)A^{ad} \big)\mathcal{T}_{\beta^n_{t_k}}\eta_i - \mathcal{T}_{\beta^n_{t_k}}\eta_i \big\rangle (\beta^{n,p}_{t_{k+1}} - \beta^{n,p}_{t_k}) \bigg| = 0.
\end{align}
It also follows from \eqref{eq: bound beta fluctuations} (and the boundedness of $\|f(u_s)\|$ thanks to the definition of the stopping time) that
\begin{equation}
  \lim_{M\to\infty} \sum_{k=1}^{M-1} \int_{t_k}^{t_{k+1}} \big\langle f(u_{s}),  \exp\big((t_{k+1}-s)A^{ad} \big) \mathcal{T}_{\beta^n_{t_k}} \big\rangle ds (\beta^{n,p}_{t_{k+1}} - \beta^{n,p}_{t_k}) = 0.
\end{equation}
The remaining terms converge to the cross-variation, i.e.
\begin{align}
 \sum_{k=1}^{M-1}\int_{t_k}^{t_{k+1}}&\bigg\langle \exp\big((t_{k+1}-s) A \big) B_s dW_s ,\mathcal{T}_{\beta^n_{t_k}}\eta_i \bigg\rangle (\beta^{n,p}_{t_{k+1}} - \beta^{n,p}_{t_k}) \\
 =& \sum_{k=1}^{M-1}\int_{t_k}^{t_{k+1}}\bigg\langle  B_s dW_s ,\exp\big((t_{k+1}-s) A^{ad} \big)\mathcal{T}_{\beta^n_{t_k}}\eta_i \bigg\rangle (\beta^{n,p}_{t_{k+1}} - \beta^{n,p}_{t_k})\\
  \to& \int_0^t \big\langle  du^n_s , \mathcal{T}_{\beta_s}\eta_i \big\rangle d\beta_s,
\end{align}
using the fact that 
\[
\exp\big((t_{k+1}-s) A^{ad} \big)\mathcal{T}_{\beta^n_{t_k}}\eta_i\to \mathbf{I},
\]
as $|t_{k+1} -s | \to 0$. The final cross-variation term can be similarly shown to converge as $M\to\infty$, i.e.
\begin{equation}
  \frac{1}{2} \sum_{k=1}^{M-1}\frac{\partial^2 \mathcal{G}_i}{\partial \beta_t^p\partial \beta_t^q}\bigg|_{\beta^n_{t_k},u^n_{t_k}} (\beta^{n,p}_{t_{k+1}} - \beta^{n,p}_{t_k})  (\beta^{n,q}_{t_{k+1}} - \beta^{n,q}_{t_k}) \to \int_0^t\frac{\partial^2 \mathcal{G}_i}{\partial \beta_t^p\partial \beta_t^q}\bigg|_{\beta_s^n,u^n_s}  d\beta^{n,p}_s d\beta^{n,q}_s.
\end{equation}

\textit{Step 3: Combining these results}

Combining the results from Steps 1 and 2, we have demonstrated that the formal Ito Lemma in \eqref{ItoGRot} is in fact rigorously justified. The proof worked because the unbounded operator $A$ could always be moved to the smooth manifold parameterized by $\big\lbrace \mathcal{T}_{\gamma}u_* \big\rbrace_{\gamma}$. Now the phase $\beta_t$ was constructed to be precisely such that the formal expansion $d\mathcal{G}_i(t)$ in \eqref{ItoGRot} is zero. Since this formal expansion corresponds to the rigorous methods of Steps 1 and 2, it must necessarily be the case that $\mathcal{G}_i(t) = 0$ for all times less than $\xi_n$. We can then take $n\to\infty$, so that $\xi_n \to \tau$ (if $\tau$ is finite).

\subsection{Approximated variational phase SDE}

Having considered the variational phase, the moving/adapted frame leads to considerable algebraic complexity in the calculations. Yet, a more localized-in-time tracking is possible in a fixed frame, the approximated variational phase. We consider again a splitting of the form
\begin{equation*}
u(\bar{y},t)=\mathcal{T}_{\beta_t}u_*(\bar{y})+w(\bar{y},t)
\end{equation*}
Differently from before, this time we require $w(\cdot,t)\in R(I-P)$ where $R(I-P)$ is the range of the projection $P$ defined in \eqref{ProjectRotat}. In fact, in this way we can consider at a fixed time $\tilde{t}>0$ the linearization operator at $u_*$ that is $\mathcal{L}_{*}$ (see \eqref{linearizationRotat}). From Section~\ref{expBoundAsyRot} follows that for the strongly continuous semigroup $\exp\left(\mathcal{L}_{*}t\right)$ generated by $\mathcal{L}_{*}$ we have
\begin{equation}
||\exp\left(\mathcal{L}_{*}t\right)w(\cdot,\tilde{t})||_{H^2}\leq C\exp(-\omega t)||w(\cdot,\tilde{t})||_{H^2}
\end{equation}
where $0<\omega<b$ for the spectral bound $b>0$ in Section~\ref{SpectAsyRot}. To get $w(\cdot,t)\in R(I-P)$ we have to impose \begin{equation*}
Pw(\cdot,t)=0
\end{equation*}
and this condition is equivalent, see \eqref{ProjectRotat}, to require $\left\langle w(\cdot,t),\eta^i\right\rangle=0$ for $i=0,1,2$. Therefore, we impose the condition 
\begin{equation}
\label{VariationalPhaseRotatingApprox}
0=\left\langle u(\cdot,t)-\mathcal{T}_{\beta_t}u_*,\eta^i\right\rangle
\end{equation}
for $i=0,1,2$. For this reason, we define the approximated variational phase $\beta_{t}$ as the stochastic process satisfying the following conditions
\begin{enumerate}
\item  It is continuous for all $t<\tau$ where $\tau$ is the stopping time \begin{equation*}
\tau=\inf \left\{t \geq 0: \operatorname{det}\left(\mathcal{M}\left(u_{t}, \beta_{t}\right)\right)=0\right\}
\end{equation*}where the matrix $\mathcal{M}(z, \boldsymbol{\alpha})$ is defined in \eqref{matrixDerRotatingVarApprox}.
\item  For all $t<\tau$, it solves the identity\begin{equation}\label{variationalConditionApprox}
0=\mathcal{G}(u_t,\beta_t)
\end{equation}
where 

$$
\begin{aligned}
\mathcal{G}: H^2(\mathbb{R}^2,\mathbb{R}^N)&\times \mathbb{R}^3  \longrightarrow\mathbb{R}^3  \\
\left(z,\gamma\right)\mapsto\mathcal{G}_i(z, \gamma) &=\left\langle z-\mathcal{T}_{\gamma}u_*, \mathcal{T}_{\gamma}\eta^i\right\rangle
\end{aligned}
$$
for $i=0,1,2$ where $\eta^i$ are as defined in \eqref{functionsEta}.
\end{enumerate}

We define also the $3 \times 3$ square matrix $\mathcal{M}(z, \boldsymbol{\alpha})$ with elements
\begin{equation}
\label{matrixDerRotatingVarApprox}
\mathcal{M}_{i j}(z, \gamma)=-\partial_{\gamma_{j} }\mathcal{G}_{i}(z, \gamma)
\end{equation}
 Using the differentiation rules, we can compute the elements of the matrix more explicitly 
$$
\begin{aligned}
\mathcal{M}_{ij}(z, \alpha)&=\left\langle \partial_{\gamma^j}\mathcal{T}_{\gamma}u_*, \eta^{i}\right\rangle
\end{aligned}
$$
We now define the matrix $\mathcal{M}^*(z,\gamma)$ with elements
\begin{equation}\label{MstarApprox}
\mathcal{M}_{pq}^*(z,\gamma)=\left\langle \mathcal{T}_{\gamma}\partial_p u_*, \eta^{q}\right\rangle
\end{equation}
and we notice that we have
\begin{equation}\label{formMRotApprox}
\mathcal{M}(z,\gamma)=\begin{bmatrix} 1 & 0\\
0 & -R_{-\gamma^0} 
\end{bmatrix}\mathcal{M}^*(z,\gamma)
\end{equation}
From the implicit function theorem we have that we can always solve \eqref{variationalConditionApprox} for some neighborhood of $\left(u_{t}, \beta_{t}\right)$ as long as the matrix $\mathcal{M}\left(u_{t}, \beta_{t}\right)$ is invertible. Therefore, assuming that $t<\tau$ we always have a local solution for $\beta_{t}$ in terms of $u_{t}$. In fact, for $t=0$ we have $u_0=u_*$ and therefore $\beta_0=0$ is the solution to \eqref{variationalCondition} and
\begin{equation*}
\mathcal{M}(u_*,0)=\begin{bmatrix} 1  & 0\\
0 &- I
\end{bmatrix}
\end{equation*}
is invertible. We notice that with this approximated variational phase we cannot control the invertibility of $\mathcal{M}\left(u_{t}, \beta_{t}\right)$ through $\|u_t-\mathcal{T}_{\beta_t}u_*\|$. We want now to understand the dynamics of $\beta_{t}$  determining (at least formally for the moment) an $\mathrm{SDE}$ of the form
\begin{equation}
\label{approxVarSDERotating}
d \beta_{t}=\mathcal{V}\left(u_{\mathrm{t}}, \beta_{t}\right) d t+\varepsilon \mathcal{Y}\left(t, u_{t}, \beta_{t}\right) d W_{t}
\end{equation}
for functions $\mathcal{V}: \mathbb{R}_{+}\times H^2(\mathbb{R}^2,\mathbb{R}^N) \times \mathbb{R}^3 \rightarrow \mathbb{R}^3$ and $\mathcal{Y}: \mathbb{R}_{+}\times  H^2(\mathbb{R}^2,\mathbb{R}^N)\times \mathbb{R}^3\rightarrow \mathcal{L}\left(H^2(\mathbb{R}^2,\mathbb{R}^N), \mathbb{R}^3\right)$ to be determined below. We will now use formal computations to determine the terms of this equation and in a second moment we will show that the resulting equation has a solution that really satisfies \eqref{VariationalPhaseRotatingApprox}. To satisfy \eqref{variationalCondition} we impose  $d \mathcal{G}\left(u_{t}, \beta_{t}\right)=0$. Applying formally It\^{o}'s formula we obtain
\begin{equation}
\label{ItoGRotApprox}
d \mathcal{G}_{i}=\left\langle d u_{t}, \eta^{i}\right\rangle+\sum_{j=0}^{2} \partial_{\beta_{t}^{j}} \mathcal{G}_{i} d\beta_{t}^{j}+\frac{1}{2} \sum_{j, k=0}^{2}  \partial_{\beta_{t}^{j} \beta_{t}^{k} }\mathcal{G}_{i}d \beta_{t}^{j} d \beta_{t}^{k}
\end{equation}
where the respective covariations of the processes are denoted as $d \beta_{t}^{j} d \beta_{t}^{k}$ and $d u_{t} d \beta_{t}^{j} .$ If we have $u_{t}$ in $ H^4_{eucl}(\mathbb{R}^2,\mathbb{R}^N)$, using the fact that $D \Delta u_* +\omega_{*} \partial_{\psi}u_*+f\left(u_*\right)=0$  we have 
$$
\begin{aligned}
\left\langle d u_{t}, \eta^{i}\right\rangle &=\left\langle\left\{ D  \Delta u_t +\omega_{*} \partial_{\psi} u_{t}+f\left(u_{t}\right)\right\} d t+\varepsilon B\left(t, u_{t}\right) d W_{t},\eta^{i}\right\rangle \\
&= \left\langle \left\{D \Delta u_t +\omega_{*} \partial_{\psi} u_{t}-D \Delta\mathcal{T}_{\beta_t}u_*-\omega_{*} \partial_{\psi} \mathcal{T}_{\beta_t}u_*+f\left(u_{t}\right)-f\left(\mathcal{T}_{\beta_t}u_*\right)\right.\right.\\
&\left.\left.+\mathcal{T}_{\beta_{t}} \left(\nabla u_*\right)^T R_{-\beta_t^0+\frac{\pi}{2}}\begin{bmatrix}\beta_t^1 \\
 \beta_t^2\end{bmatrix}\right\} d t+\varepsilon B\left(t, u_{t}\right) d W_{t},\eta^{i}\right\rangle \\
&=\left\langle u_{t}-\mathcal{T}_{\beta_t}u_*, D \Delta \eta^{i}-\omega_{*}\partial_{\psi}\eta^{i}\right\rangle d t+\left\langle f\left(u_{t}\right)-f\left(\mathcal{T}_{\beta_t}u_*\right), \eta^{i}\right\rangle d t\\
&+\varepsilon\left\langle B\left(t, u_{t}\right) d W_{t}, \eta^{i}\right\rangle
+\left\langle\mathcal{T}_{\beta_{t}}\left(\nabla u_*\right)^T R_{-\beta_t^0+\frac{\pi}{2}}\begin{bmatrix}\beta_t^1 \\
 \beta_t^2\end{bmatrix}, \eta^{i}\right\rangle dt.
\end{aligned}
$$Notice that if $u_{t}$ is not in $ H^4_{eucl}(\mathbb{R}^2,\mathbb{R}^N)$ the expression in the last line of the previous equation is still well-defined as we are assuming that $\eta^i\in H^4_{eucl}(\mathbb{R}^2,\mathbb{R}^N)$ for $i=0,1,2$. Matching the stochastic terms (the coefficients of $\left.d W_{t}\right)$ in \eqref{ItoGRot}, we obtain that
$$
\varepsilon\left\langle B\left(t, u_{t}\right) d W_{t}, \eta^{i}\right\rangle-\varepsilon \sum_{j=0}^2 \mathcal{M}_{i j}\left(u_{t}, \beta_t\right) \mathcal{Y}_{j}\left(t, u_{t}, \beta_t\right) d W_{t}=0
$$
Inverting the previous equation, we determine that the linear operator $\mathcal{Y}_{j}\left(t, u_{t}, \beta_t\right)$ must be such that for each $z \in H^2(\mathbb{R}^2,\mathbb{R}^N)$,
$$
\begin{aligned}
\mathcal{Y}_{i}\left(t, u_{t}, \beta_{t}\right) \cdot z &=\sum_{j=0}^2 \mathcal{N}_{i j}\left(u_{t}, \beta_t\right)\left\langle B\left(t, u_{t}\right) z, \eta^{j}\right\rangle \text { where } \\
\mathcal{N}\left(u_{t}, \beta_{t}\right) &=\mathcal{M}\left(u_{t}, \beta_{t}\right)^{-1} \text { and } \mathcal{N}\left(u_{t}, \beta_{t}\right)=\left(\mathcal{N}_{i j}\left(u_{t}, \beta_{t}\right)\right)_{1 \leq i, j \leq m}
\end{aligned}
$$
denoting as $\mathcal{M}\left(u_{t}, \beta_t\right)^{-1}$ the matrix inverse of $\mathcal{M}\left(u_{t}, \beta_t\right)$. The inverse matrix exists as  from the definition of the stopping time $\tau$  we have that $\mathcal{M}\left(u_{t}, \beta_t\right)$ is invertible for $t<\tau$. We thus obtain that the covariations are of the form (using the standard theory for stochastic integrals with respect to infinite dimensional Wiener Processes, see Chapter 4.3 in \cite{daprato_zabczyk_2014}),
$$
\begin{aligned}
d \beta_{t}^{j} d \beta_{t}^{k} &=\varepsilon^{2} \sum_{p, q=0}^2 \mathcal{N}_{j p}\left(u_{t}, \beta_t\right) \mathcal{N}_{k q}\left(u_{t}, \beta_t\right)\left\langle B^{ad}\left(t, u_{t}\right)\eta^{p}, B^{ad}\left(t, u_{t}\right) \eta^{q}\right\rangle d t \\
\end{aligned}
$$
We can now compare the drift terms (all of the coefficients of $d t$ terms) in \eqref{ItoGRotApprox} to formally compute $\mathcal{V}\left(u_{t}, \beta_t\right)$ and we find that
$$
\begin{aligned}
&-\sum_{j=0}^2 \mathcal{M}_{i j}\left(u_{t}, \beta_t\right) \mathcal{V}_{j}\left(u_{t}, \beta_t\right)\\
&+\frac{\varepsilon^{2}}{2} \sum_{j, k, p, q=0}^2\left\langle \partial_{\beta_t^j\beta_t^k}\mathcal{T}_{\beta_t}u_*, \eta^i\right\rangle \mathcal{N}_{j p}\left(u_{t}, \beta_t\right) \mathcal{N}_{k q}\left(u_{t}, \beta_t\right)\left\langle B^{ad}\left(t, u_{t}\right) \eta^{p}, B^{ad}\left(t, u_{t}\right) \eta^{q}\right\rangle \\
&+\left\langle u_{t}-\mathcal{T}_{\beta_t}u_*, D\Delta\eta^{i}-\omega_* \partial_{\psi}\eta^{i}\right\rangle\\
&+\left\langle f\left(u_{t}\right)-f\left(\mathcal{T}_{\beta_t}u_*\right), \eta^{i}\right\rangle+\left\langle\mathcal{T}_{\beta_t} \left(\nabla u_*\right)^T R_{-\beta_t^0+\frac{\pi}{2}}\begin{bmatrix}\beta_t^1 \\
 \beta_t^2\end{bmatrix},\eta^{i}\right\rangle=0 .
\end{aligned}
$$
Inverting the matrix $\mathcal{M}\left(u_{t}, \beta_t\right)$, we are able to  find that for $0 \leq r \leq 2$,
$$
\begin{aligned}
&\mathcal{V}_{r}\left(t, u_{t}, \beta_t\right)=\sum_{i=0}^2 \mathcal{N}_{r i}\left(u_{t}, \beta_t\right)\left\{\right. \\
&+\frac{\varepsilon^{2}}{2} \sum_{j, k, p, q=0}^2\left\langle \partial_{\beta_t^j\beta_t^k}\mathcal{T}_{\beta_t}u_*, \eta^i\right\rangle\mathcal{N}_{j p}\left(u_{t}, \beta_{t}\right) \mathcal{N}_{k q}\left(u_{t}, \beta_t\right)\left\langle B^{ad}\left(t, u_{t}\right) \eta^{p}, B^{ad}\left(t, u_{t}\right) \eta^{q}\right\rangle \\
&+\left\langle u_{t}-\mathcal{T}_{\beta_t}u_*, D\Delta \eta^{i}-\omega_* \partial_{\psi} \eta^{i}\right\rangle+\left\langle f\left(u_{t}\right)-f\left(\mathcal{T}_{\beta_t}u_*\right), \eta^{i}\right\rangle\\
&+\left\langle\mathcal{T}_{\beta_{t}} \left(\nabla u_*\right)^T R_{-\beta_t^0+\frac{\pi}{2}}\begin{bmatrix}\beta_t^1 \\
 \beta_t^2\end{bmatrix},  \eta^{i}\right\rangle\}
\end{aligned}
$$
Now defining $\sigma:\{0,1,2\}\rightarrow \{0,1,2\}$ such that $\sigma(0)=0, \sigma(1)=2$ and $\sigma(2)=0$ and $\lambda_0=0, \lambda_1=\omega_*$ and $\lambda_2=-\omega_*$ we have
$$
\begin{aligned}
&\left\langle u_{t}-\mathcal{T}_{\beta_t}u_*, D\Delta\eta^{i}-\omega_* \partial_{\psi}\eta^{i}\right\rangle+\left\langle f'\left(u_*\right) \cdot\left(u_{t}-\mathcal{T}_{\beta_t}u_*\right),\eta^{i}\right\rangle\\
=&\left\langle u_{t}-\mathcal{T}_{\beta_t}u_*, \left(D\Delta\eta^{i}-\omega_*\partial_{\psi}\eta^{i}+ f'(u_*)^T\eta^{i}\right)\right\rangle\\
=&\lambda_i\left\langle u_{t}-\mathcal{T}_{\beta_t}u_*,\eta^{\sigma(i)}\right\rangle 
\end{aligned}
$$ since  $\mathcal{L}_{*}^{ad}\eta^i=\lambda_i\eta^{\sigma(i)}$ see \eqref{modifiedEigen}. We thus find that
$$
\begin{aligned}
&\mathcal{V}_{r}\left(t, u_{t}, \beta_t\right)=\sum_{i=0}^2 \mathcal{N}_{r i}\left(u_{t}, \beta_t\right)\left\{\right. \\
&+\frac{\varepsilon^{2}}{2} \sum_{j, k, p, q=0}^2\left\langle \partial_{\beta_t^j\beta_t^k}\mathcal{T}_{\beta_t}u_*, \eta^i\right\rangle\mathcal{N}_{j p}\left(u_{t}, \beta_{t}\right) \mathcal{N}_{k q}\left(u_{t}, \beta_t\right)\left\langle B^{ad}\left(t, u_{t}\right) \eta^{p}, B^{ad}\left(t, u_{t}\right) \eta^{q}\right\rangle \\
&+\lambda_i\left\langle u_{t}-\mathcal{T}_{\beta_{t}}u_*, \eta^{\sigma(i)}\right\rangle+\left\langle f\left(u_{t}\right)-f\left(\mathcal{T}_{\beta_t}u_*\right)-f^{\prime}\left(u_*\right)\left(u_t-\mathcal{T}_{\beta_t}u_*\right), \eta^{i}\right\rangle\\
&+\left\langle\mathcal{T}_{\beta_{t}} \left(\nabla u_*\right)^T R_{-\beta_t^0+\frac{\pi}{2}}\begin{bmatrix}\beta_t^1 \\
 \beta_t^2\end{bmatrix},  \eta^{i}\right\rangle\}
\end{aligned}
$$

Notice that having that $0=\mathcal{G}_{i}\left(u_{t}, \beta_{t}\right) =\left\langle u_t-\mathcal{T}_{\beta_t}u_*, \eta^{i}\right\rangle$ for every $i=0,1,2$ we obviously have $0=\lambda_i\left\langle u_{t}-\mathcal{T}_{\beta_t}u_*,\eta^{\sigma(i)}\right\rangle$ and therefore it follows that

$$
\begin{aligned}
&\mathcal{V}_{r}\left(t, u_{t}, \beta_t\right)=\sum_{i=0}^{2} \mathcal{N}_{r i}\left(u_{t}, \beta_t\right)\left\{\right. \\
&+\frac{\varepsilon^{2}}{2} \sum_{j, k, p, q=0}^{2}\left\langle \partial_{\beta_t^j\beta_t^k}\mathcal{T}_{\beta_t}u_*, \eta^i\right\rangle\mathcal{N}_{j p}\left(u_{t}, \beta_{t}\right) \mathcal{N}_{k q}\left(u_{t}, \beta_t\right)\left\langle B^{ad}\left(t, u_{t}\right) \eta^{p}, B^{ad}\left(t, u_{t}\right) \eta^{q}\right\rangle \\
&+\left\langle f\left(u_{t}\right)-f\left(\mathcal{T}_{\beta_t}u_*\right)-f^{\prime}\left(u_*\right)\left(u_t-\mathcal{T}_{\beta_t}u_*\right), \eta^{i}\right\rangle+\left\langle\mathcal{T}_{\beta_{t}} \left(\nabla u_*\right)^T R_{-\beta_t^0+\frac{\pi}{2}}\begin{bmatrix}\beta_t^1 \\
 \beta_t^2\end{bmatrix},  \eta^{i}\right\rangle\}
\end{aligned}
$$Moreover, we notice that 
\begin{equation*}
\begin{aligned}
&\sum_{r=0}^{2} \partial_{\beta_t^r}\mathcal{T}_{\beta_t}u_*\sum_{i=0}^{2} \mathcal{N}_{r i}\left(u_{t}, \beta_t\right)\left\langle\mathcal{T}_{\beta_{t}} \left(\nabla u_*\right)^T R_{-\beta_t^0+\frac{\pi}{2}}\begin{bmatrix}\beta_t^1 \\
 \beta_t^2\end{bmatrix},  \eta^{i}\right\rangle=\\
 &\mathcal{T}_{\beta_t} \left(\nabla u_*\right)^T R_{-\beta^0_t+\frac{\pi}{2}}\begin{bmatrix}\beta_t^1 \\
\beta_t^2 \\
\end{bmatrix}
\end{aligned}
\end{equation*}

Therefore we define 
$$
\begin{aligned}
&\widetilde{\mathcal{V}}_{r}\left(t, u_{t}, \beta_t\right)=\sum_{i=0}^{2} \mathcal{N}_{r i}\left(u_{t}, \beta_t\right)\left\{\right. \\
&+\frac{\varepsilon^{2}}{2} \sum_{j, k, p, q=0}^{2}\left\langle \partial_{\beta_t^j\beta_t^k}\mathcal{T}_{\beta_t}u_*, \eta^i\right\rangle\mathcal{N}_{j p}\left(u_{t}, \beta_{t}\right) \mathcal{N}_{k q}\left(u_{t}, \beta_t\right)\left\langle B^{ad}\left(t, u_{t}\right) \eta^{p}, B^{ad}\left(t, u_{t}\right) \eta^{q}\right\rangle \\
&+\left\langle f\left(u_{t}\right)-f\left(\mathcal{T}_{\beta_t}u_*\right)-f^{\prime}\left(u_*\right)\left(u_t-\mathcal{T}_{\beta_t}u_*\right), \eta^{i}\right\rangle\}
\end{aligned}
$$
and now we have for $w_t=u_t-\mathcal{T}_{\beta_t}u_*$

$$
\begin{aligned}
dw_t=&
\left[\mathcal{L}_*w_t+f\left(u_{t}\right)-f\left(\mathcal{T}_{\beta_t}u_*\right)-f^{\prime}(u_*)w_t\right] d t+\varepsilon B\left(t,u_{t}\right) d W_{t}\\
 &-\sum_{i=0}^{2} \partial_{\beta_t^i}\mathcal{T}_{\beta_t}u_* \mathcal{Y}_{i}\left(t, u_{t}, \beta_t\right) dW_t-\sum_{i=0}^{2}\partial_{\beta_t^i}\mathcal{T}_{\beta_t}u_* \tilde{\mathcal{V}}_{i}\left(t, u_{t}, \beta_t\right)dt \\
 &-\frac{\varepsilon^{2}}{2} \sum_{j, k=0}^{2} \partial_{\beta_t^j \beta_t^k}\mathcal{T}_{\beta_t}u_* d\beta^j_t d\beta^k_t
\end{aligned}
$$

\subsection{Rigorous Proof of the Approximated Variational Phase SDE}

We notice that the approximated variational phase SDE has locally Lipschitz coefficients $\mathcal{Y}\left(t, u_{t}, \beta_{t}\right)$ and $\mathcal{V}\left(t, u_{t}, \beta_{t}\right)$  in $\beta_t$ for a neighbourhood of $0$ and $\mathcal{Y}\left(t, u_{t}, \beta_{t}\right)$ is Hilbert-Schmidt. We therefore have a solution of the SDE. We can show that the solution satisfies \eqref{VariationalPhaseRotatingApprox}. We have, in fact,
$$
\begin{aligned}
dw_t=&\left[\mathcal{L}_*w_t+f\left(u_{t}\right)-f\left(\mathcal{T}_{\beta_t}u_*\right)-f^{\prime}(u_*)w_t+\mathcal{T}_{\beta_t} \left(\nabla u_*\right)^T R_{-\beta^0_t+\frac{\pi}{2}}\begin{bmatrix}\beta_t^1 \\
\beta_t^2 \\
\end{bmatrix}\right] d t\\
&+\varepsilon B\left(t,u_{t}\right) d W_{t}
-\sum_{i=0}^2\partial_{\beta_t^i}\mathcal{T}_{\beta_t}u_*d\beta_t^i
-\sum_{i,j=0}^2\partial_{\beta_t^i\beta_t^j}\mathcal{T}_{\beta_t}u_*d\beta_t^i d\beta_t^j\\
=&\left[\mathcal{L}_*w_t+f\left(u_{t}\right)-f\left(\mathcal{T}_{\beta_t}u_*\right)-f^{\prime}(u_*)w_t+\mathcal{T}_{\beta_t} \left(\nabla u_*\right)^T R_{-\beta^0_t+\frac{\pi}{2}}\begin{bmatrix}\beta_t^1 \\
\beta_t^2 \\
\end{bmatrix}\right] d t\\
&+\varepsilon B\left(t,u_{t}\right) d W_{t}
 -\sum_{i=0}^2 \partial_{\beta_t^i}\mathcal{T}_{\beta_t}u_* \mathcal{Y}_{i}\left(t, u_{t}, \beta_t\right) dW_t\\
 &-\sum_{i=0}^{2}  \partial_{\beta_t^i}\mathcal{T}_{\beta_t}u_* \mathcal{V}_{i}\left(t, u_{t}, \beta_t\right)dt -\frac{\varepsilon^{2}}{2} \sum_{j, k=0}^{2}  \partial_{\beta_t^j\beta_t^k}\mathcal{T}_{\beta_t}u_*d\beta_t^j d\beta_t^k\\
 =&\left(\mathcal{L}_{\beta_t} w_{t}+\mathcal{K}_{t}\right) d t+\varepsilon \widetilde{B}\left(t, u_{t}, \beta_t\right) d W_{t}
\end{aligned}
$$

where
$$
\begin{aligned}
\widetilde{B}(s, z, \boldsymbol{\alpha}): & \mathbb{R}^{+} \times H \times \mathbb{R}^{3} \rightarrow \mathcal{L}(H, H) \\
\widetilde{B}(s, z, \boldsymbol{\alpha})=& B(s, z)-\sum_{j=0}^{2} \partial_{\alpha^j}\mathcal{T}_{\alpha}u_* \mathcal{Y}_{j}(s, z, \boldsymbol{\alpha}) \\
&\mathcal{K}_{t}:=f\left(u_{t}\right)-f\left(\mathcal{T}_{\beta_t}u_*\right)- f^{\prime}\left(u_*\right) \cdot w_{t}-\sum_{i=0}^{2} \partial_{\beta_t^i}\mathcal{T}_{\beta_t}u_* \mathcal{V}_{i}\left(t, u_{t}, \beta_t\right) -\\
&\frac{\varepsilon^{2}}{2} \sum_{j, k, p, q=0}^{2}\partial_{\beta_t^j \beta_t^k}\mathcal{T}_{\beta_t}u_* \mathcal{N}_{j p}\left(u_{t}, \beta_t\right) \mathcal{N}_{k q}\left(u_{t}, \beta_t\right)\left\langle B^{ad}\left(t, u_{t}\right) \eta^{p}, B^{ad}\left(t, u_{t}\right) \eta^{q}\right\rangle
\end{aligned}
$$The solution for $w_{t}$, written in mild form, satisfies for $t \in\left[t_{k}, t_{k+1}\right]$,
$$
w_{t}=\exp(\mathcal{L}_*t) w_{t_{k}}+\int_{t_{k}}^{t}\exp(\mathcal{L}_*(t-s)) \mathcal{K}_{s} d s+\varepsilon \int_{t_{k}}^{t} \exp(\mathcal{L}_*(t-s))\widetilde{B}\left(s, u_{s}, \boldsymbol{\beta}_{s}\right) d W_{s},
$$
where $\exp(\mathcal{L}_*t)$ is the semigroup generated by $\mathcal{L}_{*} .$  We now define the stopping times
$$
\begin{array}{r}
\xi_{n}:=\inf \left\{t \in[0, \tau]: \operatorname{det}\left(\mathcal{M}\left(u_{t}, \beta_t\right)\right)=n^{-1} \text { or }\left\|\int_{0}^{t} B\left(s, u_{s}\right) d W_{s}\right\| \geq n\right. \\
\text { or } \left.\sup _{0\leq i \leq 2}\left\|\int_{0}^{t} \mathcal{Y}_{i}\left(s, u_{s}, \boldsymbol{\beta}_{s}\right) d W_{s}\right\| \geq n \text { or } \sup _{0\leq i \leq 2}\left|\beta_{t}^{i}\right| \geq n\right\}
\end{array}
$$
It may be seen that the sequence $\left(\xi_{n}\right)_{n \geq 1}$ is nondecreasing, and that $\lim _{n \rightarrow \infty} \xi_{n}=\tau$ a.s.
for some partition $\Pi=(t_k)_{k=1}^M$ of $[0,T]$. Now for any $r \geq 0$, it must be that
$$
\begin{aligned}
&\lim _{s \rightarrow 0} \sup _{r \leq z \leq t \leq r+s}(t-z)^{-1}\left\langle\exp(\mathcal{L}_*(t-z))  w_{z}, \eta^i\right\rangle =\\
&\lim _{s \rightarrow 0} \sup _{r \leq z \leq t \leq r+s}(t-z)^{-1}\left\langle w_{z}, \exp(\mathcal{L}_*(t-z)) ^{ad} \eta^i\right\rangle \\
& \rightarrow\left\langle\mathcal{L}_{*}^{ad} \eta^{i}, w_{r}\right\rangle
\end{aligned}
$$
and, therefore,
$$
\begin{aligned}
\sum_{k=0}^{M-1}\left\langle w_{t_{k+1}\wedge \xi_{n}}\right.&\left.-w_{t_{k}\wedge\xi_{n}}, \eta^i\right\rangle \\
& \rightarrow \int_{0}^{\xi_{n} \wedge t}\left[\left\langle\mathcal{L}_{*}^{ad} \eta^i, w_{s}\right\rangle+\left\langle\mathcal{K}_{s}, \eta^i\right\rangle\right] d s+\varepsilon \int_{0}^{\xi_{n} \wedge t}\left\langle\eta^i, \widetilde{B}(s) d W_{s}\right\rangle \\
& \rightarrow \int_{0}^{\xi_{n} \wedge t}\left[\lambda_i\left\langle w_s, \eta^{\sigma(i)}\right\rangle+\left\langle\mathcal{K}_{s}, \eta^i\right\rangle\right] d s+\varepsilon \int_{0}^{\xi_{n} \wedge t}\left\langle\eta^i, \widetilde{B}(s) d W_{s}\right\rangle =0
\end{aligned}
$$

\subsection{Stability results - Variational Phase}

In this section we use the variational phase to determine exponential concentration results. More precisely, we show that the probability
of the solution of the SPDE leaving a close neighborhood of the manifold of rotating solutions over an exponentially long period of time (in $\varepsilon^{-2}$ the magnitude of the noise diffusivity), is exponentially unlikely. Our basic plan follows that of \cite{maclaurin2020metastability}: we discretize time, and linearize the dynamics about the deterministic rotating wave at the start of the time period. Since the exponential decay of the semigroup generated by the linear operator may not be felt immediately (i.e. the constant $C$ may be much greater than $1$), we make sure to choose the time interval to be sufficiently long that the effects of this constant $C$ is dominated by the exponential decay $\exp(-bt)$.  Because it is necessary that we exploit the exponential decay of the linearization, it is essential that we linearize about the nearest rotating pattern. Therefore it is most natural to take the mild solution of the SPDE relative to the frame that is co-rotating with the nearest pattern, at the start of the time interval. In this section, we will additionally assume
\begin{assumption}
\label{integrabilityAssStochConv}
For any $t>s$ and any $x \in \mathcal{H}, \exp((t-s) \mathcal{L}_*) B(s, x)$ is a Hilbert-Schmidt operator, with Hilbert-Schmidt norm upperbounded by
$$
\int_{0}^{T_{0}}(t-s)^{-2 \xi} \sup _{x \in H}\left\|\exp((t-s) \mathcal{L_*}) B(s, x)\right\|_{H S}^{2} d s<C_{H S}
$$
for some $\xi \in(0,1 / 2)$ and constant $C_{H S}<\infty$, where
$$
T_{0}=\frac{\log( 4 C)}{b}
$$\end{assumption}
First of all, we define the stopping time $\tau_1=\inf\{t\geq 0 : ||w_t||=\kappa\}$ for a generic $\kappa>0$. We will add some requirements on $\kappa$ later. 

The following theorem gives the scaling of the first hitting time for leaving the manifold of rotating solutions: it parallels analogous bounds for stochastic oscillators \cite{giacomin2018small,bressloff2018variational}. Recall that $v_t = u_t - \mathcal{T}_{\beta_t} u_*$.
\begin{theorem}
Under Assumptions \ref{reg_f_rot}, \ref{Reg_stochRot}, \ref{regularity_rotWaves}, \ref{NegDefMat}, \ref{SpectAsyRot}, \ref{integrabilityAssStochConv} let $p \in[0,2) .$ There exists a constant $\tilde{C}>0$ (independent of the choice of $p$ ) and $\varepsilon_{(p)}>0$ such that for all $\varepsilon \in\left(0, \varepsilon_{(p)}\right)$, and all $\kappa \in\left[\varepsilon_{(p)}^{p}, \bar{\kappa}\right]$ and all $T>0$ we have
\begin{equation}\label{eq: main stability result}
\mathbb{P}\left(\sup _{t \in[0, T]}\left\|v_{t}\right\|>\kappa\right) \leq \tilde{C}T \exp \left(-\tilde{C} \varepsilon^{-2} \kappa^{2}\right)
\end{equation}
for a constant $\tilde{C} >0$ independent of $\varepsilon$ and $\kappa$.
\end{theorem}

\proof
In this proof, we adapt some of the ideas from the proof of Theorem 5.1 in \cite{maclaurin2020metastability} to the case of stochastic rotating waves with noise term vanishing outside a small neighbourhood of the rotating wave initial condition. We must choose a time-discretization that is long enough that the exponential decay of the linearization about the wave is strong enough to damp down the initial `difference' between the nearest translated wave and the solution at time $t_a$. Thanks to a Taylor Expansion, the remainder of the drift term is quadratic in the amplitude, and therefore negligible as long as $\kappa$ in \eqref{eq: main stability result} is small enough. Then all that remains is to bound an exponential martingale - and we can use existing results from SPDEs to do this.

Let consider a time discretization into intervals of length $\Delta t=b^{-1} \log \left(4 C^{-1}\right)$, and we write $t_{a}=a \Delta t$. Additionally, we simplify the notation writing $u_{a}=u_{t_{a}}$ and $\beta_{a}=\beta_{t_{a}}$ and we define the event
$$
\mathcal{A}_{a}=\left\{\left\|v_{a}\right\| \leq \kappa /(2 C)\right\} \cap\left\{\left\|v_{a+1}\right\|>\kappa /(2 C) \text { or } \sup _{t \in\left[t_{a}, t_{a+1}\right]}\left\|v_{t}\right\|>\kappa\right\}
$$
Writing $R=\lfloor T / \Delta t\rfloor $ and observing that $\left\|w_{0}\right\|=0$, 
we have 
\begin{equation}\label{eq: event decomposition time discretization}
\left\{\sup _{t \in[0, T]}\left\|v_{t}\right\|>\kappa\right\}\subset \bigcup_{a=0}^R\mathcal{A}_{a}
\end{equation}
and, therefore, using monotonicity and sub-additivity of probability measures, one finds that
$$
\mathbb{P}\left(\sup _{t \in[0, R \Delta t]}\left\|v_{t}\right\|>\kappa\right) \leq \sum_{a=0}^{R} \mathbb{P}\left(\mathcal{A}_{a}\right)
$$
Next define $w_{a,t} = u_t - \mathcal{T}_{\beta_a} u_*$, for $t\in [t_a, t_{a+1}]$. Although the theorem is formulated in terms of the variational phase $v_t$, we will use the triangle inequality to write, for $t\in [t_a, t_{a+1}]$,
\begin{equation}
\| v_t \| \leq \| w_{a,t} \| + \| v_t - w_{a,t} \|,
\end{equation}
and we separately bound each term on the right hand side. 

Define $\mathcal{L}_a$ to be the linearization of the deterministic dynamics about $\mathcal{T}_{\beta_a}u_*$, in a frame that is co-rotating about the Euclidean point $(\beta^1_a , \beta^2_a)$, with initial rotation equal to $\beta^0_a$. We thus see that $\mathcal{L}_a\cdot z = \Delta \mathcal{T}_{\beta_a}\cdot z +\omega_*\mathcal{T}_{\beta_a}\partial_{\psi}\mathcal{T}_{-\beta_a}+ f'\big( \mathcal{T}_{\beta_a} u_* \big)$, and let $\mathcal{L}_a^*$ be the adjoint. 

We now change variable, writing the SPDE to be relative to the reference frame that is rotating about the point with Euclidean coordinates $(\beta^1_a,\beta^2_a)$. For clarity, we write the solution of the SPDE relative to the off-center co-rotating reference frame as $y_{a,t}$: i.e. $y_{a,t} := \mathcal{S}_a(t)u_t$, where $u_t$ is the mild solution of the original SPDE in the static reference frame \eqref{react_diff_comoving_pert}, and we have written $\mathcal{S}_a(t)$ to be the linear operator that rotates about the Euclidean point $(\beta^1_a,\beta^2_a)$ by an angle $\beta^0_a + \omega^* (t-t_a)$. In other words, $\mathcal{S}_a(t)$ is the semigroup generated by $\omega_*\mathcal{T}_{\beta_a}\partial_{\psi}\mathcal{T}_{-\beta_a}$. We thus find that $y_{a,t}$ must be the mild solution of the following SPDE,
\begin{equation}\label{eq: yt spde}
dy_{a,t} = \big\lbrace \Delta y_{a,t} +\omega_*\mathcal{T}_{\beta_a}\partial_{\psi}\mathcal{T}_{-\beta_a}y_{a,t} + f(y_{a,t}) \big\rbrace dt + \hat{B}_a(t,y_{a,t})dW_t .
\end{equation}
where
\begin{equation}
\hat{B}_a(t,y_{a,t}) := \mathcal{S}_a(t)B(t,\mathcal{S}_a(-t)y_{a,t}),
\end{equation}
and $\mathcal{S}_a(-t)$ is the inverse of $\mathcal{S}_a(t)$. Indeed, to verify that the mild solution of \eqref{eq: yt spde} is consistent with the mild solution to the original SPDE in \eqref{react_diff_comoving_pert}, we can exploit the fact that the diffusion operator commutes with the translation / rotation operator, which means that their semigroups must commute as well. This allows us to write the mild solution of \eqref{eq: yt spde} as
\begin{multline}
y_{a,t} = \exp\big((t-t_a)\Delta \big) \mathcal{S}_a(t - t_a)y_{a,t_a} + \int_{t_a}^t \exp\big((t-s)\Delta \big) \mathcal{S}_a(t - s)f(y_{a,s}) ds\\ + \int_{t_a}^t \exp\big((t-s)\Delta \big) \mathcal{S}_a(t - s)\hat{B}_a(s,y_{a,s})dW_s.
\end{multline}
Since $f$ is `local', meaning that it commutes with the rotation operator, one obtains that
\begin{align}
\mathcal{S}_a(-t)y_{a,t} =& \exp\big((t-t_a)\Delta \big) \mathcal{S}_a( - t_a)y_{a,t_a} + \int_{t_a}^t \exp\big((t-s)\Delta \big) \mathcal{S}_a( - s)f(y_{a,s}) ds\nonumber \\ &+ \int_{t_a}^t \exp\big((t-s)\Delta \big) \mathcal{S}_a( - s)\hat{B}_a(s,y_{a,s})dW_s \\
=& \exp\big((t-t_a)\Delta \big)u_a + \int_{t_a}^t \exp\big((t-s)\Delta \big) f(u_s) ds + \int_{t_a}^t \exp\big((t-s)\Delta \big)B(s,u_s)dW_s,
\end{align}
and we see that indeed, as required, $\mathcal{S}_a(-t)y_{a,t}$ corresponds to the mild solution of the SPDE in the static reference frame (i.e. \eqref{react_diff_comoving_pert}).

Lets now define $w_{a,t}$, for $t\in [t_a, t_{a+1}]$ to be the amplitude, i.e. $w_{a,t} := y_{a,t} - \mathcal{T}_{\beta_a}u_*$ so that
\begin{equation}
dw_{a,t} = \big\lbrace \mathcal{L}_a  w_{a,t} + \mathcal{H}_t \big\rbrace dt + \hat{B}_a(t,y_{a,t}) dW_t
\end{equation}
where
$$
\begin{aligned}
&\mathcal{H}_{s}:= f\left(\mathcal{T}_{\beta_a}u_*+w_{a,s}\right)-f\left(\mathcal{T}_{\beta_a}u_*\right)-f^{\prime}\left(\mathcal{T}_{\beta_a} u_*\right)w_{a,s}.
\end{aligned}
$$
Next we take the mild solution over $[t_a, t_{a+1}]$, 
\begin{equation}\label{variaConst_transverseSPDE}
\begin{aligned}
w_{a,t}=&\exp\left(\mathcal{L}_a\left(t-t_{a}\right)\right) w_{a,t_a}+\int_{t_{a}}^{t} \exp\left(\mathcal{L}_a\left(t-s\right)\right) \mathcal{H}_{s} d s\\
&+\varepsilon \int_{t_{a}}^{t} \exp\left(\mathcal{L}_a\left(t-s\right)\right) B(s, u_{s}) d W_{s}.
\end{aligned}
\end{equation}
It can be shown that, since $f'\big( \mathcal{T}_{\beta_a}u_*\big)$ is linear and bounded, the above mild solution is consistent with the mild solution obtained by taking the semigroup generated by the linear operators with spatial derivatives only. Finally the triangle inequality implies that
\begin{equation}\label{variaConst_transverseSPDE 2}
\begin{aligned}
\| w_{a,t} \| \leq &\| \exp\left(\mathcal{L}_a\left(t-t_{a}\right)\right) w_{a,t_a} \| +\| \int_{t_{a}}^{t} \exp\left(\mathcal{L}_a\left(t-s\right)\right) \mathcal{H}_{s} d s \| \\
&+\varepsilon \| \int_{t_{a}}^{t} \exp\left(\mathcal{L}_a\left(t-s\right)\right) B(s, u_{s}) d W_{s} \|.
\end{aligned}
\end{equation}
Our goal is to prove that there exist constants $\tilde{C}_1,\tilde{C}_2$ such that
\begin{multline}\label{eq: norm va 3}
\big\lbrace \| v_a \| \leq \kappa / (2C) \text{ and either }\| v_{a+1} \| > \kappa / (2C)\text{ or }\sup_{s\in [t_a, t_{a+1}]} \| v_{s}\|^2 > \kappa    \big\rbrace \subseteq  \\ \big\lbrace  \varepsilon \tilde{C}_1  \sup_{s\in [t_a,  t_{a+1}]} \sup_{0\leq i \leq 2}  \bigg| \int_{t_a}^{s} \mathcal{Y}_i(r,u_r,\beta_r) dW_r \bigg| \geq \kappa   \big\rbrace \cup \\ \bigg\lbrace  \varepsilon  \tilde{C}_2 \sup_{t\in [t_a, t_{a+1}]} \| \int_{t_{a}}^{t} \exp\left(\mathcal{L}_a\left(t-s\right)\right) \hat{B}_a(s, y_{a,s}) d W_{s} \| \geq \kappa  \bigg\rbrace .
\end{multline}
The probabilities of the events onthe right hand side can then be bounded using Chernoff's Inequality. We start by bounding the term with $\mathcal{H}_s$. Indeed a second order Taylor Expansion implies that
\begin{equation}
\mathcal{H}_s = f''(z)\cdot w_{a,s} \cdot w_{a,s},
\end{equation}
for some $z$ in the convex hull of $u_s$ and $\mathcal{T}_{\beta_a}u_*$. By assumption the second Frech\'et Derivative of $f$ is bounded. Also the operator norm of the semigroup 
$ \exp\left(\mathcal{L}_a\left(t-t_a\right)\right)$ is uniformly bounded by $C$, i.e. 
\begin{equation}
 \sup_{a\in \mathbb{Z}^+}\sup_{t\in [t_a, t_{a+1}]}\| \exp\left(\mathcal{L}_a\left(t-t_a\right)\right) \| \leq C.
\end{equation}
We thus find that there must exist  a constant $\bar{C} > 0$ such that for all $t \in\left[t_{a}, t_{a+1} \right]$, and recalling that $t_{a+1} - t_a := \Delta t$, 
\begin{equation}\label{Lemma Taylor}
\begin{aligned}
&\| \int_{t_{a}}^{t} \exp\left(\mathcal{L}_a\left(t-s\right)\right) \mathcal{H}_{s} d s \| \leq \bar{C}\Delta t \sup _{s \in\left[t_{a}, t\right]}\left\|w_{a,s}\right\|^{2} 
\end{aligned}
\end{equation}
Let $\bar{\kappa} > 0$ (this is a uniform upper bound for $\kappa$) be such that
\begin{equation} \label{eq: bar kappa upper bound}
\bar{C}\Delta t \bar{\kappa}^2 \leq \bar{\kappa} / (8C).
\end{equation}
We now claim that, as long as $\kappa \in (0, \bar{\kappa})$ (where $\bar{\kappa}$ satisfies \eqref{eq: bar kappa upper bound}), it must be that
\begin{multline} \label{eq: to show w a t a increments}
\big\lbrace \| w_{a,t_a} \| \leq \kappa / (2C) \text{ and }\| w_{a,t_{a+1}} \| > 3\kappa / (8C) \big\rbrace\\ \subseteq  \bigg\lbrace \varepsilon \sup_{t\in [t_a, t_{a+1}]} \| \int_{t_{a}}^{t} \exp\left(\mathcal{L}_a\left(t-s\right)\right) \hat{B}_a(s, y_{a,s}) d W_{s} \| > \kappa / (8C) \bigg\rbrace .
\end{multline}
To see \eqref{eq: to show w a t a increments}, we first prove that if $\| w_{a,t_a} \| \leq \kappa / (2C)$ and
\begin{equation}\label{eq: upperbound first stochastic integral}
\varepsilon \sup_{t\in [t_a, t_{a+1}]} \| \int_{t_{a}}^{t} \exp\left(\mathcal{L}_a\left(t-s\right)\right) \hat{B}_a(s, y_{s}) d W_{s} \| \leq \kappa / (8C),
\end{equation}
then necessarily
\begin{equation} \label{eq: to prove w a s intermediate bound}
\sup_{s\in [t_a, t_{a+1}]} \| w_{a,s} \| \leq 3\kappa / 4.
\end{equation}
Suppose for a contradiction that \eqref{eq: to prove w a s intermediate bound} does not hold: then there must exist $\zeta \in [t_a, t_{a+1}]$ such that $\| w_{a,\zeta} \| = 3\kappa / 4$. Take $\zeta$ to be the smallest time $t$ greater than $t_a$ such that $\| w_{a,t} \| = 3\kappa /4$. Then clearly $\sup_{t \in [t_a,\zeta]}\| w_{a,t}\| \leq 3\kappa / 4$, and since $\kappa \leq \bar{\kappa}$ (which satisfies \eqref{eq: bar kappa upper bound}), it follows from \eqref{Lemma Taylor} that
\[
\| \int_{t_{a}}^{\zeta} \exp\left(\mathcal{L}_a\left(\zeta-s\right)\right) \mathcal{H}_{s} d s \| \leq \frac{\kappa}{8C}.
\]
The action of the semigroup on the orthogonal space has norm upperbounded by $C$, and therefore
\[
\left\| \exp\left(\mathcal{L}_a\left(\zeta-t_{a}\right)\right) w_{a,t_a} \right\| \leq  \| w_{a,t_a}\| C \leq  \frac{\kappa}{2}
\]
We thus obtain from \eqref{variaConst_transverseSPDE 2} that $\| w_{a,\zeta}\| < 3\kappa / 4$ (since we can take $C \geq 1$, since it is by definition an upperbound), a contradiction. We can therefore conclude that if $\| w_{a,t_a} \| \leq \kappa / (2C)$ and also \eqref{eq: upperbound first stochastic integral} holds, then necesarily \eqref{eq: to prove w a s intermediate bound} must also hold, as required.

We can now show \eqref{eq: to show w a t a increments}. Indeed, thanks to the exponential decay bound in \eqref{linearizationOperRotatGamma}, 
\begin{equation}\label{eq: orthogonal space decay}
\left\| \exp\left(\mathcal{L}_a\left(t_{a+1}-t_{a}\right)\right) w_{a,t_a} \right\| \leq C\| w_{a,t_a} \| \exp(-b(t_{a+1}-t_a)) \leq \| w_{a,t_a} \| / 4,
\end{equation}
using the definition $\Delta t=b^{-1} \log \left(4 C^{-1}\right)$. We thus find from \eqref{variaConst_transverseSPDE 2} and the above identities that, if \eqref{eq: upperbound first stochastic integral} holds, then
\begin{equation}
\| w_{a,t_{a+1}} \| \leq \frac{\kappa}{8C} + \frac{\kappa}{8C} + \frac{\kappa}{8C},
\end{equation}
which obviously implies that $\| w_{a,t_{a+1}} \| \leq \frac{3\kappa}{8C}$. We can thus conclude that \eqref{eq: to show w a t a increments} is true (we have proved its contrapositive).

It remains to control the difference between $v_s$ and $w_{a,s}$. Since the manifold $\lbrace \mathcal{T}_{\gamma} u_* \rbrace$ depends smoothly on $\gamma$, there must exist a constant $\hat{C} > 0$ such that
\begin{equation}\label{eq: vs bound}
\sup_{s\in [t_a,t_{a+1}]}\| v_s \| \leq \hat{C}\big(\sup_{s\in [t_a,t_{a+1}]}\| w_{a,s} \| + \sup_{s\in [t_a,t_{a+1}]}\big| \beta_s - \beta_{a} \big|  \big).
\end{equation}
In turn, the SDE for $\beta_s$ implies that there exist constants $C_1,C_2,C_3$ such that
\begin{multline}
     \sup_{s\in [t_a,t_{a+1}]}\big| \beta_s - \beta_{a} \big| \leq C_1 \varepsilon^2 \Delta t + C_2 \Delta t \sup_{s\in [t_a,t_{a+1}]}\| v_s \|^2 + C_3 \varepsilon  \sup_{s\in [t_a, t_{a+1}]} \sup_{0\leq i \leq 2}  \bigg| \int_{t_a}^{s} \mathcal{Y}_i(r,u_r,\beta_r) dW_r \bigg| .
\end{multline}
It is then straightforward to show that, using similar methods, \eqref{eq: norm va 3} holds. It remains to demonstrate that the probability of the events on the right hand side of \eqref{eq: norm va 3} is exponentially small in $\varepsilon^2$.

In summary, and noting the event decomposition in \eqref{eq: event decomposition time discretization}, we may conclude that for constants $C_1,C_2 > 0$ that are independent of $\varepsilon$ and $\kappa$ (which one could in principle infer from the above proof) and for $\kappa \in (0,\bar{\kappa})$,
\begin{equation}
\mathbb{P}\left\{\sup _{t \in[0, T]}\left\|v_{t}\right\|>\kappa\right\} \leq C_1(R+1) \exp\big( -C_2 \kappa^2 \varepsilon^{-2} \big),
\end{equation}
and we recall that $R$ is the number of time intervals. Thus $R = O(T)$ - since the width of the time intervals is also chosen independently of $\varepsilon$, and we have thus proved the main result.
Indeed it follows from \cite{StochConvBound} and Assumption \ref{integrabilityAssStochConv} that for some constant $\tilde{C} > 0$,
\begin{equation}
\mathbb{P}\big( \varepsilon  C_3 \sup_{t\in [t_a, t_{a+1}]} \| \int_{t_{a}}^{t} \exp\left(\mathcal{L}_a\left(t-s\right)\right) B(s, u_{s}) d W_{s} \| > \kappa / (8C) \big) \leq \exp\big( - \tilde{C} (\kappa \varepsilon^{-1})^2 \big).
\end{equation}
One similarly shows, using Chernoff's Inequality (see for instance \cite{NeuralFields_MACLAURIN}), that
\begin{equation}
 \mathbb{P}\big( C_2 \varepsilon  \sup_{s\in [t_a, \eta \wedge t_{a+1}]} \sup_{0\leq i \leq 2}  \bigg| \int_{t_a}^{s} \mathcal{Y}_i(r,u_r,\beta_r) dW_r \bigg| \geq \kappa / (8C) \big) \leq  \exp\big( - \hat{C} (\kappa \varepsilon^{-1})^2 \big),
\end{equation}
for some constant $\hat{C} > 0$.
\endproof

\subsection{Stability results - Approximated Variational Phase}

In this section we will use the approximated variational phase we constructed to show a stability result for stochastic rotating waves, which is easier to obtain but requires stronger assumptions, which make it effectively valid locally-in-time with respect to the fixed reference frame. In this section, we will additionally assume
\begin{assumption}\label{integrabilityAssStochConv}
For any $t>s$ and any $x \in \mathcal{H}, \exp((t-s) \mathcal{L}_*) B(s, x)$ is a Hilbert-Schmidt operator, with Hilbert-Schmidt norm upperbounded by
$$
\int_{0}^{T_{0}}(t-s)^{-2 \xi} \sup _{x \in H}\left\|\exp((t-s) \mathcal{L_*}) B(s, x)\right\|_{H S}^{2} d s<C_{H S}
$$
for some $\xi \in(0,1 / 2)$ and constant $C_{H S}<\infty$, where
$$
T_{0}=\frac{\log( 4 C)}{b}
$$\end{assumption}
First of all, we define the stopping time $\tau_1=\inf\{t\geq 0 : ||w_t||=\kappa\}$ for a generic $\kappa>0$.We will add some requirements on $\kappa$ later. We recall the approximated variational phase SDE \eqref{approxVarSDERotating} and its solution $\beta_t$. Notice that for all $\xi>0$ there exists $\delta>0$ such that if $||u_t-u_*||<\delta $ than $|\beta_t|<\xi$, see Lemma 2.5 in \cite{nonlinearstab_Rotwaves}. Notice also that it exists a $0<\bar{\xi}$ small enough such that for $|\beta_t|<\bar{\xi}$ the matrix $M(t,\beta_t)$ is invertible and has all eigenvalues bigger than a constant $C_{eig}(\bar{\xi})>0$. Therefore, the inverse  $N(t,\beta_t)=M^{-1}(t,\beta_t)$ has eigenvalues smaller than $1/C_{eig}(\bar{\xi})$. We will make an additional requirement on $\bar{\xi}$: \begin{equation}
\label{C_4Cond}
C_4\bar{\xi}<\frac{\kappa}{32C}
\end{equation}
for a constant $C_4>0$ that will be defined later. For such a $\bar{\xi}$ we will assume in the rest of this section
\begin{assumption}
\label{NoiseDies}
For $\delta>0$ smaller than $\varepsilon$ in Section~\ref{nonlStabRotWavesTHM} and such that if $\|u_t-u_*\|<\delta$ then $|\beta_t|<\frac{\bar{\xi}}{2}$, we have that $ \|z-u_*\|\geq\delta_1$ for $\delta_1$ such that $\delta_1<\max{\{\delta, \frac{\bar{\xi}}{4C}\}}$ implies
$$B(t,z)=0$$
\end{assumption} 
 In this case, we have that at every time $t>0$ we have $N(t,\beta_t)=M^{-1}(t,\beta_t)$ has eigenvalues smaller than $1/C_{eig}(\bar{\xi})$. To show this we just have to show that $|\beta_t|<\bar{\xi}$. First of all, we notice that this is the case when $\|u_t-u_*\|<\delta_1<\delta$. In fact, we have chosen $\delta$ such that $|\beta_t|<\bar{\xi}$. When, instead, $\|u_t-u_*\|\geq\delta_1$ the stochastic part of the equation for $u_t$ becomes zero. Therefore, $u_t$ follows deterministic dynamics since when $\|u_t-u_*\|=\delta_1$ until when we will have $\|u_t-u_*\|<\delta_1<\delta$ again (if this is the case). We now simplify the notation restarting the time from $0$ when we have $\|u_t-u_*\|=\delta_1$. From 5) in Theorem \ref{nonlStabRotWavesTHM} we have that for every $t\geq 0$ when the equation follows deterministic dynamics it exists $\beta_{\infty}$ such that  $$|\begin{bmatrix}1&0\\
0& R_{\omega_*t}
\end{bmatrix}\beta_t-\beta_{\infty}|\leq C\exp(-\frac{b}{2} t)\|u_0-u_*\|\leq C\delta_1$$for every time $t\geq 0$. From the previous inequality it follows
$$
\begin{aligned}
|\begin{bmatrix}1&0\\
0& R_{\omega_*t}
\end{bmatrix}\beta_t-\beta_0|&\leq |\begin{bmatrix}1&0\\
0& R_{\omega_*t}
\end{bmatrix}\beta_t-\beta_{\infty}|+|\beta_0-\beta_{\infty}|\\
&\leq 2C\delta_1<\frac{\bar{\xi}}{2}
\end{aligned}
$$We notice again that we have $\|u_0-u_*\|=\delta_1<\delta$ and thus $|\beta_0|<\frac{\bar{\xi}}{2}$. Therefore, we obtain
$$
\begin{aligned}
|\beta_t|&=|\begin{bmatrix}1&0\\
0& R_{\omega_*t}
\end{bmatrix}\beta_t|<|\beta_0|+|\begin{bmatrix}1&0\\
0& R_{\omega_*t}
\end{bmatrix}\beta_t-\beta_0|
<\bar{\xi}
\end{aligned}
$$This is very important for the following arguments.

Now remember that $w_t=u_t-\mathcal{T}_{\beta_t}u_*$ and we have the following equation
$$
\begin{aligned}
dw_t=&\left[\mathcal{L}_*w_t+f\left(\mathcal{T}_{\beta_t}u_*+w_t\right)-f\left(\mathcal{T}_{\beta_t}u_*\right)-f^{\prime}\left(u_*\right)w_t\right] d t+\varepsilon B\left(t,u_{t}\right) d W_{t}\\
&-\sum_{i=0}^2\partial_{\beta_t^i}\mathcal{T}_{\beta_t}u_*d\beta_t^i
-\sum_{i,j=0}^2\partial_{\beta_t^i\beta_t^j}\mathcal{T}_{\beta_t}u_*d\beta_t^i d\beta_t^j\\
=&\left[\mathcal{L}_*w_t+f\left(\mathcal{T}_{\beta_t}u_*+w_t\right)-f\left(\mathcal{T}_{\beta_t}u_*\right)-f^{\prime}\left(u_*\right)w_t+\mathcal{T}_{\beta_t} \left(\nabla u_*\right)^T R_{-\beta^0_t+\frac{\pi}{2}}\begin{bmatrix}\beta_t^1 \\
\beta_t^2 \\
\end{bmatrix}\right] d t\\
&+\varepsilon B\left(t,u_{t}\right) d W_{t}-\sum_{i=0}^{2} \partial_{\beta_t^i}\mathcal{T}_{\beta_t}u_* \mathcal{Y}_{i}\left(t, u_{t}, \beta_t\right) dW_t\\
&-\sum_{i=0}^{2} \partial_{\beta_t^i}\mathcal{T}_{\beta_t}u_* \mathcal{V}_{i}\left(t, u_{t}, \beta_t\right)dt \\
 &-\frac{\varepsilon^{2}}{2} \sum_{j, k=0}^{2} \partial_{\beta_t^i\beta_t^j}\mathcal{T}_{\beta_t}u_*d\beta^j_t d\beta_t^k
 \\
=&\left[\mathcal{L}_*w_t+f\left(\mathcal{T}_{\beta_t}u_*+w_t\right)-f\left(\mathcal{T}_{\beta_t}u_*\right)-f^{\prime}\left(u_*\right)w_t\right] d t+\varepsilon B\left(t,u_{t}\right) d W_{t}\\
 &-\sum_{i=0}^{2} \partial_{\beta_t^i}\mathcal{T}_{\beta_t}u_* \mathcal{Y}_{i}\left(t, u_{t}, \beta_t\right) dW_t-\sum_{i=0}^{2} \partial_{\beta_t^i}\mathcal{T}_{\beta_t}u_* \widetilde{\mathcal{V}}_{i}\left(t, u_{t}, \beta_t\right)dt \\
 &-\frac{\varepsilon^{2}}{2} \sum_{j, k=0}^{2} \partial_{\beta_t^i\beta_t^j}\mathcal{T}_{\beta_t}u_*d\beta^j_t d\beta_t^k
\end{aligned}
$$

We can now prove the following theorem
\begin{theorem}
Under Assumptions \ref{reg_f_rot}, \ref{Reg_stochRot}, \ref{regularity_rotWaves}, \ref{NegDefMat}, \ref{SpectAsyRot}, \ref{integrabilityAssStochConv} let $p \in[0,2) .$ There exists a constant $C>0$ (independent of the choice of $p$ ) and $\varepsilon_{(p)}>0$ such that for all $\varepsilon \in\left(0, \varepsilon_{(p)}\right)$, and all $\kappa \in\left[\varepsilon_{(p)}^{p}, \bar{\kappa}\right]$ and all $T>0$ we have
$$
\mathbb{P}\left(\sup _{t \in[0, T]}\left\|w_{t}\right\|>\kappa\right) \leq C( T) \exp \left(-C \varepsilon^{-2} \kappa^{2}\right)
$$
for a constant $C(T)>0$ independent from $\varepsilon$ and $\kappa$.
\end{theorem}

\proof
In this proof, we adapt some of the ideas from the proof of Theorem 5.1 in \cite{maclaurin2020metastability} to the case of stochastic rotating waves with noise term vanishing outside a small neighbourhood of the rotating wave initial condition.

Let consider a time discretization into intervals of length $\Delta t=\omega^{-1} \log \left(4 C^{-1}\right)$, and we write $t_{a}=a \Delta t$. Additionally, we simplify the notation writing $u_{a}=u_{t_{a}}$ and $\beta_{a}=\beta_{t_{a}}$ and we define the event
$$
\mathcal{A}_{a}=\left\{\left\|w_{a}\right\| \leq \kappa /(2 C)\right\} \cap\left\{\left\|w_{a+1}\right\|>\kappa /(2 C) \text { or } \sup _{t \in\left[t_{a}, t_{a+1}\right]}\left\|w_{t}\right\|>\kappa\right\}
$$
Writing $R=\lfloor T / \Delta t\rfloor $ and observing that $\left\|w_{0}\right\|=0$, 
we have 
$$
\left\{\sup _{t \in[0, T]}\left\|w_{t}\right\|>\kappa\right\}\subset \bigcup_{a=0}^R\mathcal{A}_{a}
$$
and, therefore, using monotonicity and sub-additivity of probability measures, one have
$$
\mathbb{P}\left(\sup _{t \in[0, R \Delta t]}\left\|w_{t}\right\|>\kappa\right) \leq \sum_{a=0}^{R} \mathbb{P}\left(\mathcal{A}_{a}\right)
$$

Now we consider again the mild solution
\begin{equation}\label{variaConst_transverseSPDE}
\begin{aligned}
w_{t}=&\exp\left(\mathcal{L}_*\left(t-t_{a}\right)\right) w_{a}+\int_{t_{a}}^{t} \exp\left(\mathcal{L}_*\left(t-s\right)\right) \mathcal{H}_{s} d s\\
&+\varepsilon \int_{t_{a}}^{t} \exp\left(\mathcal{L}_*\left(t-s\right)\right)\widetilde{B}\left(s, u_{s}, \beta_{s}\right) d W_{s}
\end{aligned}
\end{equation}
where
$$
\begin{aligned}
&\widetilde{B}(s, z, \alpha):  \mathbb{R}^{+} \times H^2(\mathbb{R}^2,\mathbb{R}^N) \times \mathbb{R}^{3} \rightarrow \mathcal{L}(H, H) \\
&\widetilde{B}(s, z, \alpha)= B(s, z)-\sum_{j=0}^{2} \partial_{\alpha^j}\mathcal{T}_{\alpha}u_* \mathcal{Y}_{j}(s, z, \alpha) \\
&\mathcal{H}_{s}:= f\left(\mathcal{T}_{\beta_s}u_*+w_s\right)-f\left(\mathcal{T}_{\beta_s}u_*\right)-f^{\prime}\left(u_*\right)w_s-\sum_{i=0}^{2} \partial_{\beta_t^i}\mathcal{T}_{\beta_t}u_* \widetilde{\mathcal{V}}_{i}\left(s, u_{s}, \beta_{s}\right) \\
&-\frac{\varepsilon^{2}}{2} \sum_{j, k, p, q=0}^{2} \partial_{\beta_t^j \beta_t^k}\mathcal{T}_{\beta_t}u_*\mathcal{N}_{j p}\left(u_{s}, \beta_{s}\right) \mathcal{N}_{k q}\left(u_{s}, \beta_{s}\right)\left\langle B^{ad}\left(s, u_{s}\right) \eta^{p}, B^{ad}\left(s, u_{s}\right) \eta^{q}\right\rangle
\end{aligned}
$$
and we show the following lemma
\begin{lemma}
There exist constants $C_{1}, C_{2}, C_{3},C_4$ such that for all $t \in\left[t_{a}, t_{a+1} \wedge \eta\right]$, where $t_{a}=\sup \left\{t_{b}\right.$ :
$\left.t_{b} \leq t\right\}$
$$
\begin{aligned}
&\int_{t_{a}}^{t} \exp\left(\mathcal{L}_*\left(t-s\right)\right) \mathcal{H}_{s} d s \leq \varepsilon^{2} C_{1}+C_{3} \sup _{s \in\left[t_{a}, t\right]}\left\|w_{s}\right\|^{2} +C_4\sup _{s \in\left[t_{a}, t\right]}|\beta_t| \sup _{s \in\left[t_{a}, t\right]}\left\|w_{s}\right\|\\
&+\varepsilon C_{2} \sup _{s \in\left[t_{a}, \eta \wedge t_{a+1}\right]} \sup _{0\leq i \leq 2}\left|\int_{t_{a}}^{s} \mathcal{Y}_{i}\left(r, u_{r}, \boldsymbol{\beta}_{r}\right) d W_{r}\right|
\end{aligned}
$$
where $\beta_t$ is the approximated variational phase.
\end{lemma}
\proof From the exponential decay bound of the semigroup $\exp(t\mathcal{L}_*)$, the definition of $P$ and the Cauchy-Schwartz inequality we have
$$
\begin{aligned}
\left\|\exp(\mathcal{L}_*(t-s))\right\|_{\mathcal{L}} & \leq\left\|\exp(\mathcal{L}_*(t-z)))(I-P)\right\|_{\mathcal{L}}+\left\|P\right\|_{\mathcal{L}} \\
& \leq C+\sum_{i=0}^{2}\left\|\partial_i u_*\right\|\left\|\eta^{i}\right\|
\end{aligned}
$$Aditionally, we have a  constant $\bar{C}$ such that
\begin{equation}\label{BoundVY}
\begin{gathered}
\sup _{ t \in [0,T], \text{ } \|u_t-\mathcal{T}_{\beta_t}u_*\|\leq \kappa}\sup _{0\leq i \leq 2} \left|\widetilde{\mathcal{V}}_{i}(t, u_t, \beta_t)\right| \leq \bar{C} \\
\sup _{ t \in [0,T], \text{ } \|u_t-\mathcal{T}_{\beta_t}u_*\|\leq \kappa}\sup _{0\leq i \leq 2}\left\|\mathcal{Y}_{i}(t, u_t, \beta_t)\right\|_{H S} \leq \bar{C} .
\end{gathered}
\end{equation}
This follows trivially from the bound on the eigenvalues of $\mathcal{N}(t,\beta_t)$ given by \ref{NoiseDies}, the Cauchy-Schwartz inequality and the inequality
\begin{equation}\label{inequalityFF}
\left\|f\left(\mathcal{T}_{\gamma}u_*+w\right)-f\left(\mathcal{T}_{\gamma}u_*\right)-f^{\prime}\left(u_*\right)w\right\|\leq C_f\left(\left\|w\right\|+|\gamma|\right)\| w\|
\end{equation}that can be proved exactly in the same way as follows.

Employing the triangle inequality, we obtain
$$
\begin{aligned}
&\left\|f\left(\mathcal{T}_{\gamma}u_*+w\right)-f\left(\mathcal{T}_{\gamma}u_*\right)-f^{\prime}\left(u_*\right)w\right\|\leq\\
&\left\|f\left(\mathcal{T}_{\gamma}u_*+w\right)-f\left(\mathcal{T}_{\gamma}u_*\right)-f^{\prime}\left(\mathcal{T}_{\beta_t}u_*\right)w\right\|+\left\|f^{\prime}\left(\mathcal{T}_{\beta_t}u_*\right)w_t-f^{\prime}\left(u_*\right)w_t\right\|
\end{aligned}
$$
Using the second order Taylor expansion, we know there exists $\lambda \in[0,1]$ such that, writing $\bar{u}=\lambda \mathcal{T}_{\gamma}u_*+\left(1-\lambda\right) \left(\mathcal{T}_{\gamma}u_*+w\right)$, 
$$
f\left(\mathcal{T}_{\gamma}u_*+w\right)-f\left(\mathcal{T}_{\gamma}u_*\right)-f^{\prime}\left(\mathcal{T}_{\gamma}u_*\right) \cdot w_{a}=f^{\prime \prime}\left(\bar{u}\right) \cdot w \cdot w
$$The assumed boundedness of the second derivative implies that for some constant $C_1>0$,
$$
\left\|f\left(\mathcal{T}_{\gamma}u_*+w\right)-f\left(\mathcal{T}_{\gamma}u_*\right)-f^{\prime}\left(\mathcal{T}_{\gamma}u_*\right) w\right\| =\left\| f^{\prime \prime}\left(\bar{u}\right) \cdot w \cdot w\right\|\leq C_1\left\|w\right\|^{2}
$$
Using Taylor expansions again, there exists $\lambda \in[0,1]$ such that, writing $\tilde{u}=\lambda \mathcal{T}_{\gamma}u_*+\left(1-\lambda\right) u_*$ we have
$$
f^{\prime}(\mathcal{T}_{\gamma}u_*)-f^{\prime}(u_*)=f^{\prime \prime}(\tilde{u})\cdot \left(\mathcal{T}_{\gamma}u_*-u_*\right)
$$
and the boundedness of the first and second Fr\'echet derivatives of $f$ imply that there exists constants $C_2>0$ and $C_3>0$ such that
$$
\begin{aligned}
\left\|f^{\prime}\left(\mathcal{T}_{\gamma}u_*\right)w-f^{\prime}\left(u_*\right)w\right\|&=\| f^{\prime \prime}(\tilde{u})\cdot \left(\mathcal{T}_{\gamma}u_*-u_*\right)\|\\ 
\leq C_2 \| \left(\mathcal{T}_{\gamma}u_*-u_*\right)\|\leq C_3 |\gamma|
&
\end{aligned}
$$
, it follows
$$\begin{aligned}
&\left\|f\left(\mathcal{T}_{\beta_t}u_*+w_t\right)-f\left(\mathcal{T}_{\beta_t}u_*\right)-f^{\prime}\left(u_*\right)w\right\|\leq \\
&C_f\left(\sup _{t \in\left[t_{a}, t_{a+1}\right]}\left\|w_t\right\|+\sup _{t \in\left[t_{a}, t_{a+1}\right]}|\beta_t|\right)\sup _{t \in\left[t_{a}, t_{a+1}\right]}\| w_t\|
\end{aligned}$$
From this last observation, Cauchy-Schwartz inequality and the bounds \eqref{BoundVY}, the statement of the lemma follows directly.
\endproof

Now, we require that $\varepsilon_{(p)}$ (defined in the statement of the theorem) satisfies
$$
0<\varepsilon_{(p)}^{2} C_{1} \leq \varepsilon_{(p)}^{p} /(16 C)
$$
Observe that this is always possible since by assumption $p<2$. Moreover, we remark that we have 
$$
\sup _{t \in\left[0,T\right]}|\beta_t|\leq \bar{\xi}
$$
where $\bar{\xi}>0$ such that $C_4\bar{\xi}<1 / (32C)$ from Assumption~\ref{NoiseDies}. Observe that $C_4(\bar{\xi})$ can be assumed to be decreasing in $\bar{\xi}$ as the eigenvalues of $N(t,\beta_t)$ can be assumed to be decreasing in $\bar{\xi}$. Therefore, such a $\bar{\xi}$ always exists. Additionally, we add the additional requirement on $\bar{\kappa}$ 
$$
C_{3} \bar{\kappa}^{2} \leq \bar{\kappa} /(16 C)
$$
Since $\kappa \leq \bar{\kappa}$, we will always have that $C_{3} \kappa^{2} \leq \frac{\kappa}{16 c}$. With all this observations, we can now show the following lemma.

\begin{lemma}

$$
\begin{aligned}
\mathcal{A}_{a} & \subseteq\left\{\left\|w_{a}\right\| \leq \kappa /(2 C)\right\} \cap\left\{\mathcal{B}_{a} \cup \mathcal{C}_{a}\right\} \text { where } \\
\mathcal{B}_{a} &=\left\{\varepsilon C_{2} \sup _{s \in\left[t_{a}, \eta \wedge t_{a+1}\right]} \sup _{0\leq i \leq 2}\left|\int_{t_{a}}^{s \wedge \eta} \mathcal{Y}_{i}\left(r, u_{r}, \beta_{r}\right) d W_{r}\right| \geq \frac{\kappa}{16 c}\right\} \\
\mathcal{C}_{a} &=\left\{\varepsilon \sup _{t \in t_{a+1} \wedge \eta}\left\|\int_{t_{a}}^{t \wedge \eta} \exp\left(\mathcal{L}_*\left(t-s\right)\right) \widetilde{B}\left(s, u_{s}, \beta_{s}\right) d W_{s}\right\| \geq \frac{\kappa}{16 c}\right\}
\end{aligned}
$$\end{lemma}
\proof
Using the mild solution formulation \eqref{variaConst_transverseSPDE} for $w_{t}$, the triangle inequality, the exponential decay bound for the semigroup and the previous lemma we obtain, for all $t \in\left[t_{a}, \eta \wedge t_{a+1}\right]$,
$$
\begin{aligned}
\left\|w_{t}\right\| \leq &\left\|\exp\left(\mathcal{L}_*\left(t-t_a\right)\right) w_{a}\right\|+\left\|\int_{t_{a}}^{t} \exp\left(\mathcal{L}_*\left(t-s\right)\right) \mathcal{H}_{s} d s\right\|\\
&+\varepsilon\left\|\int_{t_{a}}^{t} \exp\left(\mathcal{L}_*\left(t-s\right)\right) \widetilde{B}\left(s, u_{s}, \boldsymbol{\beta}_{s}\right) d W_{s}\right\| \\
\leq & C \exp \left(-\omega\left(t-t_{a}\right)\right)\left\|w_{a}\right\|+\varepsilon^{2} C_{1}+\varepsilon C_{2} \sup _{s \in\left[t_{a}, \eta \wedge t_{a+1}\right]} \sup _{0\leq i \leq 2}\left|\int_{t_{a}}^{s} \mathcal{Y}_{i}\left(r, u_{r}, \boldsymbol{\beta}_{r}\right) d W_{r}\right| \\
&+C_{3} \kappa^{2}+C_4\sup _{t \in\left[t_{a}, t_{a+1}\right]}|\beta_t| \kappa+\varepsilon\left\|\int_{t_{a}}^{t} \exp\left(\mathcal{L}_*\left(t-s\right)\right) \widetilde{B}\left(s, u_{s}, \beta_{s}\right) d W_{s}\right\|
 \\
\leq & C \exp \left(-\omega\left(t-t_{a}\right)\right)\left\|w_{a}\right\|+C_4\sup _{t \in\left[t_{a}, t_{a+1}\right]}|\beta_t|\sup _{s \in\left[t_{a}, t\right]}\left\|w_{s}\right\|+\varepsilon^{2} C_{1}\\
&+\varepsilon C_{2} \sup _{s \in\left[t_{a}, \eta \wedge t_{a+1}\right]} \sup _{0\leq i \leq 2}\left|\int_{t_{a}}^{s} \mathcal{Y}_{i}\left(r, u_{r}, \boldsymbol{\beta}_{r}\right) d W_{r}\right| \\
&+C_{3} \kappa^{2}+C_4\xi \kappa+\varepsilon\left\|\int_{t_{a}}^{t} \exp\left(\mathcal{L}_*\left(t-s\right)\right) \widetilde{B}\left(s, u_{s}, \beta_{s}\right) d W_{s}\right\|
\end{aligned}
$$
Now,since $t-t_{a} \leq \omega^{-1} \log (4 C)$, $\left\|w_{a}\right\| \leq \kappa /(2 C)$ we have $ C\exp\left(-\omega\left(t-t_{a}\right)\right)\left\|w_{a}\right\| \leq \kappa /(8 C) .$ Furthermore, from the previous observations, we have $C_{3} \kappa^{2} \leq \kappa /(32C) $  and, since by definition $\kappa>\varepsilon_{(p)}^{p}$,  that $\varepsilon^{2} C_{1} \leq \kappa /(16C)$. Moreover, we remark again that, by assumption \ref{NoiseDies}, $C_4\bar{\xi}<\kappa/(32C)$. We thus find that
$$
\begin{aligned}
\left\|w_{a+1}\right\| \leq& \frac{\kappa}{4 C}+\varepsilon C_{2} \sup _{s \in\left[t_{a}, \eta \wedge t_{a+1}\right]} \sup _{\leq i \leq m}\left|\int_{t_{a}}^{s} \mathcal{Y}_{i}\left(r, u_{r}, \boldsymbol{\beta}_{r}\right) d W_{r}\right| \\
&+\varepsilon\left\|\int_{t_{a}}^{t_{a+1}} \exp\left(\mathcal{L}_*\left(\Delta t\right)\right)x) \widetilde{B}\left(s, u_{s}, \boldsymbol{\beta}_{s}\right) d W_{s}\right\|
\end{aligned}
$$
Therefore, if $\left\|w_{a+1}\right\|>\frac{\kappa}{2 \mathrm{c}}$, then the event $\mathcal{B}_{a}$ or $\mathcal{C}_{a}$ must happen. The remaining possibility in the event $\mathcal{A}_{a}$ is $\sup _{t \in\left[t_{a}, t_{a+1}\right]}\left\|w_{t}\right\|>\kappa .$ Observe that, since $\left\|U_{a}\left(t-t_{a}\right) w_{a}\right\| \leq C\left\|w_{a}\right\| \leq C \kappa /(2 C)=\kappa/2$,
$$
\begin{array}{r}
\sup _{t \in\left[t_{a}, t_{a+1} \wedge \eta\right]}\left\|w_{t}\right\| \leq \kappa / 2+\varepsilon^{2} C_{1}+\varepsilon C_{2} \sup _{s \in\left[t_{a}, \eta \wedge t_{a+1}\right]} \sup _{i \leq m}\left|\int_{t_{a}}^{s} \mathcal{Y}_{i}\left(r, u_{r}, \beta_{r}\right) d W_{r}\right| \\
+C_{3} \kappa^{2}+C_4\bar{\xi}\kappa+\varepsilon \sup _{t \in\left[t_{a}, t_{a+1}\right]}\left\|\int_{t_{a}}^{t \wedge \eta} \exp\left(\mathcal{L}_*\left(t-s\right)\right) \widetilde{B}\left(s, u_{s}, \beta_{s}\right) d W_{s}\right\| \\
\leq \kappa / 2+\kappa / 8+\varepsilon C_{2} \sup _{s \in\left[t_{a}, \eta \wedge t_{a+1}\right]} \sup _{0\leq i \leq 2}\left|\int_{t_{a}}^{s} \mathcal{Y}_{i}\left(r, u_{r}, \beta_{r}\right) d W_{r}\right| \\
+\varepsilon \sup _{t \in\left[t_{a}, t_{a+1}\right]}\left\|\int_{t_{a}}^{t \wedge \eta}\exp\left(\mathcal{L}_*\left(t-s\right)\right) \widetilde{B}\left(s, u_{s}, \beta_{s}\right) d W_{s}\right\|
\end{array}
$$
since $\varepsilon^{2} C_{1} \leq \kappa /(16 C), C_{3} \kappa^{2} \leq \kappa /(32C)$ and $C_4\xi \kappa \leq\kappa /(32C) $. We again notice that if $\sup _{t \in\left[t_{a}, t_{a+1} \wedge \eta\right]}\left\|w_{t}\right\| \geq \kappa$,
then $\mathcal{B}_{a}$ or $\mathcal{C}_{a}$ must happen.
\endproof

The result now follows from the following lemma

\begin{lemma}
For small enough $\varepsilon$ there exists a constant $C>0$ such that
$$
\begin{aligned}
&\sup _{a \geq 0} \mathbb{P}\left(\mathcal{B}_{a}\right) \leq \exp \left(-C \varepsilon^{-2} \kappa^{2}\right) \\
&\sup _{a \geq 0} \mathbb{P}\left(\mathcal{C}_{a}\right) \leq \exp \left(-C \varepsilon^{-2} \kappa^{2}\right)
\end{aligned}
$$\end{lemma}
\proof
We prove the first inequality taking an exponential moment of the stochastic integral. That is, we start with the estimate
$$
\begin{aligned}
&\mathbb{P}\left(C_{2} \sup _{s \in\left[t_{a}, \eta \wedge t_{a+1}\right]} \sup _{0\leq i \leq 2}\left|\int_{t_{a}}^{s \wedge \eta} \mathcal{Y}_{i}\left(r, u_{r}, \beta_{r}\right) d W_{r}\right| \geq \frac{\kappa}{16 \varepsilon C}\right) =\\
&\mathbb{P}\left(C_{2} \sup _{s \in\left[t_{a}, \eta \wedge t_{a+1}\right]} \sup _{0\leq i \leq 2}\left|\int_{t_{a}}^{s \wedge \eta} \sum_{j=0}^{2} \mathcal{N}_{i j}\left(u_{t}, \beta_t\right)\left\langle B\left(t, u_{t}\right)  d W_{r}, \eta^{j}\right\rangle\right| \geq \frac{\kappa}{16 \varepsilon C}\right) \leq\\
&\mathbb{P}\left(C_{2}M \sup _{s \in\left[t_{a}, \eta \wedge t_{a+1}\right]} \sup _{0\leq j \leq 2}\left|\int_{t_{a}}^{s \wedge \eta} \left\langle B\left(t, u_{t}\right)  d W_{r}, \eta^{j}\right\rangle\right| \geq \frac{\kappa}{16 \varepsilon C}\right)\leq \\
&\sum _{j=0}^2\mathbb{P}\left(C_{2}M \sup _{s \in\left[t_{a}, \eta \wedge t_{a+1}\right]} \left(\int_{t_{a}}^{s \wedge \eta} \left\langle B\left(t, u_{t}\right)  d W_{r}, \eta^{j}\right\rangle\right) \geq \frac{\kappa}{16 \varepsilon C}\right)+\\
&\sum _{j=0}^2\mathbb{P}\left(C_{2}M \sup _{s \in\left[t_{a}, \eta \wedge t_{a+1}\right]} \left(\int_{t_{a}}^{s \wedge \eta} -\left\langle B\left(t, u_{t}\right)  d W_{r}, \eta^{j}\right\rangle\right) \geq \frac{\kappa}{16 \varepsilon C}\right)
\end{aligned}
$$since, as we discussed before the matrixs $\mathcal{N}(t,\beta_t)$ has eigenvalues uniformly bounded in $t$.

We can now define the stochastic process $z^j_t=\int_{t_{a}\wedge \eta}^{s \wedge \eta} \left\langle B\left(t, u_{t}\right)  d W_{r}, \eta^{j}\right\rangle$ that is  a martingale by the Optional Stopping theorem (see \cite{Pro1992}).  Additionally, we define, for a constant $\alpha>0$, the stochastic process $e^j_t=\exp(\alpha z^j_t)=\exp(\alpha \int_{t_{a}\wedge \eta}^{s \wedge \eta} \left\langle B\left(t, u_{t}\right)  d W_{r}, \eta^{j}\right\rangle)$. Being $z_t$ a martingale and using Jensen's inequality (the exponential function is convex) we obtain that $e_t$ is a sub-martingale. Therefore, we have

$$
\begin{aligned}
\mathbb{P}\left( \sup _{s \in\left[t_{a},  t_{a+1}\right]}z^j_t\geq \frac{\kappa}{16C_{2}M \varepsilon C}\right)&\leq \mathbb{P}\left( \sup _{s \in\left[t_{a},  t_{a+1}\right]}e^j_t\geq\exp\left( \frac{\alpha\kappa}{16C_{2}M \varepsilon C}\right)\right)\\ \leq &\mathbb{E}\left[e^j_{t_{a+1}}\right]\exp\left( -\frac{\alpha\kappa}{16C_{2}M \varepsilon C}\right)
\end{aligned}
$$using Doob's submartingale inequality (see page 54 in \cite{revuz1991continuous}). Now we obtain, since $B(t,x)$ is uniformly bounded, we have
$$
\begin{aligned}
&\mathbb{E}\left[e^j_{t_{a+1}}\right]=\mathbb{E}\left[\exp \left(\alpha z^j_{t_{a+1}}\right)\right]=\\
&\mathbb{E}\left[\exp\left(\alpha \int_{t_{a}\wedge \eta}^{t_{a+1} \wedge \eta} \left\langle B\left(t, u_{t}\right)  d W_{r}, \eta^{j}\right\rangle-\frac{\alpha^2}{2}\int_{t_{a}\wedge \eta}^{t_{a+1} \wedge \eta} \left\langle B^{ad}\left(t, u_{t}\right)  \eta^j, B^{ad}\left(t, u_{t}\right)  \eta^j\right\rangle+\right.\right.\\
&\left.\left.\frac{\alpha^2}{2}\int_{t_{a}\wedge \eta}^{t_{a+1} \wedge \eta} \left\langle B^{ad}\left(t, u_{t}\right)  \eta^j, B^{ad}\left(t, u_{t}\right)  \eta^j\right\rangle\right)\right]\leq\\
&\mathbb{E}\left[\exp(\alpha \int_{t_{a}\wedge \eta}^{t_{a+1} \wedge \eta} \left\langle B\left(t, u_{t}\right)  d W_{r}, \eta^{j}\right\rangle)-\frac{\alpha^2}{2}\int_{t_{a}\wedge \eta}^{t_{a+1} \wedge \eta} \left\langle B^{ad}\left(t, u_{t}\right)  \eta^j, B^{ad}\left(t, u_{t}\right)  \eta^j\right\rangle\right]\times\\
&  \quad \quad \quad\quad \quad \quad \quad \quad \quad\exp \left(\frac{\alpha^2}{2}C_{cov}\right)\leq\exp \left(\frac{\alpha^2}{2}C_{cov}\right)
\end{aligned}
$$where we used in the last inequalty that the expectation  is $1$. We can prove this showing that the time derivative of the expectation is zero using Ito formula or, alternatively, showing that the content of the expectation is a martingale, using Girsanov theorem (see \cite{revuz1991continuous}). We now have

$$
\mathbb{P}\left( \sup _{s \in\left[t_{a},  t_{a+1}\right]}z^j_t\geq \frac{\kappa}{16C_{2}M \varepsilon C}\right)\leq\exp \left(\frac{\alpha^2}{2}C_{cov} -\frac{\alpha\kappa}{16C_{2}M \varepsilon C}\right)
$$Choosing now $\alpha=\frac{\kappa}{\varepsilon}\delta$ for a $\delta$ small enough we have
$$
\exp \left(\left(\frac{C_{cov}}{2} \delta^2-\frac{\delta}{16C_{2}M  C}\right)\left(\frac{\kappa}{\varepsilon}\right)^2\right)=\exp \left(-\mathfrak{c}\left(\frac{\kappa}{\varepsilon}\right)^2\right)
$$
for a positive constant $\mathfrak{c}$. The same proof works considering the martingale $-z_t^j$.

Therefore, we obtain $$
\begin{aligned}
&\mathbb{P}\left(C_{2} \sup _{s \in\left[t_{a}, \eta \wedge t_{a+1}\right]} \sup _{0\leq i \leq 2}\left|\int_{t_{a}}^{s \wedge \eta} \mathcal{Y}_{i}\left(r, u_{r}, \beta_{r}\right) d W_{r}\right| \geq \frac{\kappa}{16 \varepsilon C}\right) \leq 6\exp \left(-\mathfrak{c}\left(\frac{\kappa}{\varepsilon}\right)^2\right)
\end{aligned}
$$See also Section 5 in \cite{NeuralFields_MACLAURIN}.

Using the definition of $\widetilde{B}$,
$$
\begin{aligned}
&\mathbb{P}\left(\sup _{t \in\left[t_{a}, t_{a+1}\right]}\left\|\int_{t_{a}}^{t \wedge \eta} \exp\left(\mathcal{L}_*\left(t-s\right)\right) \widetilde{B}\left(s, u_{s}, \beta_{s}\right) d W_{s}\right\| \geq \frac{\kappa}{16 \varepsilon C}\right) \leq \\
&\mathbb{P}\left(\sup _{t \in\left[t_{a}, t_{a+1}\right]}\left\|\int_{t_{a}}^{t \wedge \eta} \exp\left(\mathcal{L}_*\left(t-s\right)\right) B\left(s, u_{s}\right) d W_{s}\right\| \geq \frac{\kappa}{32 \varepsilon C}\right) \\
&\quad+\mathbb{P}\left(\sup _{t \in\left[t_{a}, t_{a+1}\right]}\left\|\sum_{i=0}^{2} \int_{t_{a}}^{t \wedge \eta} \exp\left(\mathcal{L}_*\left(t-s\right)\right) \varphi_{\beta_{s}, i} \mathcal{Y}_{i}\left(s, u_{s}, \beta_{s}\right) d W_{s}\right\| \geq \frac{\kappa}{32 \varepsilon C}\right)
\end{aligned}
$$

The bound of the second term is obtained taking again an exponential bound using exactly the same reasoning as before. From Section 5 in \cite{StochConvBound} we know that there exists a constant $\widetilde{C}$ such that
$$
\mathbb{P}\left(\sup _{t \in\left[t_{a}, t_{a+1}\right]}\left\|\int_{t_{a}}^{t \wedge \eta} \exp\left(\mathcal{L}_*\left(t-s\right)\right) B\left(s, u_{s}\right) d W_{s}\right\| \geq \frac{\kappa}{32 \varepsilon C}\right) \leq \exp \left(-\widetilde{C}\left\{\frac{\kappa}{32 \varepsilon C}\right\}^{2}\right)
$$
For this last bound we must make use of Assumption \ref{integrabilityAssStochConv}.
\endproof
Now the last lemma also concludes the proof of the main result. 
\endproof

\subsection{Variational phases in the stationary reference frame}

Being the family of rotating waves parametrized by the special Euclidean group it is not completely trivial how the variational phase and the approximated variational phase select at every time $t>0$ a particular rotating wave in the original reference frame.  We make this more explicit in this section. Going back to the original stationary reference frame and considering the (approximated) variational phase $\beta_t$ we have
$$
\begin{aligned}
u(\bar{y},t)&=u(R_{-\omega_*t}\bar{x},t)\\
&=\mathcal{T}_{\beta_t}u_*(R_{-\omega_*t}\bar{x})+w(R_{-\omega_*t}\bar{x},t)\\
&=u_*\left(R_{-\beta_t^0}\left(R_{-\omega_*t}\bar{x}-\begin{bmatrix}
\beta_t^1\\
\beta_t^2
\end{bmatrix}\right)\right)+w(R_{-\omega_*t}\bar{x},t)
\end{aligned}
$$
using the splitting presented in the previous sections. First of all, we notice that 
$$
\|w(R_{-\omega_*t}(\cdot),t)\|=\|w(\cdot,t)\|
$$Moreover, we have
$$
u_*\left(R_{-\beta_t^0}\left(R_{-\omega_*t}\bar{x}-\begin{bmatrix}
\beta_t^1\\
\beta_t^2
\end{bmatrix}\right)\right)=u_*\left(R_{-(\beta_t^0+\omega_*t)}\left(\bar{x}-R_{\omega_*t}\begin{bmatrix}
\beta_t^1\\
\beta_t^2
\end{bmatrix}\right)\right)
$$Therefore, the rotating wave selected at every time $t$ by the (approximated) variational phase is
$$
\mathcal{T}_{\tilde{\beta}_t}u_*(\bar{x})=u_*\left(R_{-\tilde{\beta}_t^0}\left(\bar{x}-\begin{bmatrix}
\tilde{\beta}_t^1\\
\tilde{\beta}_t^2
\end{bmatrix}\right)\right)=u_*\left(R_{-(\beta_t^0+\omega_*t)}\left(\bar{x}-R_{\omega_*t}\begin{bmatrix}
\beta_t^1\\
\beta_t^2
\end{bmatrix}\right)\right)
$$and the stochastic process to be considered in the stationary frame is 
$$
\tilde{\beta}_t=\begin{bmatrix}
\tilde{\beta}_t^0\\
\tilde{\beta}_t^1\\
\tilde{\beta}_t^2
\end{bmatrix}
$$
where $$\tilde{\beta}_t^0=(\beta_t^0+\omega_*t)mod(2\pi)$$and 
$$
\begin{bmatrix}
\tilde{\beta}_t^1\\
\tilde{\beta}_t^2
\end{bmatrix}=R_{\omega_*t}\begin{bmatrix}
\beta_t^1\\
\beta_t^2
\end{bmatrix}
$$

\section{Conclusion and Outlook}
\label{sec:conclusion}

We have constructed (approximated) variational phase SODEs for stochastic rotating waves on the whole plane. These SDEs are projections of the stochastic rotating wave to the finite-dimensional manifold of rotations and translations of the deterministic rotating wave. Since the special Euclidean group is non-commutative, we had to introduce a correction term in the linear operator of the SPDE determining the fluctuation of the stochastic rotating wave transverse to the manifold generated by the group symmetries. This makes using variational phases to prove meta-stability results for stochastic rotating waves challenging. Using the (approximated) variational phase, we obtained first meta-stability results for stochastic rotating waves. In particular, we proved that for specific stochastic forcings the probability of the stochastic rotating wave leaving a neighbourhood of the manifold is very small in a large-deviation sense. In future work, we would like to extend this result for larger classes of stochastic forcings.\medskip 

Another interesting open problem is the rigorous understanding of the dynamics of spiral waves. These types of rotating waves present some additional challenges in both the deterministic and the stochastic settings. Already the deterministic nonlinear stability theory of spiral waves is very subtle~\cite{sandstede2020spiral}. Since the asymptotic state at spatial infinity for a spiral wave corresponds to a periodic orbit, we also anticipate serious additional challenges in the stochastic context as already stochastic perturbations of periodic orbits for planar SODEs yield complex phenomena such as cycling~\cite{Day3}.\medskip 

On a much broader scale, our work seems to be the first mathematical approach to intertwine a higher-dimensional symmetry occurring in pattern formation with the study of stochastic reduction and stability. Of course, there are many other patterns beyond rotating/spiral waves, which have complicated symmetry groups. For these cases, it is already challenging to find the best possible construction of a variational phase best adapted to different underlying geometries. In summary, we envision that research uncovering the interaction between symmetry, nonlinearity and stochastic forcing is a promising long-term goal. On the one hand, it bridges several mathematical areas but stability and bifurcation results are also likely to be very useful in various applications.\medskip

\textbf{Acknowledgements:} CK would like to thank the VolkswagenStiftung for support via a Lichtenberg Professorship. CK and GZ would also like to thank Wolf-J\"urgen Beyn for clarifying results about function spaces for deterministic rotating waves.\medskip

\section*{References}

\bibliographystyle{plain}
\bibliography{biblio}

\end{document}